\documentclass{article}

\usepackage{todonotes}
\usepackage{geometry}
\usepackage{amsmath}
\usepackage{amssymb}
\usepackage{amsthm}
\usepackage{bm}
\usepackage{graphicx}
\graphicspath{{figures/}}
\usepackage{color}
\usepackage{tabu}
\usepackage{booktabs}
\usepackage{caption}
\usepackage{authblk}
\usepackage{subfig}
\usepackage{cases}
\usepackage{svg}
\usepackage[section]{placeins}
\usepackage[numbers]{natbib}
\usepackage{hyperref}
\hypersetup{hidelinks}
\usepackage{cleveref}
\usepackage{siunitx}
\usepackage{multirow}

\begin{document}
\title{Capturing Shock Waves by Relaxation Neural Networks}

\author[1]{Nan Zhou}
\author[1, 2, 3, \thanks{Corresponding author: zhengma@sjtu.edu.cn}]{Zheng Ma}

\affil[1]{School of Mathematical Sciences, Shanghai Jiao Tong University, Shanghai,
    200240, China}
\affil[2]{Institute of Natural Sciences, MOE-LSC,
    Shanghai Jiao Tong University, Shanghai, 200240, China}
\affil[3]{CMA-Shanghai, Shanghai Jiao Tong University, Shanghai, 200240, China}

\date{\today}

\maketitle

\begin{abstract}
    In this paper, we put forward a neural network framework to solve the nonlinear hyperbolic systems. This framework, named relaxation neural networks(RelaxNN), is a simple and scalable extension of physics-informed neural networks(PINN). It is shown later that a typical PINN framework struggles to handle shock waves that arise in hyperbolic systems' solutions. This ultimately results in the failure of optimization that is based on gradient descent in the training process. Relaxation systems provide a smooth asymptotic to the discontinuity solution, under the expectation that macroscopic problems can be solved from a microscopic perspective. Based on relaxation systems, the RelaxNN framework alleviates the conflict of losses in the training process of the PINN framework. In addition to the remarkable results demonstrated in numerical simulations, most of the acceleration techniques and improvement strategies aimed at the standard PINN framework can also be applied to the RelaxNN framework.
\end{abstract}

\section{Introduction}

Hyperbolic systems model a series of phenomena that involve wave motion and advective transport of substances. As the preferred mode for communicating information, they arise in the fields of gas dynamics, acoustics, optics, and so on. As for linear hyperbolic systems, the structure of the solution is perspicuous by decomposing any initial value problem into a Riemann problem through characteristic vectors. However, nonlinear hyperbolic systems frequently occur in a broad spectrum of disciplines. The nonlinearity, possibly leading to the formulation of shock waves, brings tremendous difficulties to solving the systems. Finite volume methods(FVM) are oriented to solve the nonlinear hyperbolic systems. Precisely, \citet{godunov1959difference} proposed a general case of reconstruct-evolve-average(REA) algorithm for solving the nonlinear Euler equations. Based on the REA algorithm, a variety of high-resolution methods are designed for capturing shock waves nowadays such as essentially nonoscillatory(ENO) methods raised by \citet{harten1997uniformly}, weighted ENO(WENO) methods raised by \citet{jiang1996efficient}. Besides FVM, finite element methods(FEM) can also apply to nonlinear hyperbolic systems through the discontinuous Galerkin method(DG) proposed by \citet{Shu2000DG}.

Rather than directly being faced with hyperbolic systems, we can approximate the solution of hyperbolic systems smoothly by modifying the hyperbolic systems. \citet{vonneumann1950method} firstly capturing shock waves by introducing an artificial viscosity. However, adding viscous terms is intrinsically non-physical and is computationally expensive. Another way to obtain smooth asymptotics is to interpret the hydrodynamic equations from the microscopic perspective. In spirit to that, \citet{jin2010asymptotic} put forward a series of asymptotic preserving schemes for multiscale kinetic and hyperbolic equations. 

Nowadays, the outstanding performance of deep learning attracts people's attention. The idea of approximating the solution of partial differential equations by neural networks can date back to \citet{dissanayake1994neural}. The theoretical guarantee for a continuous solution is the universal approximation theorem proposed by \citet{hornik1989multilayer}. With the roaring computational power, the unexpected results obtained by Physics-informed neural networks(PINN) \cite{raissi2019physics} polish this field again.  Such a lucid framework, characterized by the combination of data-driven(e.g., initial conditions and boundary conditions) and physical constraint(e.g., the governing equations), spread widely and has been applied to many fields\cite{karniadakis2021physics}. Additionally, by taking advantage of automatic differentiation\cite{baydin2018automatic}, the PINN framework is a non-mesh method that sheds light on breaking the curse of dimensionality.

 However, several possible failure modes in PINN have been fully explored\cite{krishnapriyan2021characterizing,wang2021understanding,daw2023mitigating}. Solving the nonlinear hyperbolic systems is 
a kind of widely reported failure mode due to the discontinuity brought by shock waves \cite{chaumet2022efficient,huang2023limitations}. Naturally, improvement can be raised by the experience of traditional numerical methods. \citet{patel2022thermodynamically} proposed the control volume PINN and introduced the TVD condition into the framework. \citet{de2022weak} raised the weak PINN for capturing shock waves. Above that, \citet{chaumet2022efficient} gave a more efficient weak PINN scheme and entropy admissibility conditions for the uniqueness of the solution. \citet{liuenhancing} came up with weighted equations PINN method. Besides, introducing artificial viscosity directly or implicitly is also adopted to solve this problem \cite{raissi2019physics, michoski2020solving, CHIU2022114909}. Some of them obtained impressive results.

However, these methods make progress at the expense of losing the most attractive properties of deep learning. Introducing numerical derivatives or adopting the prior information about shock location, are making this simple framework mesh-dependent. Also, without the automatic differentiation, our training costs are getting higher and higher.  All we need is to capture the shock automatically without losing the simplicity.  Based on the relaxation systems proposed by \citet{jin1995relaxation}, we put forward the relaxation neural networks(RelaxNN) framework aiming at developing the basic method for solving hydrodynamic equations under the microscopic perspective. Besides the remarkable results we obtained, RelaxNN is a mild and scalable modification of PINN, reserving the simplicity and generality of PINN at the best effort. This means that many sampling strategies and training strategies for PINN could also be incorporated into the RelaxNN framework to acquire more precise results.

Besides the deterministic situations, many uncertainties must be considered due to insufficient knowledge about problems. After all, the closures of some conservation laws are empirical, which inevitably contain some uncertainties. Uncertainty Quantification(UQ) is a field aiming at systematically quantifying the uncertainties that propagate within our models. Along the mindset of UQ, we here consider random uncertainties for high-fidelity simulations \cite{Hu2017}, which is notoriously tough for its high-dimensional properties.  Fortunately, methods based on the PINN framework, known as deep-learning-based methods, are born to solve high-dimensional partial differential equations \cite{yu2018deep, zhang2019quantifying, raissi2019physics,weinan2020machine, zang2020weak, weinan2021algorithms, chen2021deep}.  

The paper is organized as follows. We give a brief introduction about conservation laws and relaxation systems in \autoref{section: CL and shock}, \autoref{section: relaxation systems}. \autoref{section: PINN} introduces the PINN and briefly analyzes the latent failure reason for solving nonlinear hyperbolic systems. In \autoref{section: RelaxNN}, we propose the RelaxNN and give some possible modifications in specific cases. Finally, we implement our RelaxNN framework and compare our results with the solution obtained by \citet{clawpack} in \autoref{section: examples}. We also demonstrate the potential of our method to overcome the curse of dimensionality through the lens of solving the uncertainty quantification problems at the end of the \autoref{section: examples}.

\section{The conservation laws and the shock wave}
\label{section: CL and shock}
From the fluid dynamics perspective, the conservation laws can arise from physics principles. In this article, we consider the systems of conservation laws in one space variable
\begin{equation}\label{eq: conservation laws}
    \partial_{t}\bm{u} + \partial_{x} \bm{F}(\bm{u}) = 0, \, (t,x) \in \mathbb{R}_{+} \times \mathbb{R}, \, \bm{u} \in \mathbb{R}^{n} \, ,
\end{equation}
where $\bm{F}(\bm{u}) \in \mathbb{R}^{n}$ is a vector-valued function, usually called a flux function. As for linear hyperbolic systems, the solution can be easily obtained since it can be viewed as the linear combination of the right vectors at each point in space-time. Physically, it is a superposition of waves propagating at different velocities. When it comes to nonlinear hyperbolic systems, the problem becomes intangible. The shock wave will form at the position where characteristic lines intersect. At those points, the conservation laws ~\eqref{eq: conservation laws} cannot hold in the classic sense\cite{leveque_2002}. Usually, we step back to the fundamental integral conservation laws ~\eqref{eq: integral conservation laws} in this circumstance,
\begin{equation}\label{eq: integral conservation laws}
    \frac{\mathrm{d}}{\mathrm{d} t} \int_{x_{1}}^{x_{2}} \bm{u}(x,t) \, \mathrm{d} x  = \bm{F}(\bm{u}(x_{1}, t)) - \bm{F}(\bm{u}(x_{2}, t)) \, , 
\end{equation}
for any two points $x_{1}$ and $x_{2}$ in $\mathbb{R}$. Furthermore, we pursue the weak solution of conservation law ~\eqref{eq: conservation laws} in this case. $\bm{u}(t,x)$ is the weak solution of conservation law ~\eqref{eq: conservation laws} with initial condition $\bm{u}(0,x) = \bm{u}_{0}(x)$ if
\begin{equation}
    \int_{0}^{\infty} \int_{-\infty}^{\infty} [\bm{u} \bm{\phi}_{t} + \bm{F}(\bm{u})\bm{\phi}_{x}] dxdt = - \int_{0}^{\infty} \bm{u}(0,x)\bm{\phi}(0,x) dx \, ,
\end{equation}
holds for all functions $\bm{\phi}$ $\in$ $C_{0}^{1}$(continuously differentiable with compact support).
In this paper, we mainly consider nonlinear hyperbolic systems.

\section{The relaxation systems and the relaxation schemes}\
\label{section: relaxation systems}
For the unfeasible nonlinear hyperbolic systems, the basic idea of the relaxation systems is to use a local relaxation approximation. So for any system of conservation laws, we can construct a corresponding linear hyperbolic system with a stiff source term that approximates the original system with a small dissipative correction.

Precisely, for the systems of conservation law in one space variable ~\eqref{eq: conservation laws}, \citet{jin1995relaxation} introduced the corresponding relaxation system as
\begin{equation}\label{eq: relaxation system}
\left\{
\begin{aligned}
     \frac{\partial}{\partial t} \bm{u} + \frac{\partial}{\partial x} \bm{v} &= 0 \, , \enspace  \bm{v} \in \mathbb{R}^{n} \, , \\
    \frac{\partial}{\partial t} \bm{v} + A \frac{\partial}{\partial x} \bm{u} &= -\frac{1}{\varepsilon} \left( \bm{v} - \bm{F}(\bm{u}) \right)\, , \enspace \varepsilon > 0 \, ,
\end{aligned}
\right.
\end{equation}
where 
\begin{equation}
    A = \text{Diag}\{ a_{1}, a_{2}, \cdots, a_{n} \} \, , 
\end{equation}
is a chosen positive diagonal matrix. For small $\varepsilon$, applying the Chapman-Enskog expansion\cite{chapman1990mathematical} in ~\eqref{eq: relaxation system}, we can obtain the following approximation for $\bm{u}$ as \cite{liu1987hyperbolic,chen1992hyperbolic}

\begin{equation}\label{eq: approximation for u}
    \frac{\partial}{\partial t} \bm{u} + \frac{\partial}{\partial x}\bm{F}(\bm{u}) = \varepsilon \frac{\partial}{\partial x}\left( \left( A - \bm{F}^{\prime}(\bm{u})^{2} \right) \frac{\partial}{\partial x} \bm{u} \right) , 
\end{equation}
where $\bm{F}^{\prime}(\bm{u})^{2}$ is the Jacobian matrix of the flux function $\bm{F}(\cdot)$. Obviously, \Cref{eq: approximation for u} governs the first-order behavior of the relaxation system ~\eqref{eq: relaxation system}. 

Taking appropriate numerical discretizations to the relaxation system ~\eqref{eq: relaxation system}, we will obtain the relaxing schemes to the original conservation law. The solution of relaxing schemes will accurately approximate the solution of the original equation ~\eqref{eq: conservation laws} when $\varepsilon$ is sufficiently small. The linear structure of the relaxation system benefits us in avoiding the time-consuming Riemann solvers, which is inevitable to the high-resolution methods for nonlinear hyperbolic systems. The last point we need to pay attention to is to use proper implicit time discretizations to overcome the stability constraints due to the stiffness of the systems.    

\section{Physics-informed neural networks and its failure}
\label{section: PINN}
For Cauchy problem of the conservation laws ~\eqref{eq: conservation laws} with the initial condition
\begin{equation}\label{initial condition}
    \bm{u}(0,x) = \bm{u}_{0}(x) \, .
\end{equation}
The PINNs framework for solving this problem is to approximate $\bm{u}(t,x)$ using a deep neural network $\bm{u}^{\text{NN}}_{\theta}(t,x)$ parameterized by trainable variables $\theta$, then train the network by the physics information (e.g., the governing equations) and the given data (e.g., initial conditions and boundary conditions).

The loss function for training the $\bm{u}^{\text{NN}}_{\theta}(t,x)$ is defined as
\begin{subequations}\label{PINNs loss}
\begin{align}
    \mathcal{L}_{\text{PINN}} & = \omega_{\text{PDE}} \mathcal{L}_{\text{PDE}} + \omega_{\text{IC}} \mathcal{L}_{\text{IC}} \, , \\
    \mathcal{L}_{\text{PDE}} & = \frac{1}{\vert \mathcal{T}_{r} \vert}  \sum_{\bm{x} \in \mathcal{T}_{r}} \left( \partial_{t} \bm{u}^{\text{NN}}_{\theta}(\bm{x}_{i}) + \partial_{x} \bm{F}(\bm{u}^{\text{NN}}_{\theta}(\bm{x}_{i})) \right)^{2} \, , \\ 
    \mathcal{L}_{\text{IC}} & = \frac{1}{\vert \mathcal{T}_{ic} \vert} \sum_{\bm{x} \in \mathcal{T}_{ic}} \left( \bm{u}^{\text{NN}}_{\theta}(\bm{x}_{i}) - \bm{u}_{0}(\bm{x}_{i})
    \right)^{2} \, , \\ 
\end{align}
\end{subequations}
where $\mathcal{T}_{r} \subset \Omega$, $\mathcal{T}_{ic} \subset \partial \Omega$ are the sets of "residual points" and $\omega_{\text{PDE}}, \omega_{\text{IC}}$ are the loss weights. 

\begin{figure}[htbp]
\centering
\subfloat[loss curves]
{\label{fig: losses curves burgers sine}
\includegraphics[width = 0.48\textwidth]{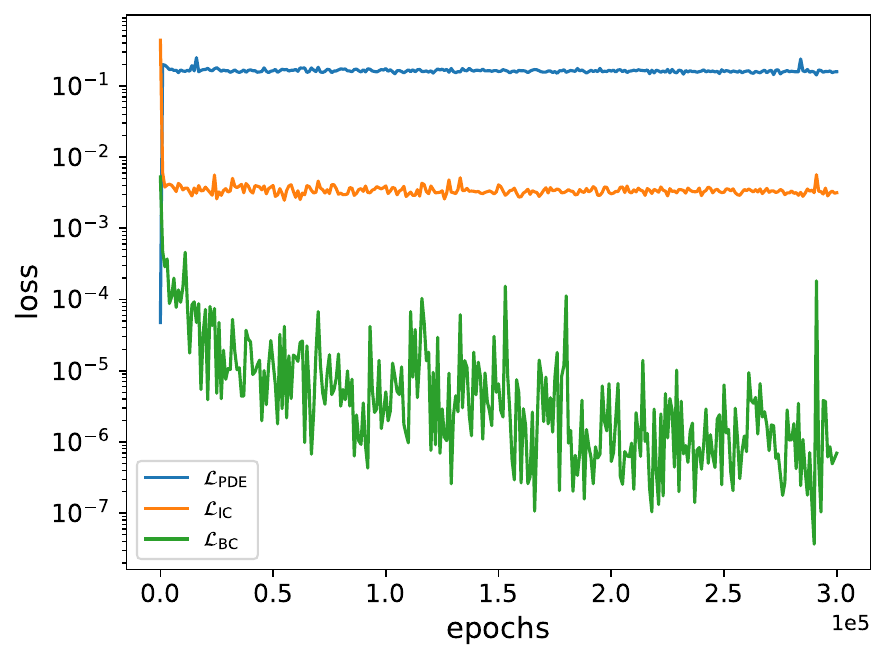}} 
\subfloat[absolute error]
{\label{fig: absolute error burgers sine}
\includegraphics[width = 0.48\textwidth]{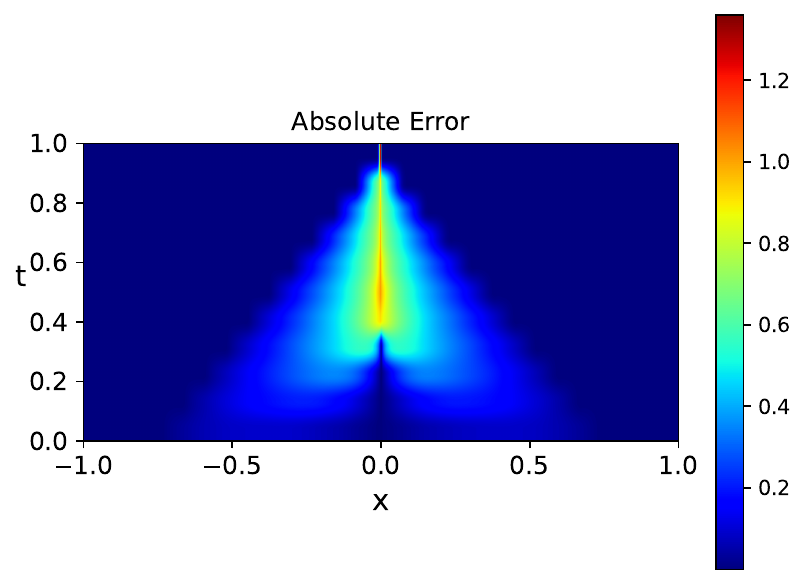}} 
\caption{\textit{Burgers' equation with Sine initial condition}: training a PINN model for approximating the solution to the Burgers' equation with sine initial condition  : (a)loss curves along the training process. (b) absolute error between the prediction of PINN and the reference solution}
\end{figure}

However, minimizing the strong form PDE loss function(Monte-Carlo approximation of $L^{2}$-norm) will not work for learning discontinuous solutions to hyperbolic conservation laws because the conservation laws ~\eqref{eq: conservation laws} near shock hold in the weak sense.  \citet{chaumet2022efficient} presents some analytic computations to explain this failure for inviscid Burgers' equation. In practice, this improper loss function constraint manifests that the PDE loss stagnates at a high level, which finally leads to the failure of gradient descent(see ~\Cref{fig: losses curves burgers sine} for an illuminating example). Also, we can find that the absolute error concentrates on $x=0$ where there is a shock formed(see ~\Cref{fig: absolute error burgers sine}).
An intuitive improvement is to adopt the weak form PDE loss function, \citet{de2022weak} present the weak PINN for nonlinear hyperbolic systems and \citet{chaumet2022efficient} come up with the more efficient weak PINN approximation based on the former work. \citet{liuenhancing} weaken the expression near shock to enhance the shock-capturing ability of PINN. 

\begin{figure}
    \centering
    \includegraphics[width=0.9\textwidth]{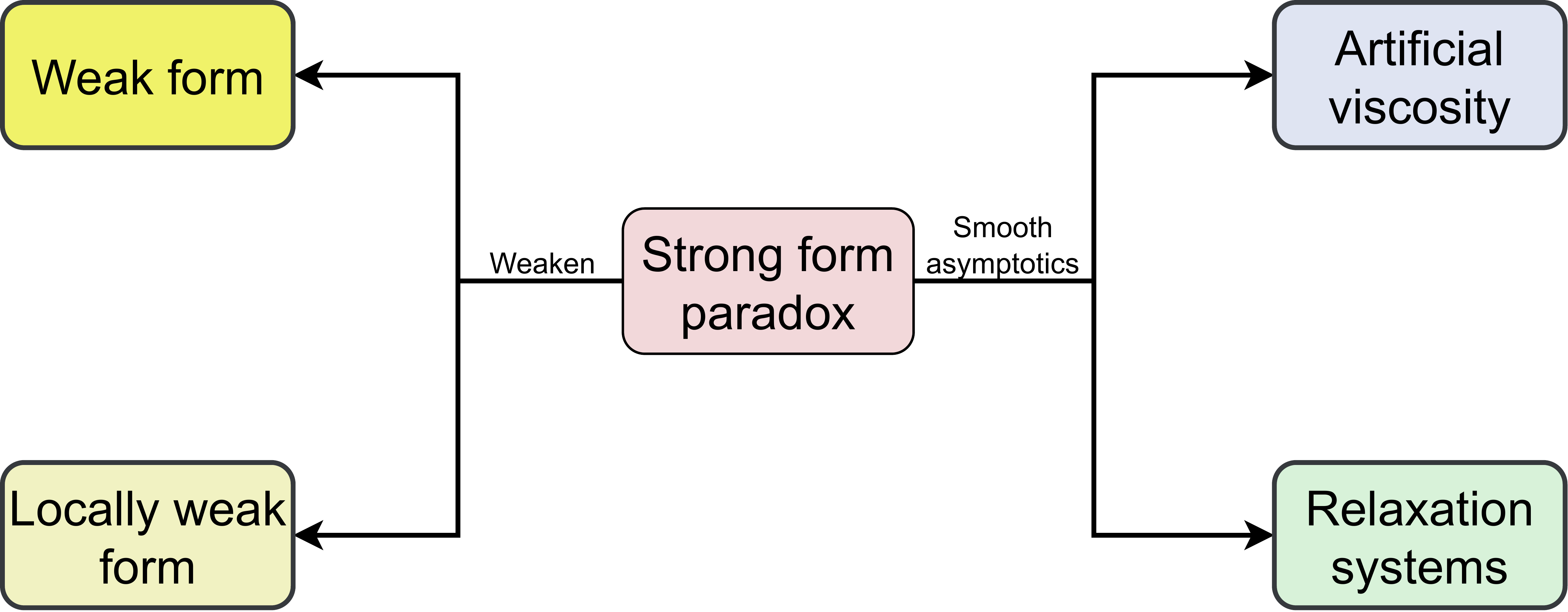}
    \caption{An overview diagram for different improvement ideas}
    \label{fig: overview diagram}
\end{figure}

Another perspective is to obtain smooth asymptotics by modifying the hyperbolic systems(See \Cref{fig: overview diagram} for an overview diagram). With this mindset, adding artificial viscosity to dissipate oscillations is a common strategy. Whether directly solving the modified system with artificial viscosity \cite{raissi2019physics,michoski2020solving} or implicitly, replace the automatical differentiation with the numerical method in PDE loss \cite{CHIU2022114909}. But actually, it is a non-physical adjustment, which is inconsistent with the mindset of physics-informed deep learning. Based on the first principles in physics, the relaxation systems ~\eqref{eq: relaxation system} is derived in spirit from the description of the hydrodynamic equations by the detailed microscopic evolution of gases in kinetic theory in the expectation that the complexity difficulties occurring at the macroscopic level can be illuminated by the hydrodynamic limit process.

\section{Relaxation neural networks}
\label{section: RelaxNN}
For the Cauchy problem of the conservation laws ~\eqref{eq: conservation laws} with the initial condition ~\eqref{initial condition}, the RelaxNN framework comprises two parts according to the relaxation systems ~\eqref{eq: relaxation system}, one approximates $\bm{u}(t,x)$ by a deep neural network $\bm{u}^{\text{NN}}_{\bm{\theta}_{1}}(t,x)$ and another approximates $\bm{v}$ by another deep neural network $\bm{v}^{\text{NN}}_{\bm{\theta}_{2}}(t,x)$, separately parameterized by trainable variables $\bm{\theta}_{1}$ and $\bm{\theta}_{2}$. If we solve the relaxation systems ~\eqref{eq: relaxation system} with PINNs' framework, the PDE loss will be 
\begin{subequations}
\begin{align}
     \mathcal{L}_{\text{PDE}} & = \mathcal{L}_{\text{residual}} + \mathcal{L}_{\text{dissipate}} \, ,  \\
     \mathcal{L}_{\text{residual}} & = \frac{1}{\vert \mathcal{T}_{r} \vert}  \sum_{\bm{x} \in \mathcal{T}_{r}} \left( \partial_{t} \bm{u}^{\text{NN}}_{\bm{\theta}_{1}}(\bm{x}_{i}) + \partial_{x} \bm{v}^{\text{NN}}_{\bm{\theta}_{2}}(\bm{x}_{i}) \right)^{2} \, ,\\
     \mathcal{L}_{\text{dissipate}} & = \frac{1}{\vert \mathcal{T}_{r} \vert}  \sum_{\bm{x} \in \mathcal{T}_{r}} \left( \varepsilon \left[\partial_{t}\bm{v}^{\text{NN}}_{\bm{\theta}_{2}}(\bm{x}_{i}) + A \partial_{x} \bm{u}^{\text{NN}}_{\bm{\theta}_{1}}(\bm{x}_{i})\right] + \bm{v}^{\text{NN}}_{\bm{\theta}_{2}}(\bm{x}_{i}) - \bm{F}(\bm{u}^{\text{NN}}_{\bm{\theta}_{1}}(\bm{x}_{i})) \right)^{2} \, .
\end{align}
\end{subequations}

In the small relaxation limit ($\varepsilon \rightarrow 0^{+}$), 
\begin{equation}
    \mathcal{L}_{\text{dissipate}} \approx \mathcal{L}_{\text{flux}} \triangleq \frac{1}{\vert \mathcal{T}_{r} \vert}  \sum_{\bm{x} \in \mathcal{T}_{r}} \left( \bm{v}^{\text{NN}}_{\bm{\theta}_{2}}(\bm{x}_{i}) - \bm{F}(\bm{u}^{\text{NN}}_{\bm{\theta}_{1}}(\bm{x}_{i}) ) \right)^{2} \, .
\end{equation}

\begin{figure}
    \centering
    \includegraphics[width = 0.96\textwidth]{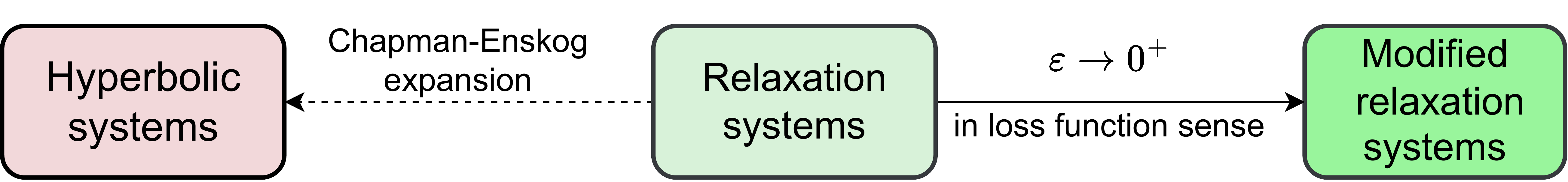}
    \caption{Brief illustration about the relationship of different systems mentioned above. Here the dashed line represents the low-order approximation.}
    \label{fig: relationship}
\end{figure}

In fact, the RelaxNN framework for the original conservation laws ~\eqref{eq: conservation laws} is the PINN framework to the following modified relaxation systems(See \Cref{fig: relationship} for the relationship of different systems)
\begin{equation}\label{eq: modified relaxation system}
\left\{
\begin{aligned}
     \frac{\partial}{\partial t} \bm{u} + \frac{\partial}{\partial x} \bm{v} &= 0 \, , \enspace  \bm{u},\bm{v} \in \mathbb{R}^{n} \, , \\ 
     \bm{v} - \bm{F}(\bm{u}) & = 0 \, .
\end{aligned}
\right.
\end{equation}

In traditional numerical methods, for the relaxation system ~\eqref{eq: relaxation system}, we must choose the initial condition as
\begin{subequations}
\begin{align}
    \bm{u}(0,x)  & = \bm{u}_{0}(0,x) ,  \\
    \bm{v}(0,x)  & = \bm{v}_{0}(0,x) \triangleq \bm
{F}(\bm{u}_{0}(0,x)) . \label{flux initial condition}
\end{align}
\end{subequations}

\begin{figure}
    \centering
    \includegraphics[width = 0.9\textwidth]{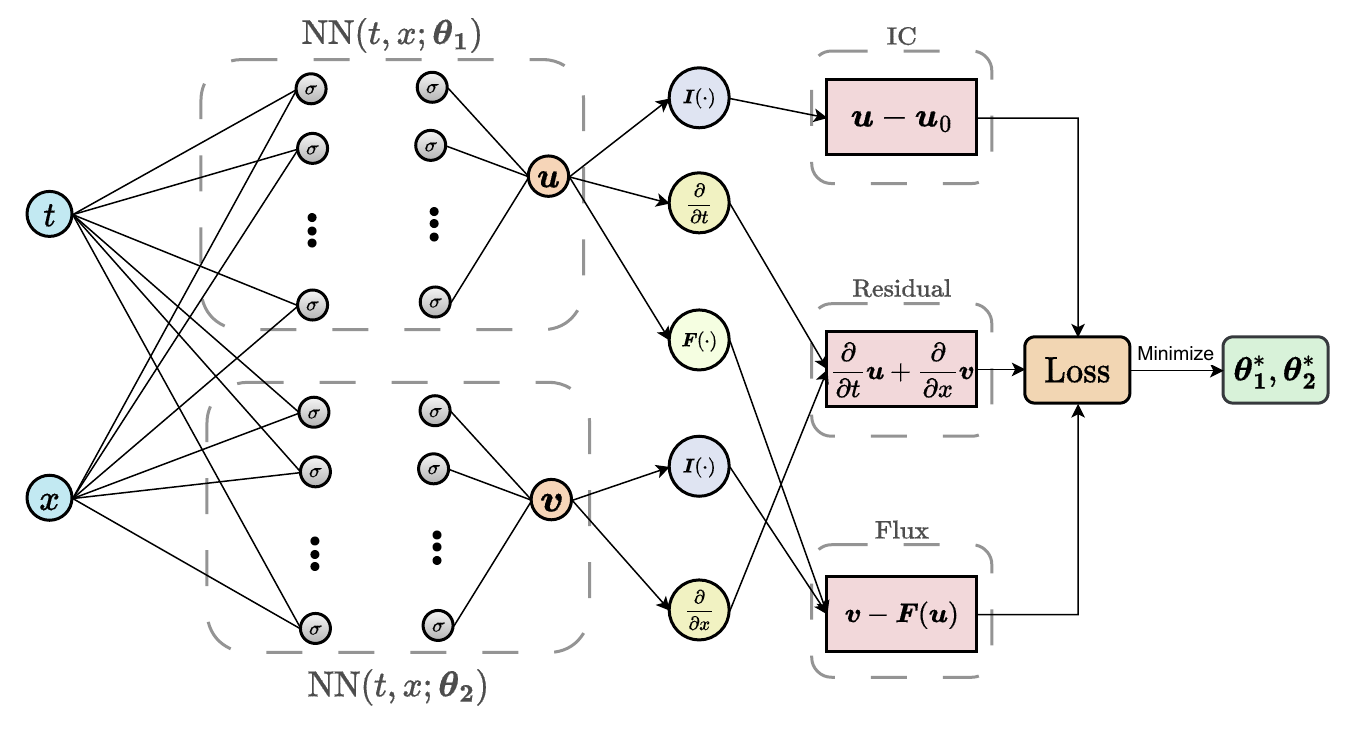}
    \caption{The schematic diagram of the RelaxNN to solve the hyperbolic systems with initial condition (IC) $\bm{u}_{0}$ and flux function $\bm{F}(\cdot)$.}
    \label{fig: schematic of RelaxNN}
\end{figure}

Choosing the special initial condition for $\bm{v}$ is to avoid the introduction of an initial layer through the relaxation system because of its stiffness. Now we do not need to impose the special initial condition penalty for the neural network $\bm{v}^{\text{NN}}_{\bm{\theta}_{2}}(t,x)$ since we replace the $\mathcal{L}_{\text{dissipate}}$ with $\mathcal{L}_{\text{flux}}$ during the optimization process. So the training loss for the relaxation neural networks will be defined as
\begin{equation}
    \mathcal{L}_{\text{RelaxNN}} = \omega_{\text{residual}} \mathcal{L}_{\text{residual}} + \omega_{\text{flux}}\mathcal{L}_{\text{flux}} + \omega_{\text{IC}}\mathcal{L}_{\text{IC}} \, .
\end{equation}
See the visualization of RelaxNN in \Cref{fig: schematic of RelaxNN}.

Sometimes, we do not need to relax the whole hyperbolic system. The relaxation system is still linear if we only relax some of the conservation laws(e.g., the conservation law of momentum or energy). Additionally, the numerical methods for the original relaxation system are not always well-balanced and thus introduce spurious waves\cite{liu2020moving}. Similar to Suliciu's relaxation and its applications \cite{suliciu1990modelling, suliciu1998thermodynamics, coquel1998relaxation}, \citet{liu2020moving} introduce the partial relaxation system which only relaxes the second equation of the Saint-Venant system of shallow water equations. In this article, we adopt and generalize the partial relaxation systems. Later we will make a more precise illustration of the specified equations we are concerned with.

\subsection{Inviscid Burgers' equation}
We first consider the inviscid Burgers' equation, the simplest hyperbolic equation with shock wave
\begin{equation}
    \partial_{t} u + \partial_{x} \left(\frac{1}{2} u^{2} \right) = 0. \label{eq: burgers equation}
\end{equation}
The modified relaxation systems for Burgers' equation are
\begin{equation}\label{eq: modified relaxation system for burgers}
\left\{
\begin{aligned}
     \frac{\partial}{\partial t} u + \frac{\partial}{\partial x} v &= 0 \, , \\
     v - \frac{1}{2}u^{2} & = 0 \, .
\end{aligned}
\right.
\end{equation}
For brevity, denote $\bm{u}_{\bm{\theta_{1}}}^{\text{NN}} = u_{\bm{\theta_{1}}}^{\text{NN}}$ approximates $u$ and $\bm{v}_{\bm{\theta_{2}}}^{\text{NN}} = v_{\bm{\theta_{2}}}^{\text{NN}}$ approximates $v$, the residual loss term and flux loss term of RelaxNN framework can be presented more precisely as
\begin{subequations}\label{burgers loss function}
\begin{align}
     \mathcal{L}_{\text{residual}} & = \frac{1}{\vert \mathcal{T}_{r} \vert}  \sum_{\bm{x} \in \mathcal{T}_{r}} \left( \partial_{t} u^{\text{NN}}_{\bm{\theta_{1}}}(\bm{x}_{i}) + \partial_{x} v^{\text{NN}}_{\bm{\theta_{2}}}(\bm{x}_{i}) \right)^{2} \, ,\\
     \mathcal{L}_{\text{flux}} & = \frac{1}{\vert \mathcal{T}_{r} \vert}  \sum_{\bm{x} \in \mathcal{T}_{r}} \left( v^{\text{NN}}_{\bm{\theta_{2}}}(\bm{x}_{i}) - \frac{1}{2}(u^{\text{NN}}_{\bm{\theta_{1}}}(\bm{x}_{i}))^{2} \right)^{2} \, .
\end{align}
\end{subequations}
\subsection{Shallow water equations}
The shallow water equations describe the evolution of the incompressible fluid(constant density) with small-amplitude waves on the surface of the sphere
\begin{equation}
\left\{
\begin{aligned}
     \partial_{t} h & + \partial_{x} \left( h u \right)  = 0 \, ,\\ 
     \partial_{t} \left( hu \right) & + \partial_{x} \left( hu^{2} + \frac{1}{2}gh^{2} \right) & = 0 \, ,
\end{aligned}
\right.
\label{eq: shallow water equations}
\end{equation}
here $h$ is the depth, $u$ is the velocity and $g$ is the gravitational constant. 
\subsubsection{Type1 fully modified relaxation systems}
Here we can relax the whole hyperbolic system with the type1 fully modified relaxation systems as 
\begin{equation}\label{eq: type1 fully modified relaxation system for swe}
\left\{
\begin{aligned}
     \frac{\partial}{\partial t} h + \frac{\partial}{\partial x} v &= 0 \, , \\
     \frac{\partial}{\partial t} (h u) + \frac{\partial}{\partial x} \varphi &= 0 \, , \\
     v - hu & = 0 \, , \\
     \varphi - (hu^{2} + \frac{1}{2}gh^{2}) & = 0 \, . 
\end{aligned}
\right.
\end{equation}
 Denote $\bm{u}_{\bm{\theta_{1}}}^{\text{NN}} =\begin{pmatrix}
    h_{\bm{\theta_{1}}}^{\text{NN}}, u_{\bm{\theta_{1}}}^{\text{NN}}
\end{pmatrix}^{\text{T}} $ approximates $\begin{pmatrix}
    h, u
\end{pmatrix}^{\text{T}}$ and $\bm{v}_{\bm{\theta_{2}}}^{\text{NN}} = \begin{pmatrix}
    v_{\bm{\theta_{2}}}^{\text{NN}}, \varphi_{\bm{\theta_{2}}}^{\text{NN}}
\end{pmatrix}^{\text{T}}$ approximates $\begin{pmatrix}
    v, \varphi
\end{pmatrix}^{\text{T}}$, the loss function of RelaxNN framework can be presented as
\begin{equation}
    \mathcal{L}_{\text{RelaxNN}}^{\text{type1}} = \omega_{\text{residual}} \mathcal{L}_{\text{residual}}^{\text{type1}} + \omega_{\text{flux}}\mathcal{L}_{\text{flux}}^{\text{type1}} + \omega_{\text{IC}}\mathcal{L}_{\text{IC}} \, ,
\end{equation}
where the residual loss term can be presented as
\begin{equation}
\begin{split}
    \mathcal{L}_{\text{residual}}^{\text{type1}} & = \frac{\omega_{r}^{m}}{\vert \mathcal{T}_{r} \vert}  \sum_{\bm{x} \in \mathcal{T}_{r}} \left( \partial_{t} h^{\text{NN}}_{\bm{\theta_{1}}}(\bm{x}_{i}) + \partial_{x} v^{\text{NN}}_{\bm{\theta_{2}}}(\bm{x}_{i}) \right)^{2} \\
    & + \frac{\omega_{r}^{p}}{\vert \mathcal{T}_{r} \vert}  \sum_{\bm{x} \in \mathcal{T}_{r}} \left( \partial_{t} (h^{\text{NN}}_{\bm{\theta_{1}}}u^{\text{NN}}_{\bm{\theta_{1}}})(\bm{x}_{i}) + \partial_{x} \varphi^{\text{NN}}_{\bm{\theta_{2}}}(\bm{x}_{i}) \right)^{2} \, ,    
\end{split}
\end{equation}
and the flux loss term can be presented as
\begin{equation}
\begin{split}
     \mathcal{L}_{\text{flux}}^{\text{type1}} & = \frac{\omega_{f}^{m}}{\vert \mathcal{T}_{r} \vert}  \sum_{\bm{x} \in \mathcal{T}_{r}} \left( v^{\text{NN}}_{\bm{\theta_{2}}}(\bm{x}_{i}) - (h^{\text{NN}}_{\bm{\theta_{1}}} u^{\text{NN}}_{\bm{\theta_{1}}})(\bm{x}_{i})\right)^{2} \, \\
     & + \frac{\omega_{f}^{p}}{\vert \mathcal{T}_{r} \vert}  \sum_{\bm{x} \in \mathcal{T}_{r}} \left( \varphi^{\text{NN}}_{\bm{\theta_{2}}}(\bm{x}_{i}) - \left(h^{\text{NN}}_{\bm{\theta_{1}}} (u^{\text{NN}}_{\bm{\theta_{1}}})^{2} + \frac{1}{2} g (h^{\text{NN}}_{\bm{\theta_{1}}})^{2} \right)(\bm{x}_{i}) \right)^{2} \, .    
\end{split}
\end{equation}
\subsubsection{Type2 partially modified relaxation systems}
However, we can also only relax the conservation law of momentum with the type2 partial relaxation systems as
\begin{equation}\label{eq: type2 partially modified relaxation system for swe}
\left\{
\begin{aligned}
     \frac{\partial}{\partial t} h + \frac{\partial}{\partial x} (h u) &= 0 \, , \\
     \frac{\partial}{\partial t} (h u) + \frac{\partial}{\partial x} \varphi &= 0 \, , \\
     \varphi - (hu^{2} + \frac{1}{2}gh^{2}) & = 0 \, . 
\end{aligned}
\right.
\end{equation}
Denote $\bm{u}_{\bm{\theta_{1}}}^{\text{NN}} =\begin{pmatrix}
    h_{\bm{\theta_{1}}}^{\text{NN}}, u_{\bm{\theta_{1}}}^{\text{NN}}
\end{pmatrix}^{\text{T}} $ approximates $\begin{pmatrix}
    h, u
\end{pmatrix}^{\text{T}}$ and $\bm{v}_{\bm{\theta_{2}}}^{\text{NN}} =\varphi_{\bm{\theta_{2}}}^{\text{NN}}$ approximates $\varphi$, the loss function of RelaxNN framework can be presented as
\begin{equation}
    \mathcal{L}_{\text{RelaxNN}}^{\text{type2}} = \omega_{\text{residual}} \mathcal{L}_{\text{residual}}^{\text{type2}} + \omega_{\text{flux}}\mathcal{L}_{\text{flux}}^{\text{type2}} + \omega_{\text{IC}}\mathcal{L}_{\text{IC}} \, ,
\end{equation}
where the residual loss term can be presented as
\begin{equation}
\begin{split}
    \mathcal{L}_{\text{residual}}^{\text{type2}} & = \frac{\omega_{r}^{m}}{\vert \mathcal{T}_{r} \vert}  \sum_{\bm{x} \in \mathcal{T}_{r}} \left( \partial_{t} h^{\text{NN}}_{\bm{\theta_{1}}}(\bm{x}_{i}) + \partial_{x} (h^{\text{NN}}_{\bm{\theta_{1}}} u^{\text{NN}}_{\bm{\theta_{1}}})(\bm{x}_{i}) \right)^{2} \\
    & + \frac{\omega_{r}^{p}}{\vert \mathcal{T}_{r} \vert}  \sum_{\bm{x} \in \mathcal{T}_{r}} \left( \partial_{t} (h^{\text{NN}}_{\bm{\theta_{1}}}u^{\text{NN}}_{\bm{\theta_{1}}})(\bm{x}_{i}) + \partial_{x} \varphi^{\text{NN}}_{\bm{\theta_{2}}}(\bm{x}_{i}) \right)^{2} \, ,    
\end{split}
\end{equation}
and the flux loss term can be presented as
\begin{equation}
\begin{split}
     \mathcal{L}_{\text{flux}}^{\text{type2}} & =  \frac{\omega_{f}^{p}}{\vert \mathcal{T}_{r} \vert}  \sum_{\bm{x} \in \mathcal{T}_{r}} \left( \varphi^{\text{NN}}_{\bm{\theta_{2}}}(\bm{x}_{i}) - \left(h^{\text{NN}}_{\bm{\theta_{1}}} (u^{\text{NN}}_{\bm{\theta_{1}}})^{2} + \frac{1}{2} g (h^{\text{NN}}_{\bm{\theta_{1}}})^{2} \right)(\bm{x}_{i}) \right)^{2} \, .    
\end{split}
\end{equation}
\subsection{Euler equations}
The Euler equations are the conservation laws governing the compressible, adiabatic, and inviscid flow in fluid dynamics. Here, we consider $1-\rm{D}$ Euler equations of gas dynamics in conservative forms
\begin{equation}
\left\{
\begin{aligned}
     \partial_{t} \rho & + \partial_{x} \left( \rho u \right)  = 0 \, ,\\ 
     \partial_{t} \left( \rho u \right) & + \partial_{x} \left( \rho u^{2} + p \right) = 0 \, , \\
     \partial_{t} \mathrm{E} & + \partial_{x} \left( u (\mathrm{E} + p ) \right) = 0  \, ,
\end{aligned}
\label{eq: Euler equations}
\right.
\end{equation}
here $\rho$ is the density, $u$ is the velocity and $p$ is the pressure for an ideal gas. To close the system of equations, we introduce the equation of state for an ideal polytropic gas
\begin{equation}
    \mathrm{E} = \frac{p}{\gamma - 1} + \frac{1}{2} \rho u^{2} \, ,
\end{equation}
with $\gamma = 1.4$ for a diatomic gas.

\subsubsection{Type1 fully modified relaxation systems}
Here we can relax the whole hyperbolic systems with the type1 fully modified relaxation systems as 
\begin{equation}\label{eq: type1 fully modified relaxation system for euler}
\left\{
\begin{aligned}
     \frac{\partial}{\partial t} \rho + \frac{\partial}{\partial x} v &= 0 \, , \\
     \frac{\partial}{\partial t} (\rho u) + \frac{\partial}{\partial x} \varphi &= 0 \, , \\
     \frac{\partial}{\partial t} (\mathrm{E}) + \frac{\partial}{\partial x} \zeta &= 0 \, , \\     
     v - (\rho u) & = 0 \, , \\
     \varphi - (\rho u^{2} + p) & = 0 \, , \\
     \zeta - \left( u(\mathrm{E}+p) \right) & = 0 \, .
\end{aligned}
\right.
\end{equation}
Denote $\bm{u}_{\bm{\theta_{1}}}^{\text{NN}} =\begin{pmatrix}
    \rho_{\bm{\theta_{1}}}^{\text{NN}}, u_{\bm{\theta_{1}}}^{\text{NN}},
    p_{\bm{\theta_{1}}}^{\text{NN}}
\end{pmatrix}^{\text{T}} $ approximates $\begin{pmatrix}
    \rho, u, p
\end{pmatrix}^{\text{T}}$ and $\bm{v}_{\bm{\theta_{2}}}^{\text{NN}} = \begin{pmatrix}
    v_{\bm{\theta_{2}}}^{\text{NN}}, \varphi_{\bm{\theta_{2}}}^{\text{NN}},
    \zeta_{\bm{\theta_{2}}}^{\text{NN}}
\end{pmatrix}^{\text{T}}$ approximates $\begin{pmatrix}
    v, \varphi, \zeta
\end{pmatrix}^{\text{T}}$, the loss function of RelaxNN framework can be presented as
\begin{equation}
    \mathcal{L}_{\text{RelaxNN}}^{\text{type1}} = \omega_{\text{residual}} \mathcal{L}_{\text{residual}}^{\text{type1}} + \omega_{\text{flux}}\mathcal{L}_{\text{flux}}^{\text{type1}} + \omega_{\text{IC}}\mathcal{L}_{\text{IC}} \, ,
\end{equation}
where the residual loss term can be presented as
\begin{equation}
\begin{split}
    \mathcal{L}_{\text{residual}}^{\text{type1}} & = \frac{\omega_{r}^{m}}{\vert \mathcal{T}_{r} \vert}  \sum_{\bm{x} \in \mathcal{T}_{r}} \left( \partial_{t} \rho^{\text{NN}}_{\bm{\theta_{1}}}(\bm{x}_{i}) + \partial_{x} v^{\text{NN}}_{\bm{\theta_{2}}}(\bm{x}_{i}) \right)^{2} \\
    & + \frac{\omega_{r}^{p}}{\vert \mathcal{T}_{r} \vert}  \sum_{\bm{x} \in \mathcal{T}_{r}} \left( \partial_{t} (\rho^{\text{NN}}_{\bm{\theta_{1}}}u^{\text{NN}}_{\bm{\theta_{1}}})(\bm{x}_{i}) + \partial_{x} \varphi^{\text{NN}}_{\bm{\theta_{2}}}(\bm{x}_{i}) \right)^{2} \\
    & + \frac{\omega_{r}^{e}}{\vert \mathcal{T}_{r} \vert}  \sum_{\bm{x} \in \mathcal{T}_{r}} \left( \partial_{t} \left(\frac{1}{\gamma - 1}\rho^{\text{NN}}_{\bm{\theta_{1}}} + \frac{1}{2}\rho^{\text{NN}}_{\bm{\theta_{1}}}(u^{\text{NN}}_{\bm{\theta_{1}}})^{2}
    \right)(\bm{x}_{i}) + \partial_{x} \zeta^{\text{NN}}_{\bm{\theta_{2}}}(\bm{x}_{i}) \right)^{2} \, ,
\end{split}
\end{equation}
and the flux loss term can be presented as
\begin{equation}
\begin{split}
     \mathcal{L}_{\text{flux}}^{\text{type1}} & = \frac{\omega_{f}^{m}}{\vert \mathcal{T}_{r} \vert}  \sum_{\bm{x} \in \mathcal{T}_{r}} \left( v^{\text{NN}}_{\bm{\theta_{2}}}(\bm{x}_{i}) - (\rho^{\text{NN}}_{\bm{\theta_{1}}} u^{\text{NN}}_{\bm{\theta_{1}}})(\bm{x}_{i})\right)^{2} \, \\
     & + \frac{\omega_{f}^{p}}{\vert \mathcal{T}_{r} \vert}  \sum_{\bm{x} \in \mathcal{T}_{r}} \left( \varphi^{\text{NN}}_{\bm{\theta_{2}}}(\bm{x}_{i}) - \left(\rho^{\text{NN}}_{\bm{\theta_{1}}} (u^{\text{NN}}_{\bm{\theta_{1}}})^{2} + ^{\text{NN}}_{\bm{\theta_{1}}} \right) (\bm{x}_{i}) \right)^{2} \\
     & + \frac{\omega_{f}^{e}}{\vert \mathcal{T}_{r} \vert}  \sum_{\bm{x} \in \mathcal{T}_{r}} \left( \zeta^{\text{NN}}_{\bm{\theta_{2}}}(\bm{x}_{i}) - \left(
     \frac{\gamma}{\gamma - 1} \rho^{\text{NN}}_{\bm{\theta_{1}}}u^{\text{NN}}_{\bm{\theta_{1}}}
     + \frac{1}{2}\rho^{\text{NN}}_{\bm{\theta_{1}}}(u^{\text{NN}}_{\bm{\theta_{1}}})^{3} \right) (\bm{x}_{i}) \right)^{2} \, .
\end{split}
\end{equation}
\subsubsection{Type2 partially modified relaxation systems}
Also, we can relax the conservation laws of momentum and energy with the type2 partially modified relaxation systems as 
\begin{equation}\label{eq: type2 partially modified relaxation system for euler}
\left\{
\begin{aligned}
     \frac{\partial}{\partial t} \rho + \frac{\partial}{\partial x} (\rho u) &= 0 \, , \\
     \frac{\partial}{\partial t} (\rho u) + \frac{\partial}{\partial x} \varphi &= 0 \, , \\
     \frac{\partial}{\partial t} (\mathrm{E}) + \frac{\partial}{\partial x} \zeta &= 0 \, , \\     
     \varphi - (\rho u^{2} + p) & = 0 \, , \\
     \zeta - \left( u(\mathrm{E}+p) \right) & = 0 \, .
\end{aligned}
\right.
\end{equation}
Denote $\bm{u}_{\bm{\theta_{1}}}^{\text{NN}} =\begin{pmatrix}
    \rho_{\bm{\theta_{1}}}^{\text{NN}}, u_{\bm{\theta_{1}}}^{\text{NN}},
    p_{\bm{\theta_{1}}}^{\text{NN}}
\end{pmatrix}^{\text{T}} $ approximates $\begin{pmatrix}
    \rho, u, p
\end{pmatrix}^{\text{T}}$ and $\bm{v}_{\bm{\theta_{2}}}^{\text{NN}} = \begin{pmatrix}
    \varphi_{\bm{\theta_{2}}}^{\text{NN}},
    \zeta_{\bm{\theta_{2}}}^{\text{NN}}
\end{pmatrix}^{\text{T}}$ approximates $\begin{pmatrix}
    \varphi, \zeta
\end{pmatrix}^{\text{T}}$, the loss function of RelaxNN framework can be presented as
\begin{equation}
    \mathcal{L}_{\text{RelaxNN}}^{\text{type2}} = \omega_{\text{residual}} \mathcal{L}_{\text{residual}}^{\text{type2}} + \omega_{\text{flux}}\mathcal{L}_{\text{flux}}^{\text{type2}} + \omega_{\text{IC}}\mathcal{L}_{\text{IC}} \, ,
\end{equation}
where the residual loss term can be presented as
\begin{equation}
\begin{split}
    \mathcal{L}_{\text{residual}}^{\text{type2}} & = \frac{\omega_{r}^{m}}{\vert \mathcal{T}_{r} \vert}  \sum_{\bm{x} \in \mathcal{T}_{r}} \left( \partial_{t} \rho^{\text{NN}}_{\bm{\theta_{1}}}(\bm{x}_{i}) + \partial_{x} (\rho^{\text{NN}}_{\bm{\theta_{1}}} u^{\text{NN}}_{\bm{\theta_{1}}})(\bm{x}_{i}) \right)^{2} \\
    & + \frac{\omega_{r}^{p}}{\vert \mathcal{T}_{r} \vert}  \sum_{\bm{x} \in \mathcal{T}_{r}} \left( \partial_{t} (\rho^{\text{NN}}_{\bm{\theta_{1}}}u^{\text{NN}}_{\bm{\theta_{1}}})(\bm{x}_{i}) + \partial_{x} \left( \rho^{\text{NN}}_{\bm{\theta_{1}}}(u^{\text{NN}}_{\bm{\theta_{1}}})^{2} + p^{\text{NN}}_{\bm{\theta_{1}}} \right)(\bm{x}_{i}) \right)^{2} \\
    & + \frac{\omega_{r}^{e}}{\vert \mathcal{T}_{r} \vert}  \sum_{\bm{x} \in \mathcal{T}_{r}} \left( \partial_{t} \left(\frac{1}{\gamma - 1}\rho^{\text{NN}}_{\bm{\theta_{1}}} + \frac{1}{2}\rho^{\text{NN}}_{\bm{\theta_{1}}}(u^{\text{NN}}_{\bm{\theta_{1}}})^{2}
    \right)(\bm{x}_{i}) + \partial_{x} \zeta^{\text{NN}}_{\bm{\theta_{2}}}(\bm{x}_{i}) \right)^{2} \, ,
\end{split}
\end{equation}
and the flux loss term can be presented as
\begin{equation}
\begin{split}
     \mathcal{L}_{\text{flux}}^{\text{type2}} & =  \frac{\omega_{f}^{p}}{\vert \mathcal{T}_{r} \vert}  \sum_{\bm{x} \in \mathcal{T}_{r}} \left( \varphi^{\text{NN}}_{\bm{\theta_{2}}}(\bm{x}_{i}) - \left(\rho^{\text{NN}}_{\bm{\theta_{1}}} (u^{\text{NN}}_{\bm{\theta_{1}}})^{2} +p^{\text{NN}}_{\bm{\theta_{1}}} \right) (\bm{x}_{i}) \right)^{2} \\
     & + \frac{\omega_{f}^{e}}{\vert \mathcal{T}_{r} \vert}  \sum_{\bm{x} \in \mathcal{T}_{r}} \left( \zeta^{\text{NN}}_{\bm{\theta_{2}}}(\bm{x}_{i}) - \left(
     \frac{\gamma}{\gamma - 1} \rho^{\text{NN}}_{\bm{\theta_{1}}}u^{\text{NN}}_{\bm{\theta_{1}}}
     + \frac{1}{2}\rho^{\text{NN}}_{\bm{\theta_{1}}}(u^{\text{NN}}_{\bm{\theta_{1}}})^{3} \right) (\bm{x}_{i}) \right)^{2} \, .
\end{split}
\end{equation}
\subsubsection{Type3 partially modified relaxation systems}
Of course, we can only relax the conservation law of energy with the type3 partially modified relaxation systems as 
\begin{equation}\label{eq: type3 partially modified relaxation system for euler}
\left\{
\begin{aligned}
     \frac{\partial}{\partial t} \rho + \frac{\partial}{\partial x} (\rho u) &= 0 \, , \\
     \frac{\partial}{\partial t} (\rho u) + \frac{\partial}{\partial x} (\rho u^{2} + p) &= 0 \, , \\
     \frac{\partial}{\partial t} (\mathrm{E}) + \frac{\partial}{\partial x} \zeta &= 0 \, , \\ 
     \zeta - \left( u(\mathrm{E}+p) \right) & = 0 \, .
\end{aligned}
\right.
\end{equation}
Denote $\bm{u}_{\bm{\theta_{1}}}^{\text{NN}} =\begin{pmatrix}
    \rho_{\bm{\theta_{1}}}^{\text{NN}}, u_{\bm{\theta_{1}}}^{\text{NN}},
    p_{\bm{\theta_{1}}}^{\text{NN}}
\end{pmatrix}^{\text{T}} $ approximates $\begin{pmatrix}
    \rho, u, p
\end{pmatrix}^{\text{T}}$ and $\bm{v}_{\bm{\theta_{2}}}^{\text{NN}} = \zeta_{\bm{\theta_{2}}}^{\text{NN}}$ approximates $\zeta$, the loss function of RelaxNN framework can be presented as
\begin{equation}
    \mathcal{L}_{\text{RelaxNN}}^{\text{type3}} = \omega_{\text{residual}} \mathcal{L}_{\text{residual}}^{\text{type3}} + \omega_{\text{flux}}\mathcal{L}_{\text{flux}}^{\text{type3}} + \omega_{\text{IC}}\mathcal{L}_{\text{IC}} \, ,
\end{equation}
where the residual loss term can be presented as
\begin{equation}
\begin{split}
    \mathcal{L}_{\text{residual}}^{\text{type3}} & = \frac{\omega_{r}^{m}}{\vert \mathcal{T}_{r} \vert}  \sum_{\bm{x} \in \mathcal{T}_{r}} \left( \partial_{t} \rho^{\text{NN}}_{\bm{\theta_{1}}}(\bm{x}_{i}) + \partial_{x} (\rho^{\text{NN}}_{\bm{\theta_{1}}} u^{\text{NN}}_{\bm{\theta_{1}}}(\bm{x}_{i}) \right)^{2} \\
    & + \frac{\omega_{r}^{p}}{\vert \mathcal{T}_{r} \vert}  \sum_{\bm{x} \in \mathcal{T}_{r}} \left( \partial_{t} (\rho^{\text{NN}}_{\bm{\theta_{1}}}u^{\text{NN}}_{\bm{\theta_{1}}})(\bm{x}_{i}) + \partial_{x} \left( \rho^{\text{NN}}_{\bm{\theta_{1}}}(u^{\text{NN}}_{\bm{\theta_{1}}})^{2} + p^{\text{NN}}_{\bm{\theta_{1}}} \right)(\bm{x}_{i}) \right)^{2} \\
    & + \frac{\omega_{r}^{e}}{\vert \mathcal{T}_{r} \vert}  \sum_{\bm{x} \in \mathcal{T}_{r}} \left( \partial_{t} \left(\frac{1}{\gamma - 1}\rho^{\text{NN}}_{\bm{\theta_{1}}} + \frac{1}{2}\rho^{\text{NN}}_{\bm{\theta_{1}}}(u^{\text{NN}}_{\bm{\theta_{1}}})^{2}
    \right)(\bm{x}_{i}) + \partial_{x} \zeta^{\text{NN}}_{\bm{\theta_{2}}}(\bm{x}_{i}) \right)^{2} \, ,
\end{split}
\end{equation}
and the flux loss term can be presented as
\begin{equation}
\begin{split}
     \mathcal{L}_{\text{flux}}^{\text{type3}} & = 
     \frac{\omega_{f}^{e}}{\vert \mathcal{T}_{r} \vert}  \sum_{\bm{x} \in \mathcal{T}_{r}} \left( \zeta^{\text{NN}}_{\bm{\theta_{2}}}(\bm{x}_{i}) - \left(
     \frac{\gamma}{\gamma - 1} \rho^{\text{NN}}_{\bm{\theta_{1}}}u^{\text{NN}}_{\bm{\theta_{1}}}
     + \frac{1}{2}\rho^{\text{NN}}_{\bm{\theta_{1}}}(u^{\text{NN}}_{\bm{\theta_{1}}})^{3} \right) (\bm{x}_{i}) \right)^{2} \, .
\end{split}
\end{equation}

\section{Numerical examples}
\label{section: examples}
To illustrate the effectiveness of RelaxNN, we conduct a series of experiments about capturing shock waves. For the systems of conservation laws with an initial condition in the finite domain : 
\begin{equation}
\left\{
\begin{aligned}
     & \partial_{t}\bm{u} + \partial_{x} \bm{F}(\bm{u}
) = 0, \quad (t,x) \in \Omega ,\\ 
     & \bm{u}(t,x) = \bm{u}_{0}(t,x), \quad (t,x) \in \partial \Omega ,
\end{aligned}
\right.
\end{equation}
we impose the boundary condition as the same as the initial condition on the assumption that waves do not propagate to the boundary in the finite period we considered. In this case, the loss function for PINN and RelaxNN will be :
\begin{subequations}\label{loss function comparison}
\begin{align}
    \mathcal{L}_{\text{PINN}} & = \omega_{\text{PDE}} \mathcal{L}_{\text{PDE}} + \omega_{\text{IC}} \mathcal{L}_{\text{IC}} + \omega_{\text{BC}} \mathcal{L}_{\text{BC}}\, , \\
    \mathcal{L}_{\text{RelaxNN}} & = \omega_{\text{residual}} \mathcal{L}_{\text{residual}} + \omega_{\text{flux}}\mathcal{L}_{\text{flux}} + \omega_{\text{IC}}\mathcal{L}_{\text{IC}} + \omega_{\text{BC}} \mathcal{L}_{\text{BC}} \, .
\end{align}
\end{subequations}
Here $\mathcal{L}_{\text{BC}}$ is defined as same as $\mathcal{L}_{\text{IC}}$ excepts $ \bm{x} \in \mathcal{T}_{bc}$, $\mathcal{T}_{bc} \subset \partial \Omega$ are the residual points sampled on the boundary.

For all experiments, training points are sampled randomly in uniform distribution. For conciseness, we roughly set $\vert \mathcal{T}_{r} \vert = 2540$,$\quad \vert \mathcal{T}_{ic} \vert = 320$ and $\quad \vert \mathcal{T}_{bc} \vert = 160$ for all experiments. Also, we set the random seed to $1$ for all packages to ensure repeatability. For the architecture of the neural network, we employ the fully connected deep neural networks (DNNs) with ``tanh'' activation functions and He uniform initialization \cite{he2015delving} for all experiments. 
For the training settings, we adopt ``Adam''\cite{kingma2015adam} Optimizer with an initial learning rate $1e-3$, decaying exponentially for every $1000$ epochs with a decay rate $0.99$ during the whole training process for all experiments.
We compare our results with high-resolution numerical methods obtained by clawpack~\cite{clawpack}. Besides the remarkable performance, the comparison among the several types of RelaxNN also confirms our hypothesis about the generation of spurious waves. Later, we will see that the mild oscillation occurs only on the type1 RelaxNN schemes for shallow water equations, as well as type1 and type2 RelaxNN schemes for Euler's equations although both RelaxNN schemes are able to capture the shock waves and acquire outstanding performance. This interesting phenomenon also embodies Occam's razor theorem.

\begin{table}
\begin{tabular}{l*{9}{c}}
\toprule
 \multirow{2}{*}{Examples}& \multirow{2}{*}{\shortstack{Relax \\ level}} &  \multicolumn{4}{c}{$\bm{u}_{\bm{\theta_{1}}}^{\text{NN}}$} &  \multicolumn{4}{c}{$\bm{v}_{\bm{\theta_{2}}}^{\text{NN}}$} \\
\cmidrule(r){3-6} \cmidrule(r){7-10}
  &   & $d_{\text{in}}$ & $N_{\text{depth}}$ & $N_{\text{width}}$ & $d_{\text{out}}$ & $d_{\text{in}}$ & $N_{\text{depth}}$ & $N_{\text{width}}$ & $d_{\text{out}}$\\
\midrule
  Burgers' Equation & &  \multirow{6}{*}{2} & 4 & 128 & 1 & \multirow{6}{*}{2} & 4 & 64 & 1  \\ 
  \cline{1-1} \cline{4-6} \cline{8-10}
  \multirow{2}{*}{Shallow Water Equations}& type1 &   & 5 & 128 & 2 &  & 5 & 128 & 2 \\
    & type2 &  & 5 & 128 & 2 &  & 5 & 64 & 1\\
  \cline{1-1} \cline{4-6} \cline{8-10}
  \multirow{3}{*}{Euler Equations} & type1 &  & 6 & 384 & 3 &  & 6 & 384 & 3 \\
    & type2 &  & 6 & 384 & 3 &  & 6 & 256 & 2 \\
    & type3 &  & 6 & 384 & 3 &  & 6 & 128 & 1 \\
\bottomrule
\end{tabular}
\caption{Networks' configuration for all examples. Here $d_{\text{in}}$, $d_{\text{out}}$, $N_{\text{depth}}$,  $N_{\text{width}}$ represent the input dimension, the output dimension, the depth and width of the hidden layers of the specified neural network. Here $d_{\text{out}}$ of $\bm{v}_{\bm{\theta_{2}}}^{\text{NN}}$ indicates that how many conservation laws we relaxed.}
\label{tab: Networks Config}
\end{table}

\subsection{Inviscid Burgers' equation}
\subsubsection{Riemann problem}
We consider the inviscid Burgers' equation in domain $ \Omega = \{ (t,x)\} = [ 0.0,1.0]\times [-0.6,0.6 ]$ with Riemann initial condition as
\begin{equation}
    u_{0}(0,x) = \left\{
    \begin{aligned}
        & 1.0 , \quad -0.6 \leq x \leq 0 \\ 
        & 0.0 , \quad 0 < x \leq 0.6
    \end{aligned}
    \right. .
\end{equation}

\begin{table}
\centering
\begin{tabular}{l*{5}{c}}
\toprule
 \multirow{2}{*}{Problems}& \multirow{2}{*}{\shortstack{Relax \\ level}} & $\omega_{\text{residual}} $ & $\omega_{\text{flux}}$ & $\omega_{\text{IC}}$ & $\omega_{\text{BC}}$ \\
 & & & & & \\
\midrule
  Riemann Problem & & 0.1 & 2.0 & 10 & 10 \\
  Sine Problem & & 0.5 & 2.0 & 5.0 & 5.0 \\
\bottomrule
\end{tabular}
\caption{\textit{Burgers' equation}: Loss weights settings}
\label{tab: w settings Burgers}
\end{table}

\begin{figure}[htbp]
\centering
\subfloat[riemann]
{\label{fig:riemann burgers total loss}
\includegraphics[width = 0.48\textwidth]{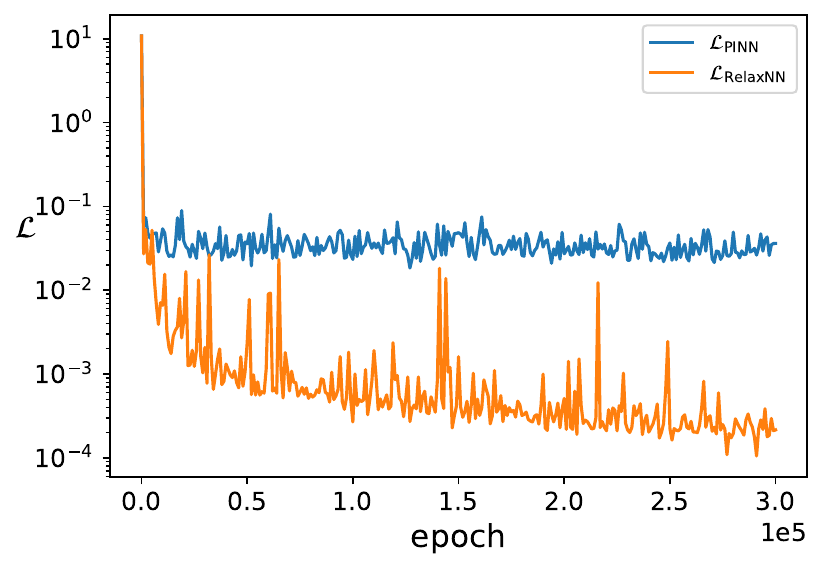}} 
\subfloat[sine]
{\label{fig:sine burgers total loss}
\includegraphics[width = 0.48\textwidth]{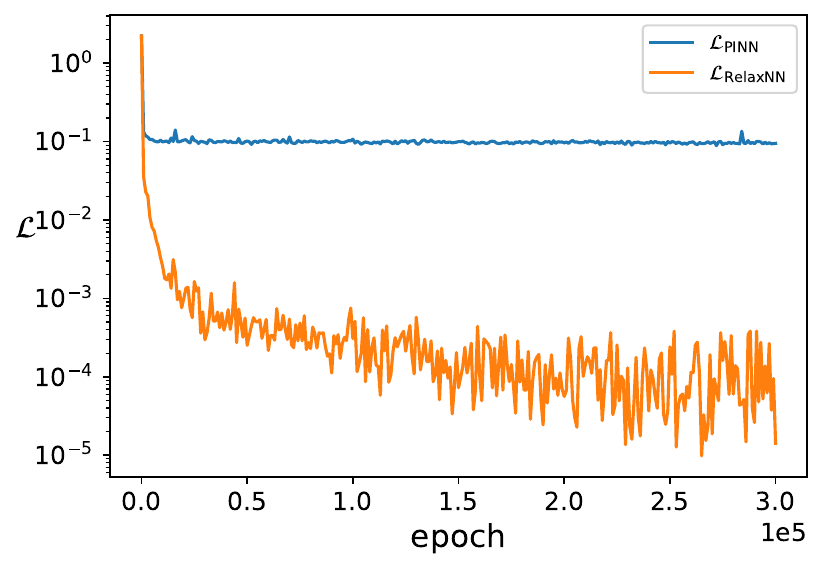}} 
\caption{\textit{Burgers' equation}: Comparison of the loss curves of RelaxNN against the one of PINN during the training process: (a)Riemann initial condition. (b) Sine initial condition.}
\end{figure}

\begin{figure}[htbp]
\centering
\subfloat[$t=0.0$]{\includegraphics[width = 0.48\textwidth]{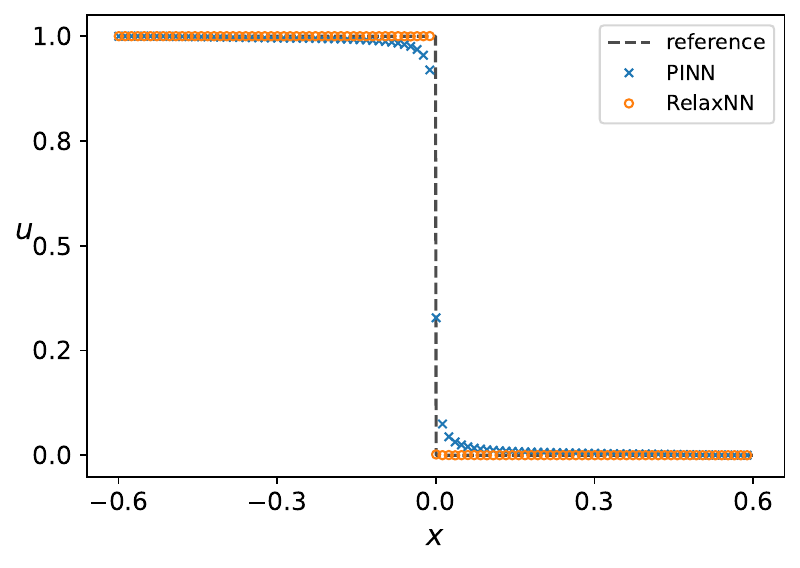}}
\subfloat[$t=0.2$]{\includegraphics[width = 0.48\textwidth]{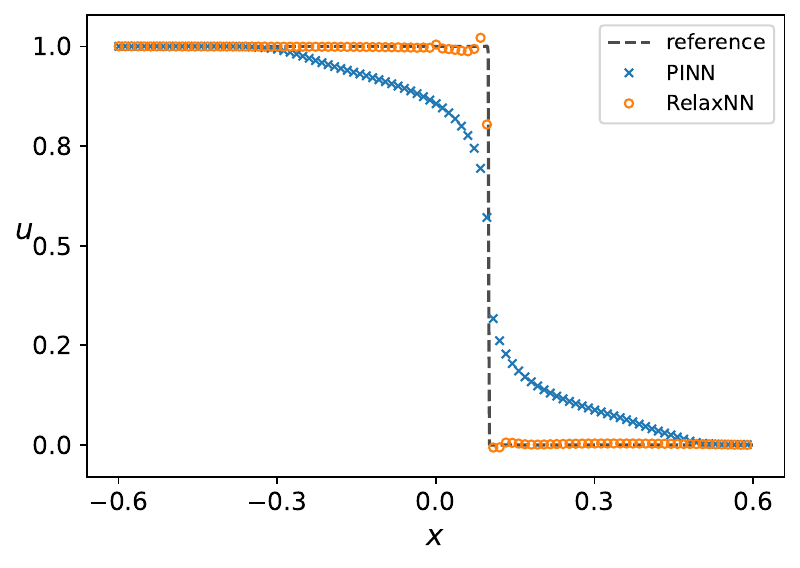}}
\\
\subfloat[$t=0.4$]{\includegraphics[width = 0.48\textwidth]{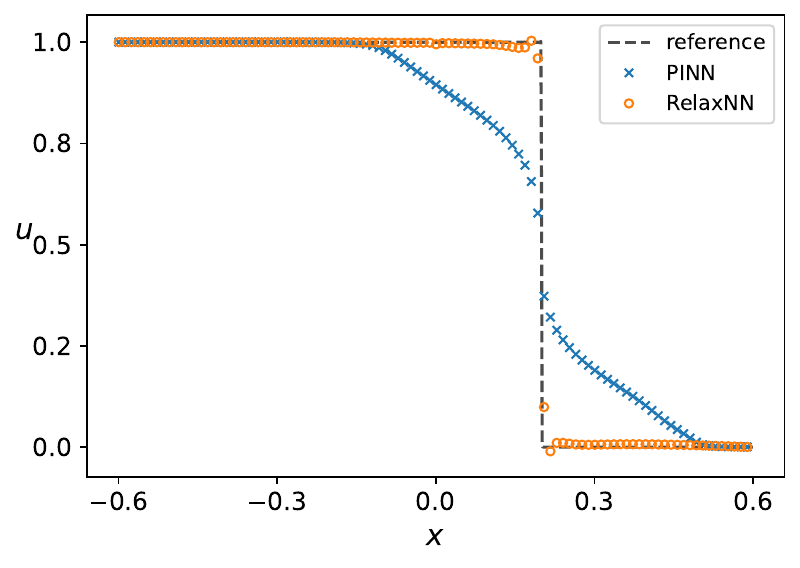}}
\subfloat[$t=0.6$]{\includegraphics[width = 0.48\textwidth]{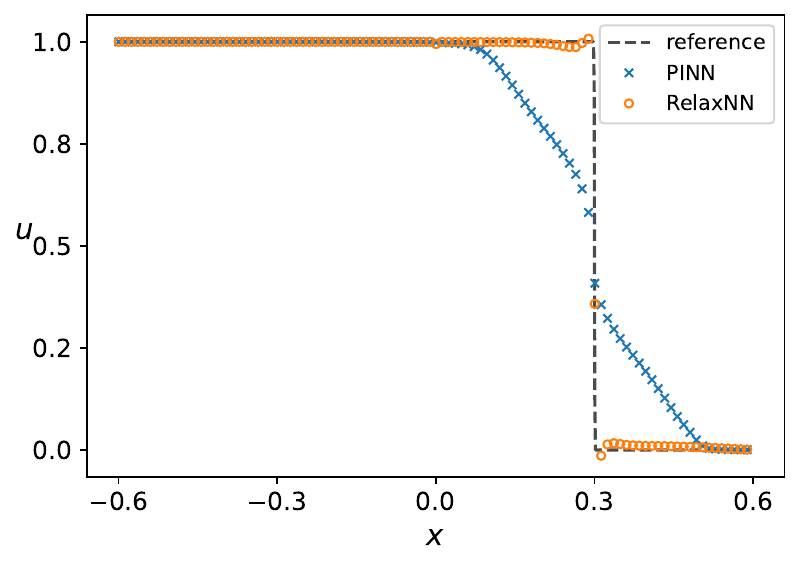}}
\caption{\textit{Burgers' equation with Riemann initial condition}: Comparison among the final epoch prediction of RelaxNN, PINN, and the reference solution spatially at different specific moments. The resulting relative $L^{2}$ error of PINN, RelaxNN are $1.64 \times 10^{-2}$, $8.29 \times 10^{-4}$. For RelaxNN, the configuration of $\bm{u}_{\bm{\theta_{1}}}^{\text{NN}}$, $\bm{v}_{\bm{\theta_{2}}}^{\text{NN}}$ are [2,128,128,128,128,1] and [2,64,64,64,64,1]. We training for 300,000 epochs and one step per epoch. Loss weights settings are shown in ~\Cref{tab: w settings Burgers}. For PINN, all settings are the same as RelaxNN without the extra neural network.}
\label{fig: slice burgers riemann}
\end{figure}

Networks' configurations are presented in \Cref{tab: Networks Config} and weights settings are shown in \Cref{tab: w settings Burgers}. With these settings, we train RelaxNN for $300,000$ epochs and one step for every epoch.  For the PINN framework, settings are the same without the extra neural network $\bm{v}_{\bm{\theta_{2}}}^{\text{NN}}$. In ~\Cref{fig:riemann burgers total loss}, We compare the total loss between the PINN and RelaxNN. In ~\Cref{fig: slice burgers riemann}, we compare the prediction of RelaxNN and PINN at the final epoch for some specific times.

\subsubsection{Sine problem}
\begin{figure}[htbp]
\centering
\subfloat[$t=0.0$]{\includegraphics[width = 0.4\textwidth]{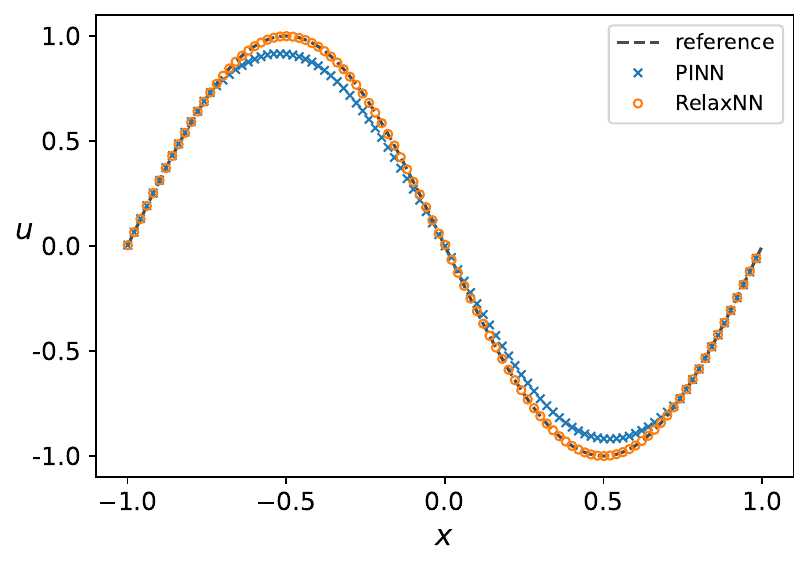}}
\subfloat[$t=0.2$]{\includegraphics[width = 0.4\textwidth]{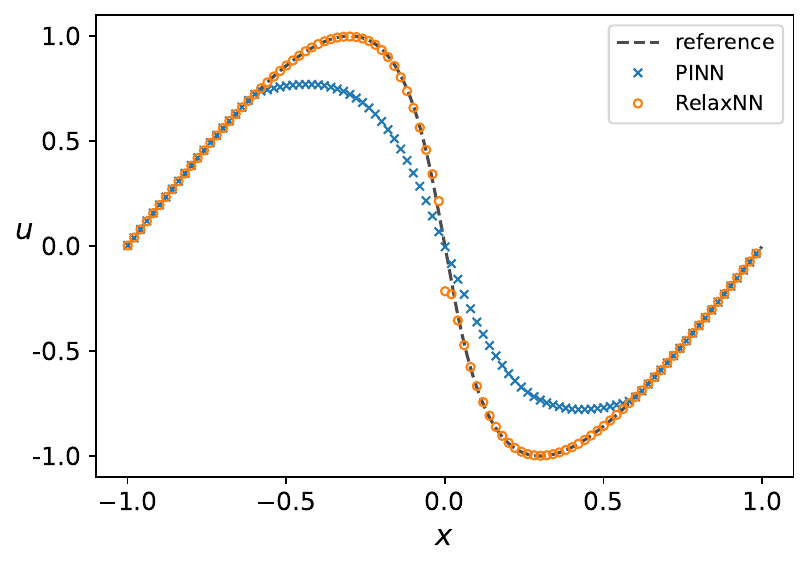}}
\\
\subfloat[$t=0.4$]{\includegraphics[width = 0.4\textwidth]{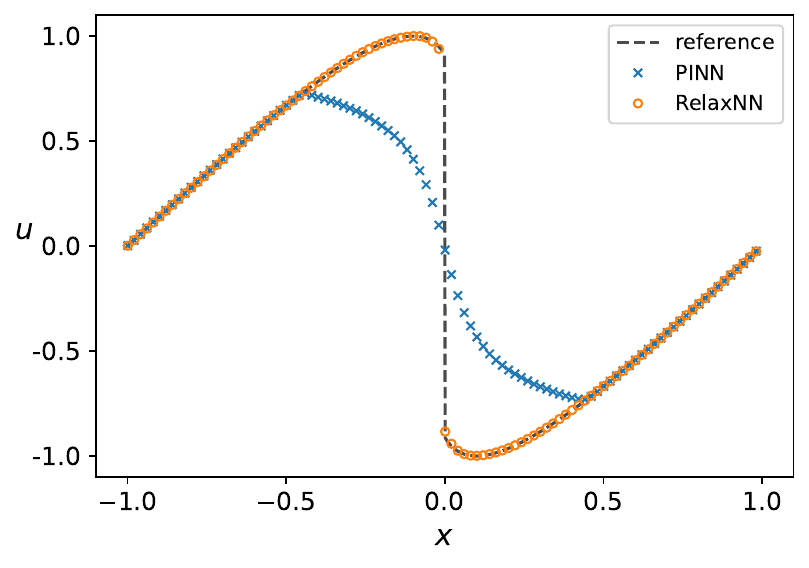}}
\subfloat[$t=0.6$]{\includegraphics[width = 0.4\textwidth]{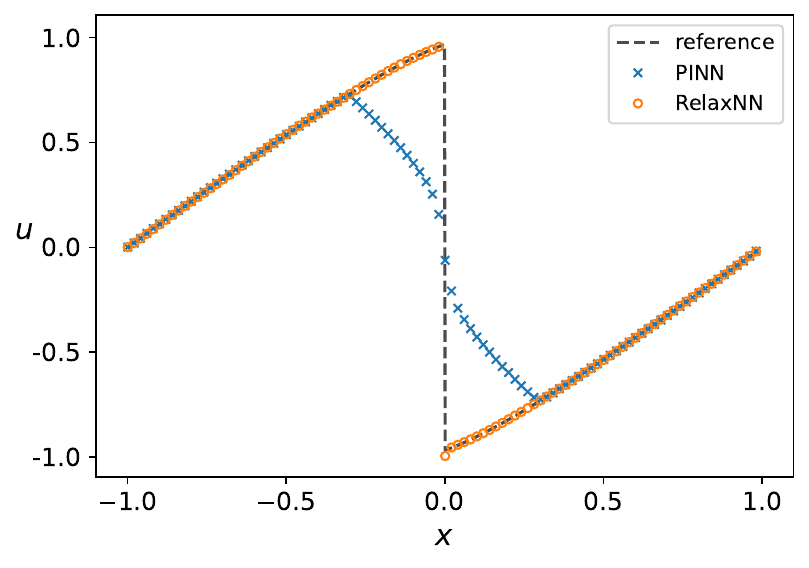}}
\caption{\textit{Burgers' equation with Sine initial condition}: Comparison among the final epoch prediction of RelaxNN, PINN, and the reference solution spatially at different specific moments. The resulting relative $L^{2}$ error of PINN, RelaxNN are $1.01 \times 10^{-1} $,$3.63 \times 10^{-4}$. For RelaxNN, the configuration of $\bm{u}_{\bm{\theta_{1}}}^{\text{NN}}$, $\bm{v}_{\bm{\theta_{2}}}^{\text{NN}}$ are [2,128,128,128,128,1] and [2,64,64,64,64,1]. We training for 300,000 epochs and one step per epoch. Loss weights settings are shown in ~\Cref{tab: w settings Burgers}. For PINN, all settings are the same as RelaxNN without the extra neural network.}
\label{fig: slice burgers sine}
\end{figure}

We consider the inviscid Burgers' equation in domain $ \Omega = \{ (t,x)\} = [ 0.0,1.0]\times [-1.0,1.0 ]$ with sine initial condition as
\begin{equation}
    u_{0}(0,x) = - \sin (\pi x), \quad  -1 \leq x \leq 1 \, .
\end{equation}
Networks' configurations are presented in \Cref{tab: Networks Config} and weights settings are shown in \Cref{tab: w settings Burgers}. With these settings, we train RelaxNN for $300,000$ epochs and one step for every epoch.  For the PINN framework, settings are the same without the extra neural network $\bm{v}^{\text{NN}}_{\beta}$. In ~\Cref{fig:sine burgers total loss}, We compare the total loss between the PINN and RelaxNN. In ~\Cref{fig: slice burgers sine}, we compare the prediction of RelaxNN and PINN at the final epoch for some specific times. 

\subsection{Shallow water equations}
\subsubsection{Dam break problem}
We consider the shallow water equations in domain $ \Omega = \{ (t,x)\} = [ 0.0,1.0]\times [-1.5,1.5]$ with dam break initial condition as
\begin{equation}
    \bm{u}_{0}(0,x) = \left( 
    \begin{array}{cc}
         1.0 \\ 
         0.0 
    \end{array}
    \right) , \quad \text{if} \quad -1.5 \leq x \leq 0.0 ;
    \quad
    \bm{u}_{0}(0,x) = \left( 
    \begin{array}{cc}
         0.5 \\
         0.0 
    \end{array}
    \right) , \quad \text{if} \quad 0.0 < x \leq 1.5 ;
\end{equation}

\begin{figure}[htbp]
\centering
\subfloat[dam break]
{\label{fig:dam swe total loss}
\includegraphics[width = 0.4\textwidth]{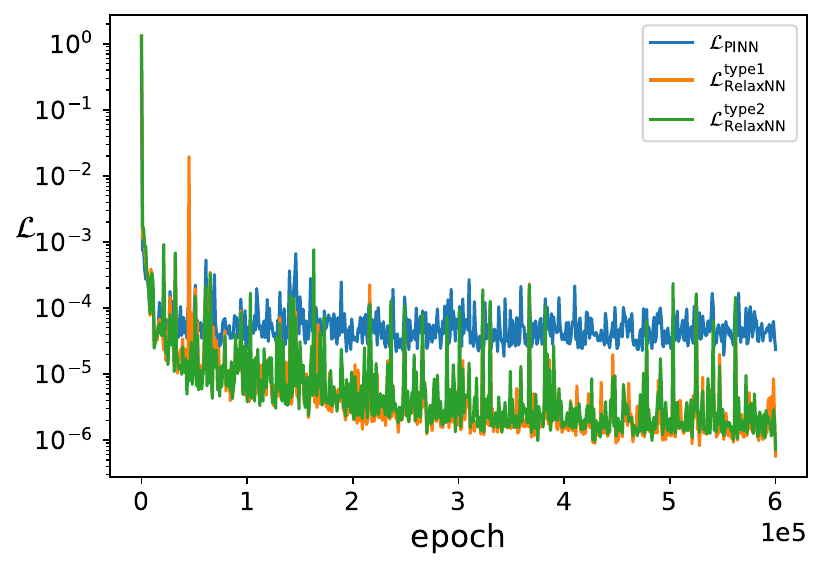}} 
\subfloat[two shock]
{\label{fig:2shock swe total loss}
\includegraphics[width = 0.4\textwidth]{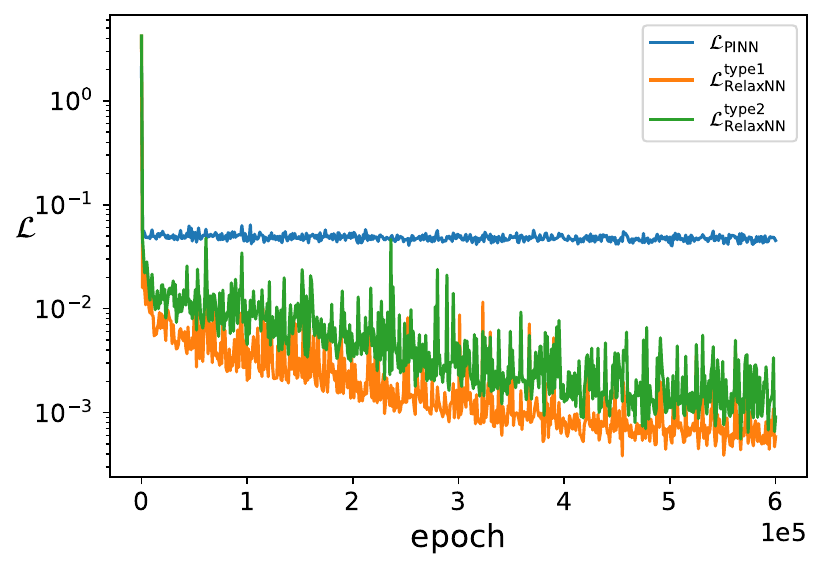}} 
\caption{\textit{Shallow water equations}: Comparison of the loss curves of RelaxNN against the one of PINN during the training process. (a)Dam break initial condition. (b)Two shock initial condition.}
\end{figure}

\begin{figure}[htbp]
\centering
\subfloat[$t=0.0$]{\includegraphics[width = 0.4\textwidth]{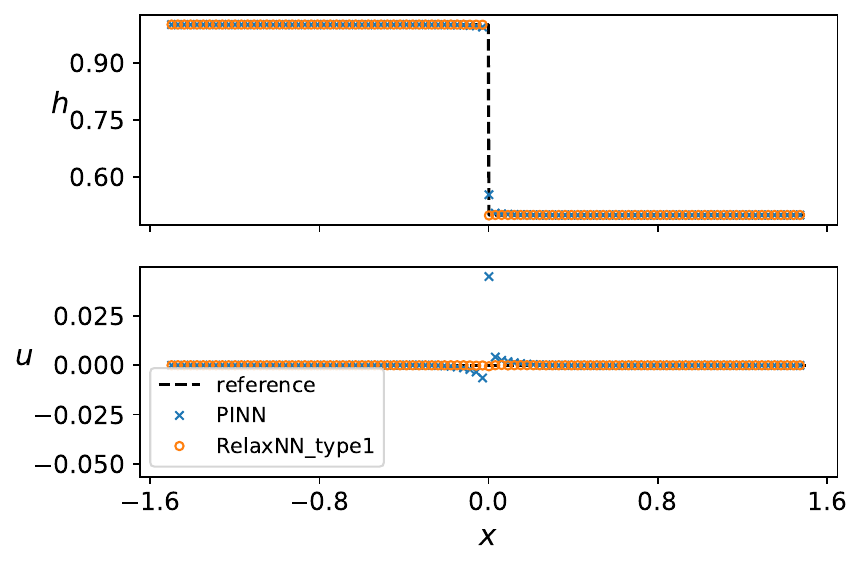}}
\subfloat[$t=0.2$]{\includegraphics[width = 0.4\textwidth]{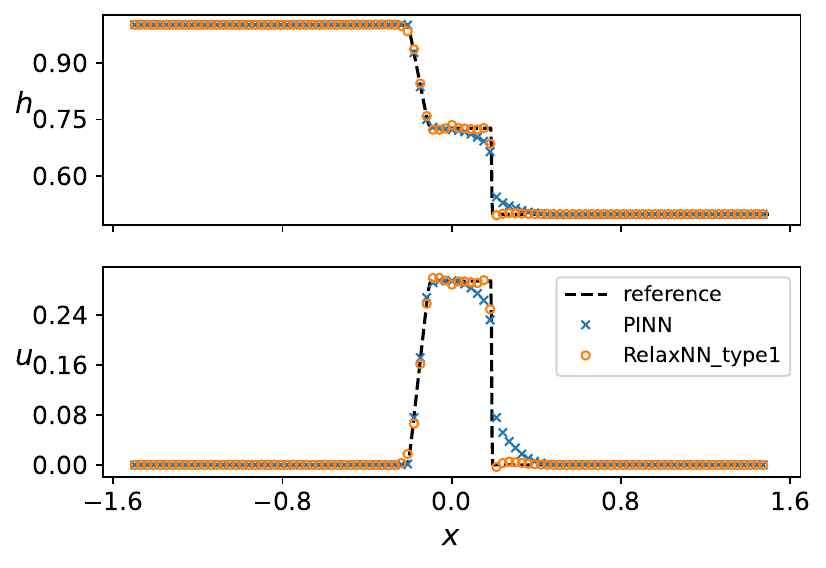}}
\\
\subfloat[$t=0.4$]{\includegraphics[width = 0.4\textwidth]{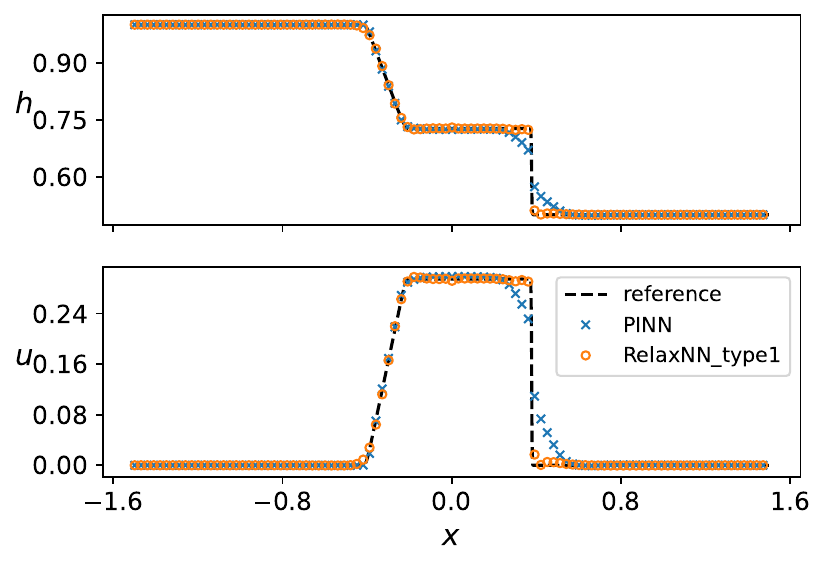}}
\subfloat[$t=0.6$]{\includegraphics[width = 0.4\textwidth]{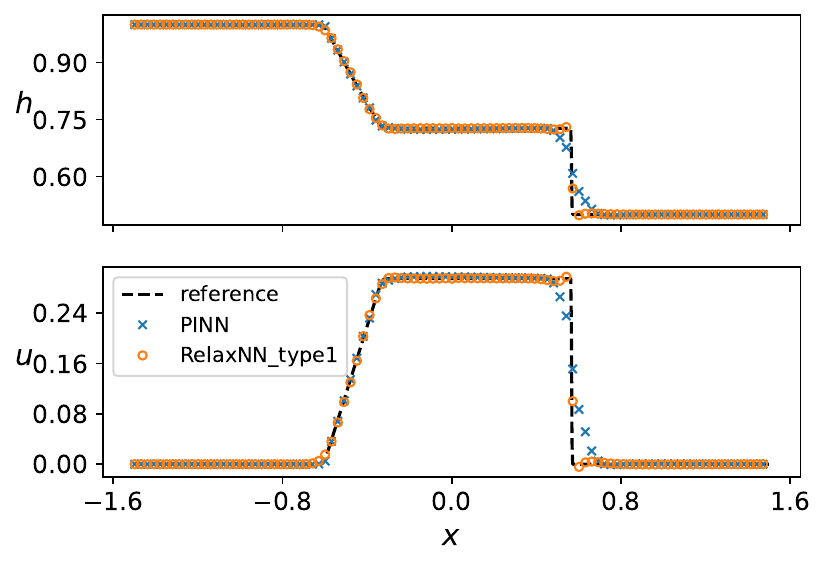}}
\caption{\textit{Shallow water equations with dam break initial condition}: Comparison among the final epoch prediction of RelaxNN(type1), PINN, and the reference solution spatially at different specific moments. The resulting relative $L^{2}$ error of PINN, RelaxNN are $4.85 \times 10^{-4}$, $9.01 \times 10^{-5}$. For RelaxNN, the configuration of $\bm{u}_{\bm{\theta_{1}}}^{\text{NN}}$, $\bm{v}_{\bm{\theta_{2}}}^{\text{NN}}$ are [2,128,128,128,128,128,2] and [2,128,128,128,128,128,2]. We training for 600,000 epochs and one step per epoch. Loss weights settings are shown in ~\Cref{tab: w settings swe}. For PINN, all settings are the same as RelaxNN without the extra neural network.}
\label{fig: slice swe dam type1}
\end{figure}

\begin{figure}[htbp]
\centering
\subfloat[$t=0.0$]{\includegraphics[width = 0.4\textwidth]{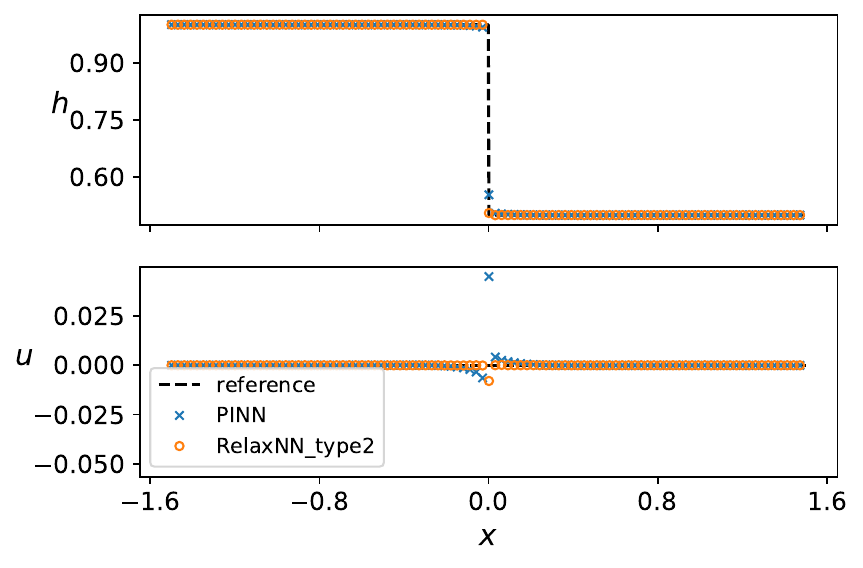}}
\subfloat[$t=0.2$]{\includegraphics[width = 0.4\textwidth]{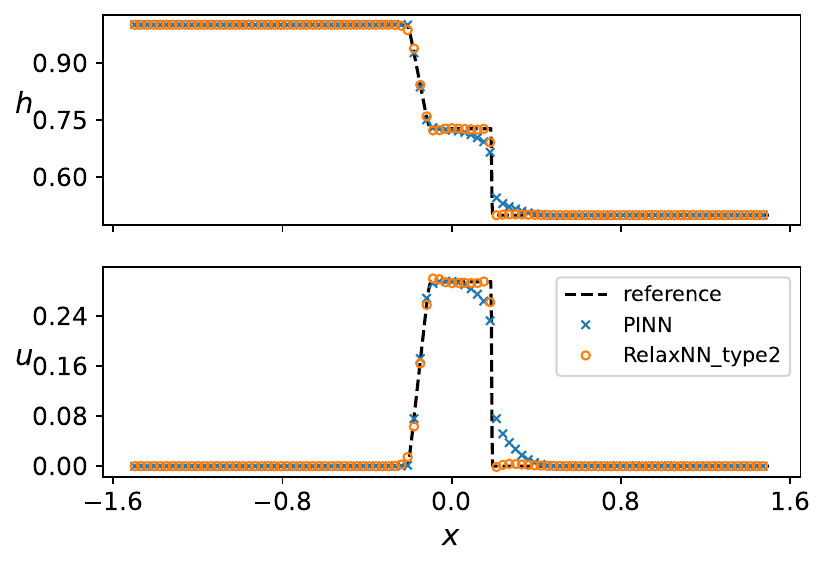}}
\\
\subfloat[$t=0.4$]{\includegraphics[width = 0.4\textwidth]{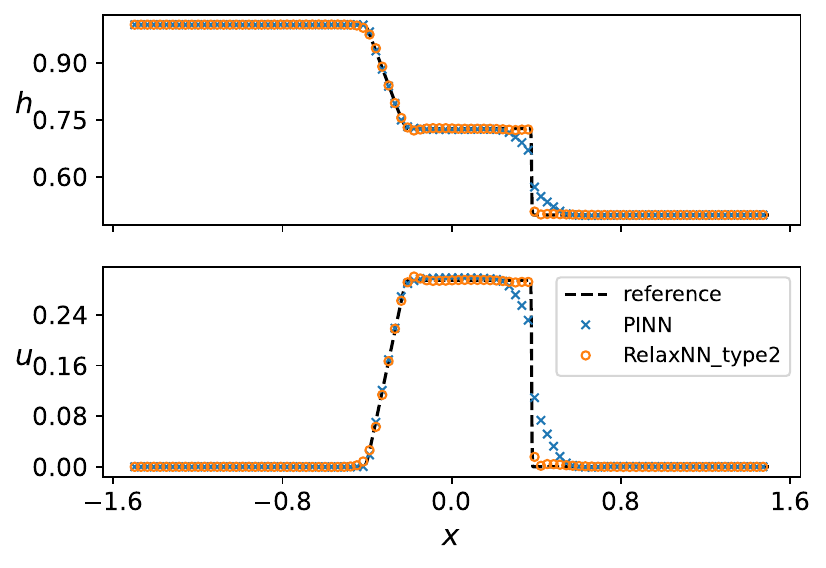}}
\subfloat[$t=0.6$]{\includegraphics[width = 0.4\textwidth]{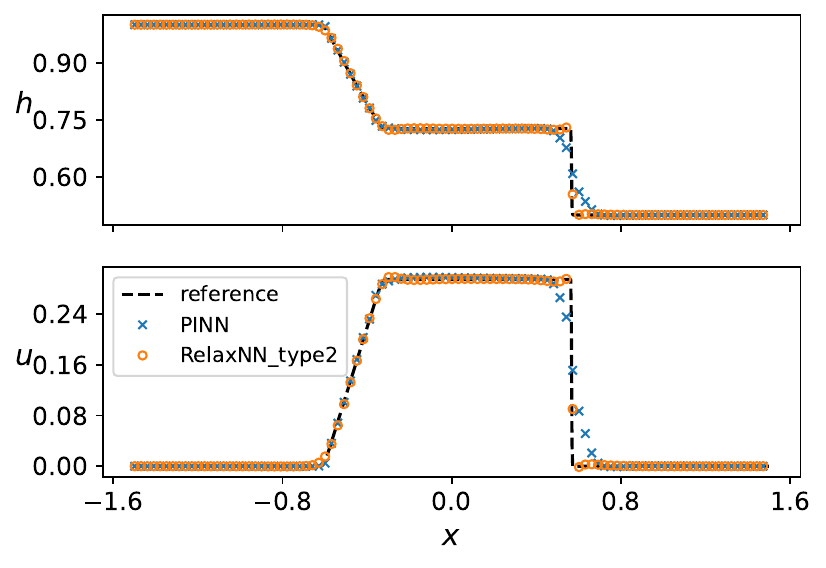}}
\caption{\textit{Shallow water equations with dam break initial condition}: Comparison among the final epoch prediction of RelaxNN(type2), PINN, and the reference solution spatially at different specific moments. The resulting relative $L^{2}$ error of PINN, RelaxNN are $4.85 \times 10^{-4}$, $7.28 \times 10^{-5}$. For RelaxNN, the configuration of $\bm{u}_{\bm{\theta_{1}}}^{\text{NN}}$, $\bm{v}_{\bm{\theta_{2}}}^{\text{NN}}$ are [2,128,128,128,128,128,2] and [2,64,64,64,64,64,1]. We training for 600,000 epochs and one step per epoch. Loss weights settings are shown in ~\Cref{tab: w settings swe}. For PINN, all settings are the same as RelaxNN without the extra neural network.}
\label{fig: slice swe dam type2}
\end{figure}

\begin{table}
\centering
\begin{tabular}{l*{7}{c}}
\toprule
 \multirow{2}{*}{Problems}& \multirow{2}{*}{\shortstack{Relax \\ level}} & \multicolumn{2}{c}{$\omega_{\text{residual}}$} & \multicolumn{2}{c}{$\omega_{\text{flux}}$} & \multirow{2}{*}{$\omega_{\text{IC}}$}& \multirow{2}{*}{$\omega_{\text{BC}}$} \\
\cmidrule(r){3-4} \cmidrule(r){5-6}
  &   & $\omega_{r}^{m}$ & $\omega_{r}^{p}$ & $\omega_{f}^{m}$ & $\omega_{f}^{p}$ &  & \\
\midrule
  Dam break & type1 & 0.01 & 0.01 & 1.0 & 1.0 & 1.0 & 1.0 \\ 
  Dam break& type2 & 0.01 & 0.01 &  & 1.0 & 1.0 & 1.0 \\
  Two shock& type1 & 0.1 & 0.1 & 1.0 & 1.0 & 1.0 & 1.0 \\ 
  Two shock& type2 & 0.1 & 0.1 &  & 1.0 & 1.0 & 1.0 \\
\bottomrule
\end{tabular}
\caption{\textit{Shallow water equations}: Loss weights settings}
\label{tab: w settings swe}
\end{table}
 Networks' configurations are presented in \Cref{tab: Networks Config} and weights settings are shown in \Cref{tab: w settings swe}. With these settings, we train RelaxNN for $600,000$ epochs and one step for every epoch.  For the PINN framework, settings are the same without the extra neural network $\bm{v}_{\bm{\theta_{2}}}^{\text{NN}}$. In ~\Cref{fig:dam swe total loss}, We compare the total loss between the PINN and RelaxNN. In ~\Cref{fig: slice swe dam type1}, we compare the prediction of type1 RelaxNN and PINN at the final epoch for some specific times. Also for type2 RelaxNN in ~\Cref{fig: slice swe dam type2}.
 
\subsubsection{Two shock problem}
\begin{figure}[htbp]
\centering
\subfloat[$t=0.0$]{\includegraphics[width = 0.4\textwidth]{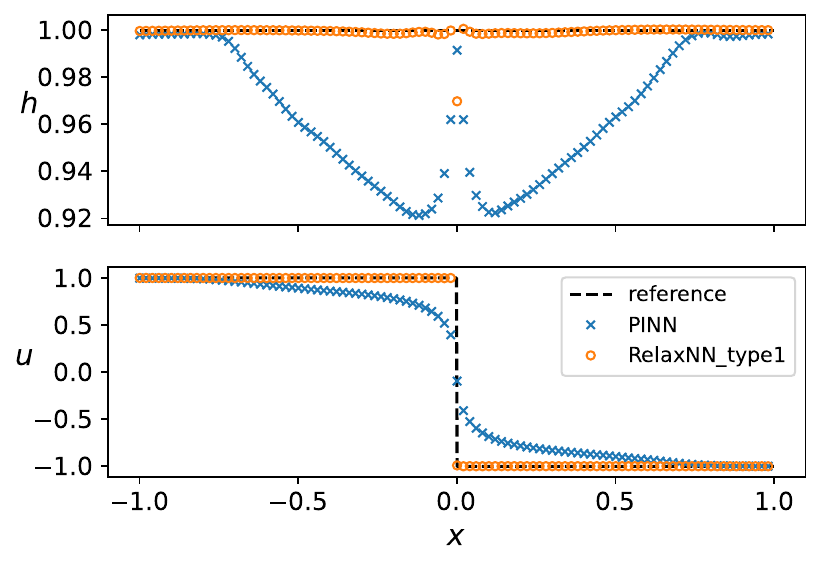}}
\subfloat[$t=0.2$]{\includegraphics[width = 0.4\textwidth]{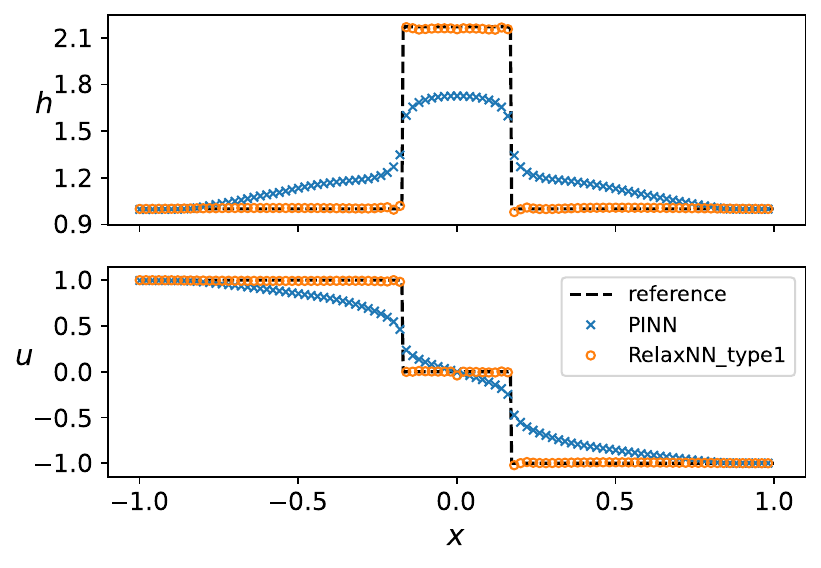}}
\\
\subfloat[$t=0.4$]{\includegraphics[width = 0.4\textwidth]{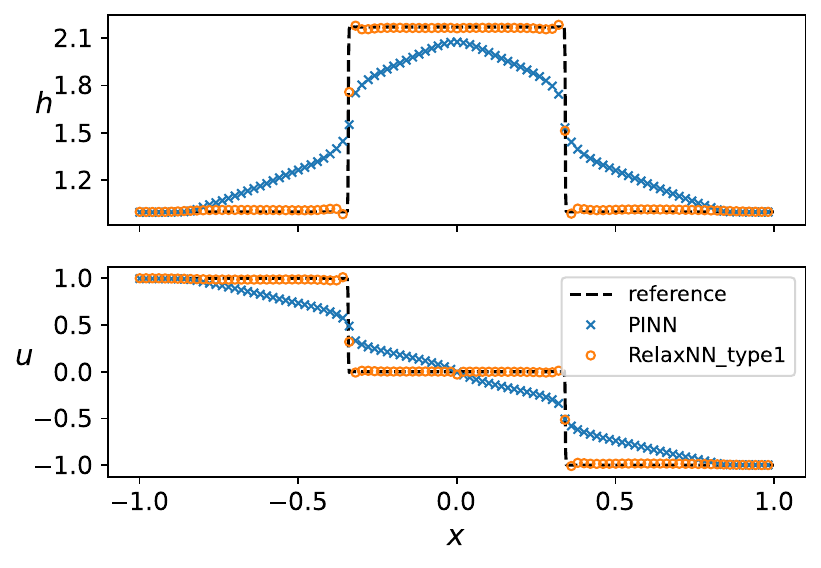}}
\subfloat[$t=0.6$]{\includegraphics[width = 0.4\textwidth]{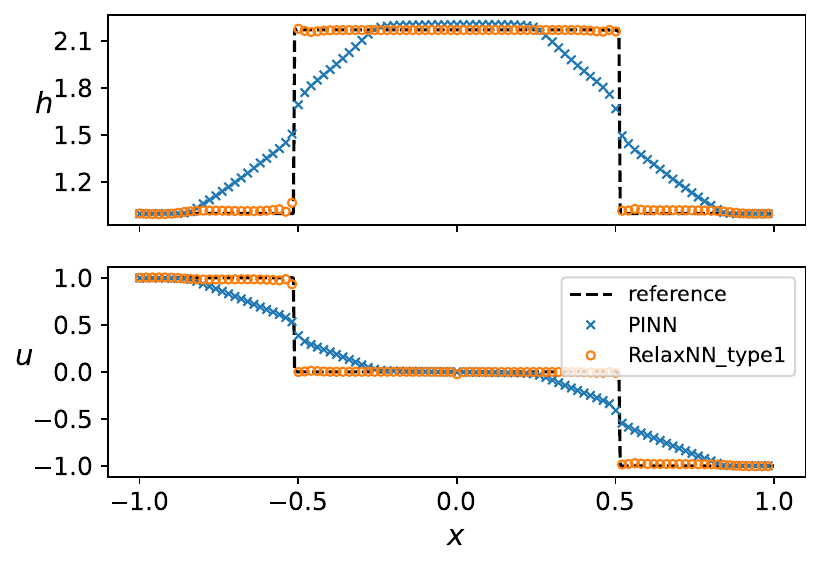}}
\caption{\textit{Shallow water equations with two shock initial condition}: Comparison among the final epoch prediction of RelaxNN(type1), PINN and the reference solution spatially at different specific moments. The resulting relative $L^{2}$ error of PINN, RelaxNN are $2.12 \times 10^{-2}$, $5.86 \times 10^{-4}$. For RelaxNN, the configuration of $\bm{u}_{\bm{\theta_{1}}}^{\text{NN}}$, $\bm{v}_{\bm{\theta_{2}}}^{\text{NN}}$ are [2,128,128,128,128,128,2] and [2,128,128,128,128,128,2]. We training for 600,000 epochs and one step per epoch. Loss weights settings are shown in ~\Cref{tab: w settings swe}. For PINN, all settings are the same as RelaxNN without the extra neural network.}
\label{fig: slice swe 2shock type1}
\end{figure}

\begin{figure}[htbp]
\centering
\subfloat[$t=0.0$]{\includegraphics[width = 0.4\textwidth]{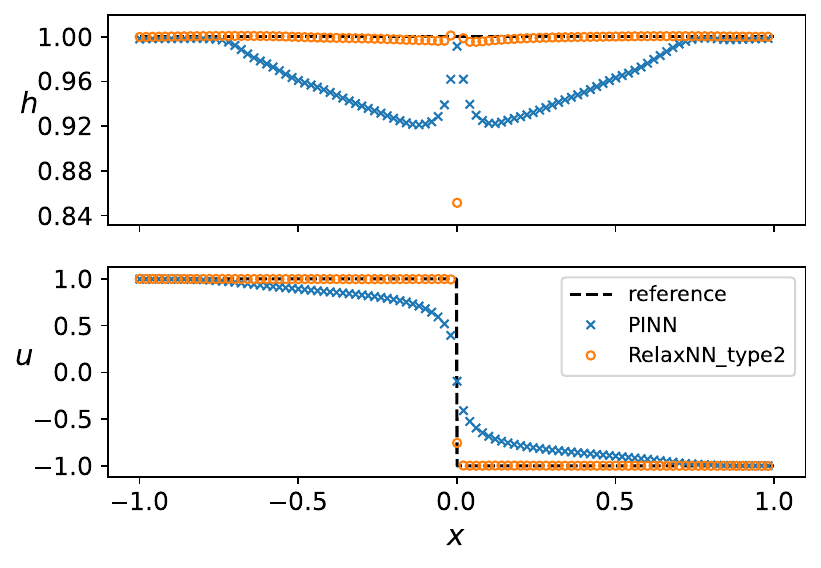}}
\subfloat[$t=0.2$]{\includegraphics[width = 0.4\textwidth]{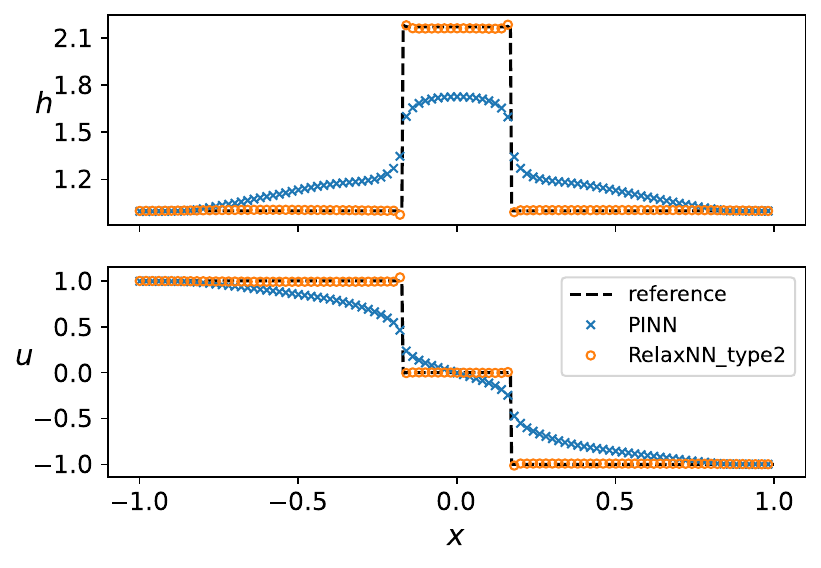}}
\\
\subfloat[$t=0.4$]{\includegraphics[width = 0.4\textwidth]{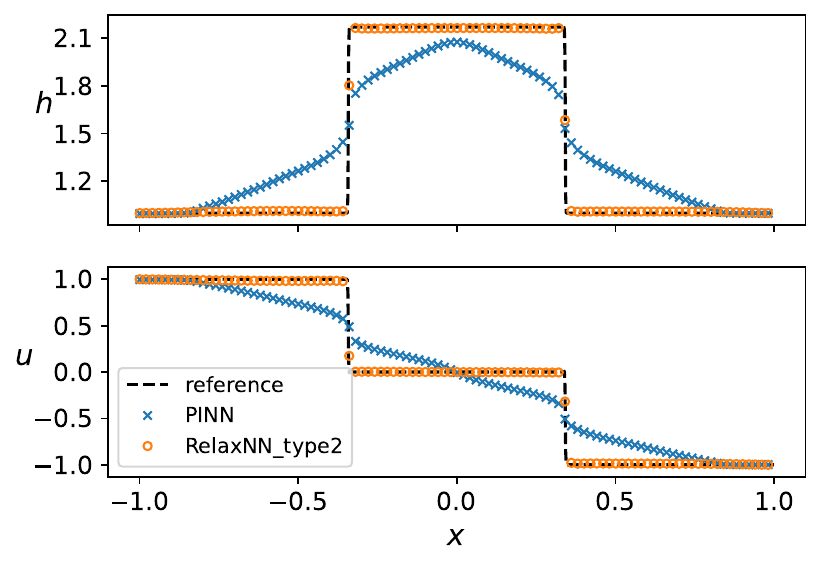}}
\subfloat[$t=0.6$]{\includegraphics[width = 0.4\textwidth]{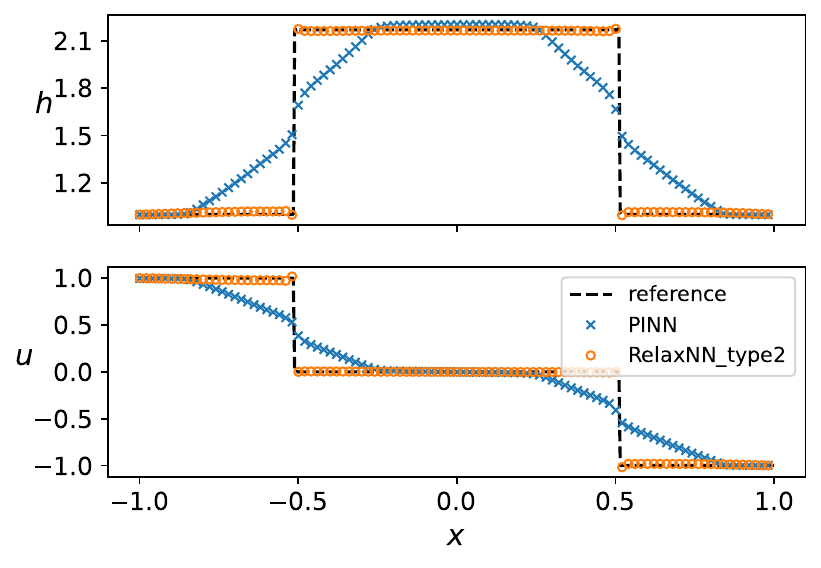}}
\caption{\textit{Shallow water equations with two shock initial condition}: Comparison among the final epoch prediction of RelaxNN(type2), PINN, and the reference solution spatially at different specific moments. The resulting relative $L^{2}$ error of PINN, RelaxNN are $2.12 \times 10^{-2}$, $2.53 \times 10^{-4}$. For RelaxNN, the configuration of $\bm{u}_{\bm{\theta_{1}}}^{\text{NN}}$, $\bm{v}_{\bm{\theta_{2}}}^{\text{NN}}$ are [2,128,128,128,128,128,2] and [2,64,64,64,64,64,1]. We training for 600,000 epochs and one step per epoch. Loss weights settings are shown in ~\Cref{tab: w settings swe}. For PINN, all settings are the same as RelaxNN without the extra neural network.}
\label{fig: slice swe 2shock type2}
\end{figure}

We consider the shallow water equations in domain $ \Omega = \{ (t,x)\} = [ 0.0,1.0]\times [-1.0,1.0]$ with two shock initial condition as
\begin{equation}
    \bm{u}_{0}(0,x) = \left( 
    \begin{array}{cc}
         1.0 \\ 
         1.0 
    \end{array}
    \right) , \quad \text{if} \quad -1.0 \leq x \leq 0.0 ;
    \quad
    \bm{u}_{0}(0,x) = \left( 
    \begin{array}{cc}
         1.0 \\
         -1.0
    \end{array}
    \right) , \quad \text{if} \quad 0.0 < x \leq 1.0 ;
\end{equation}
Networks' configurations are presented in \Cref{tab: Networks Config} and weights settings are shown in \Cref{tab: w settings swe}. With these settings, we train RelaxNN for $600,000$ epochs and one step for every epoch.  For the PINN framework, settings are the same without the extra neural network $\bm{v}_{\bm{\theta_{2}}}^{\text{NN}}$. In ~\Cref{fig:2shock swe total loss}, We compare the total loss between the PINN and RelaxNN. In ~\Cref{fig: slice swe 2shock type1}, we compare the prediction of type1 RelaxNN and PINN at the final epoch for some specific times. Also for type2 RelaxNN in ~\Cref{fig: slice swe 2shock type2}.

\subsection{Euler equations}
\subsubsection{Sod problem}
We consider the euler equations in domain $ \Omega = \{ (t,x)\} = [ 0.0,0.4]\times [-0.8,0.8]$ with sod shocktube initial condition as
\begin{equation}
    \bm{u}_{0}(0,x) = \left( 
    \begin{array}{cc}
         1.0 \\ 
         0.0 \\
         1.0
    \end{array}
    \right) , \quad \text{if} \quad x \leq 0.0 ;
    \quad
    \bm{u}_{0}(0,x) = \left( 
    \begin{array}{cc}
         0.125 \\
         0.0 \\
         0.1
    \end{array}
    \right) , \quad \text{if} \quad x > 0.0 ;
\end{equation}

\begin{figure}[htbp]
\centering
\subfloat[sod shock tube]
{\label{fig:sod euler total loss}
\includegraphics[width = 0.4\textwidth]{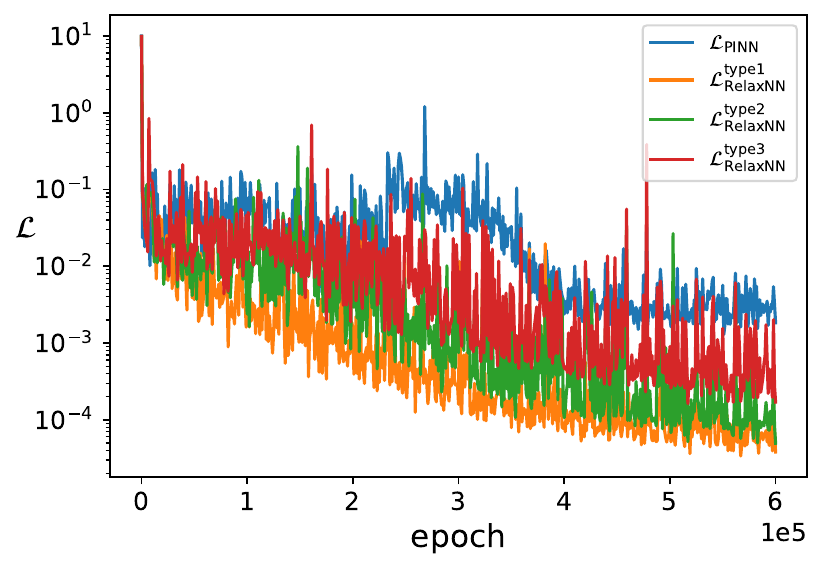}}
\subfloat[lax shock tube]
{\label{fig:lax euler total loss}
\includegraphics[width = 0.4\textwidth]{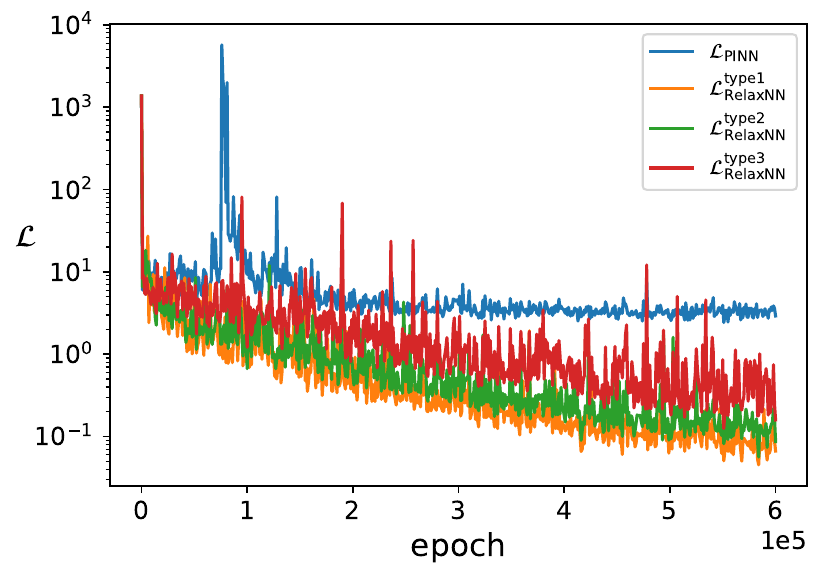}} 
\caption{loss curves during the training process of Euler equations}
\end{figure}

\begin{table}
\centering
\begin{tabular}{l*{9}{c}}
\toprule
 \multirow{2}{*}{Problems}& \multirow{2}{*}{\shortstack{Relax \\ level}} & \multicolumn{3}{c}{$\omega_{\text{residual}}$} & \multicolumn{3}{c}{$\omega_{\text{flux}}$} & \multirow{2}{*}{$\omega_{\text{IC}}$}& \multirow{2}{*}{$\omega_{\text{BC}}$} \\
\cmidrule(r){3-5} \cmidrule(r){6-8}
  &   & $\omega_{r}^{m}$ & $\omega_{r}^{p}$ & $\omega_{r}^{e}$& $\omega_{f}^{m}$ & $\omega_{f}^{p}$ & $\omega_{f}^{e}$&  & \\
\midrule
  Sod & type1 & 0.1 & 0.05 & 0.01 & 5.0 & 5.0 & 5.0 & 5.0 & 5.0 \\ 
  Sod & type2 & 0.1 & 0.05 & 0.01 &  & 5.0 & 5.0 & 5.0 & 5.0 \\
  Sod & type3 & 0.1 & 0.05 & 0.01 &  &  & 5.0 & 5.0 & 5.0 \\
  Lax & type1 & 1.0 & 0.5 & 0.1 & 100.0 & 100.0 & 10.0 & 100.0 & 100.0 \\ 
  Lax & type2 & 1.0 & 0.5 & 0.1 &  & 100.0 & 10.0 & 100.0 & 100.0 \\
  Lax & type3 & 1.0 & 0.5 & 0.1 &  &  & 10.0 & 100.0 & 100.0 \\
\bottomrule
\end{tabular}
\caption{\textit{Euler equations}: Loss weights settings.}
\label{tab: w settings euler}
\end{table}
 Networks' configurations are presented in \Cref{tab: Networks Config} and weights settings are shown in \Cref{tab: w settings euler}. With these settings, we train RelaxNN for $600,000$ epochs and one step for every epoch.  For the PINN framework, settings are the same without the extra neural network $\bm{v}_{\bm{\theta_{2}}}^{\text{NN}}$. In ~\Cref{fig:sod euler total loss}, We compare the total loss between the PINN and RelaxNN. In ~\Cref{fig: slice swe dam type1}, we compare the prediction of type1 RelaxNN and PINN at the final epoch for some specific times. Also for type2 RelaxNN in ~\Cref{fig: slice swe dam type2}.

\begin{figure}[htbp]
\centering
\subfloat[$t=0.0$]{\includegraphics[width = 0.4\textwidth]{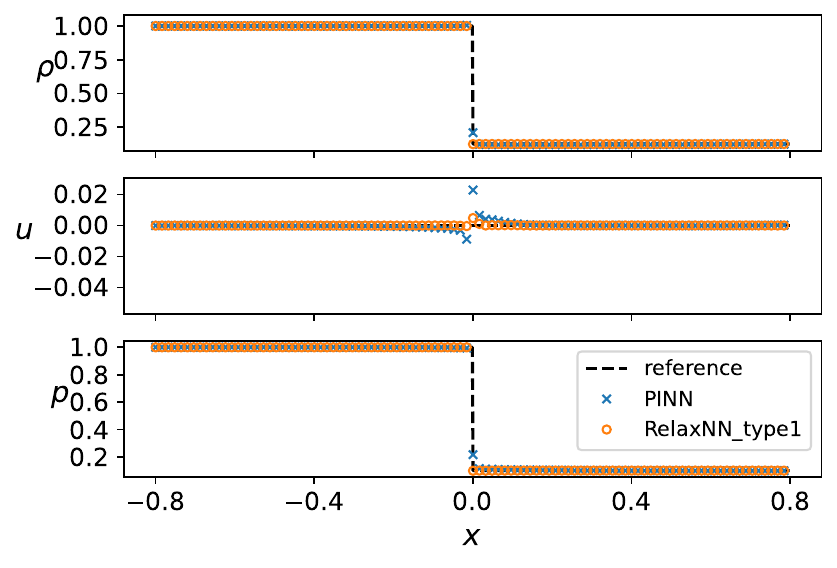}}
\subfloat[$t=0.08$]{\includegraphics[width = 0.4\textwidth]{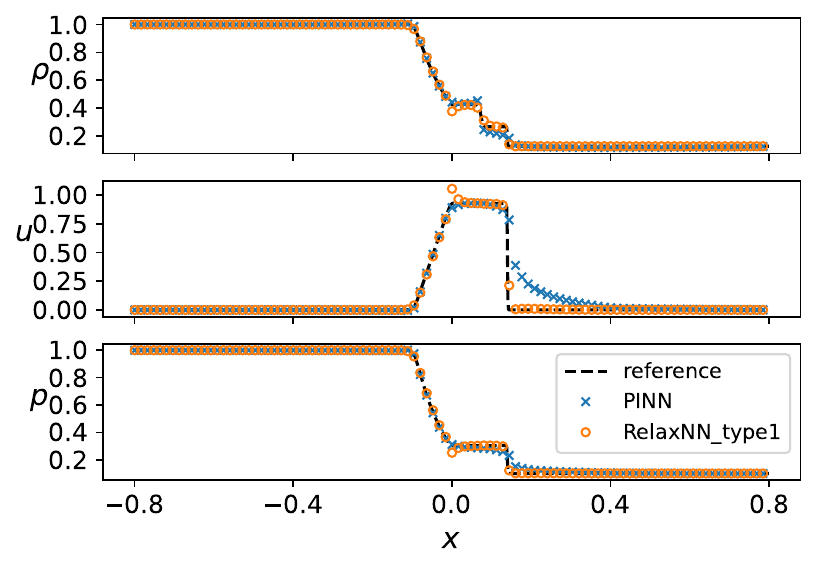}}
\\
\subfloat[$t=0.16$]{\includegraphics[width = 0.4\textwidth]{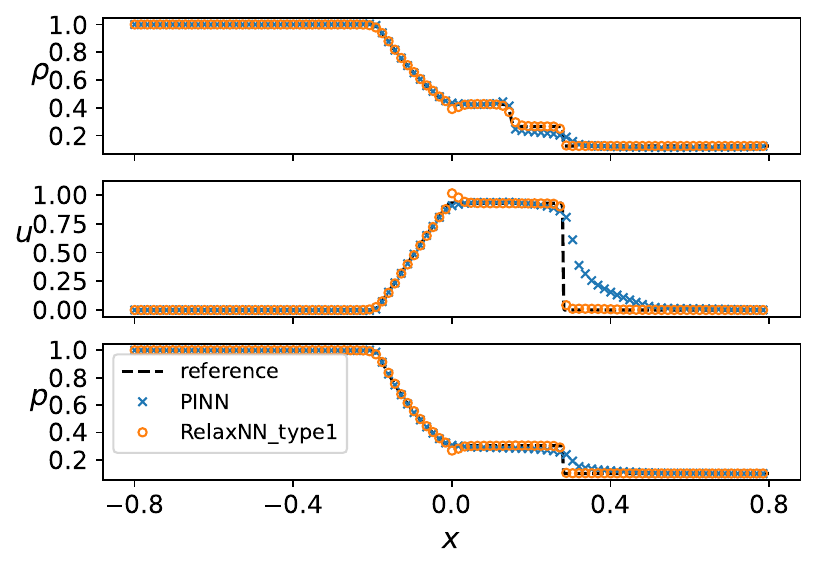}}
\subfloat[$t=0.24$]{\includegraphics[width = 0.4\textwidth]{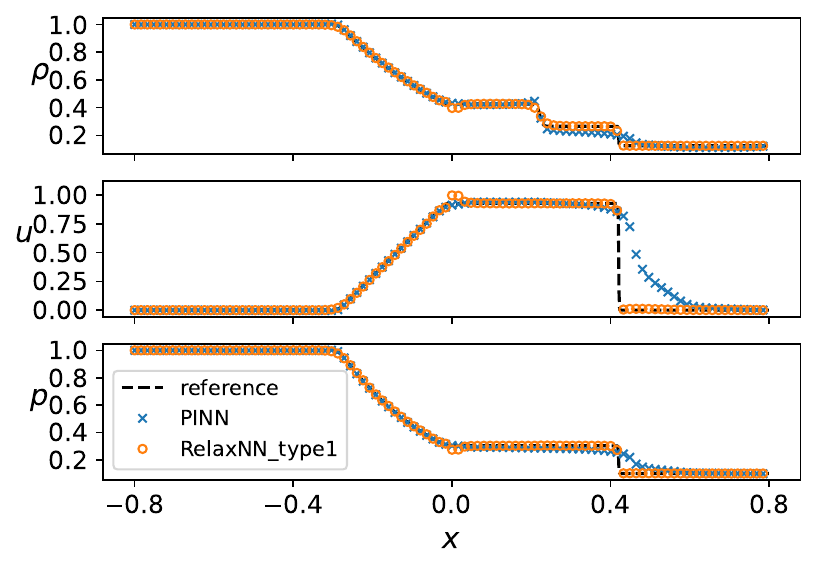}}
\\
\subfloat[$t=0.32$]{\includegraphics[width = 0.4\textwidth]{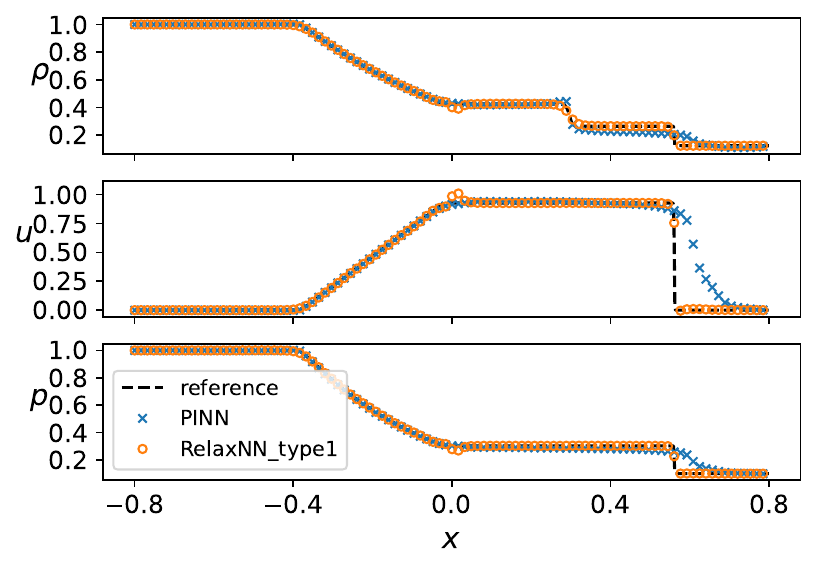}}
\subfloat[$t=0.40$]{\includegraphics[width = 0.4\textwidth]{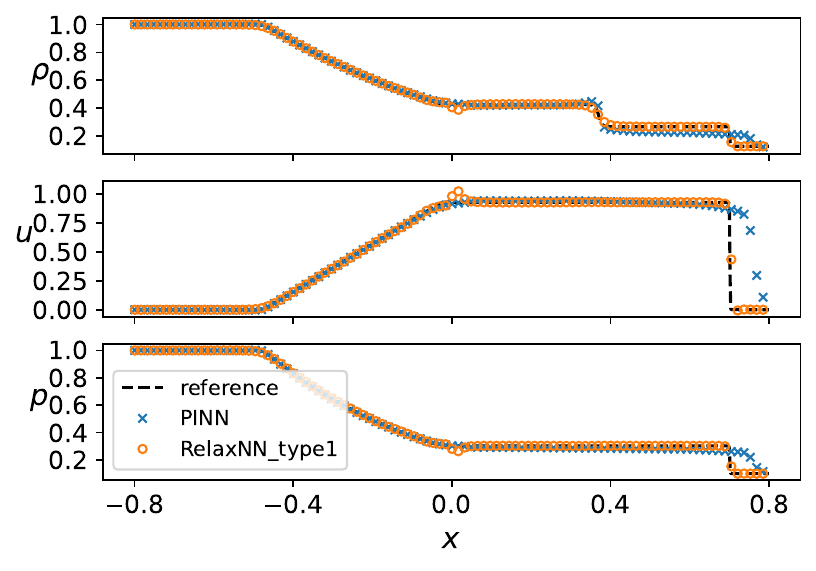}}
\caption{\textit{Euler equations with sod shock tube initial condition}: Comparison among the final epoch prediction of RelaxNN(type1), PINN, and the reference solution spatially at different specific moments. The resulting relative $L^{2}$ error of PINN, RelaxNN are $1.43 \times 10^{-2}$, $7.61 \times 10^{-4}$. For RelaxNN, the configuration of $\bm{u}_{\bm{\theta_{1}}}^{\text{NN}}$, $\bm{v}_{\bm{\theta_{2}}}^{\text{NN}}$ are [2,384,384,384,384,384,384,3] and [2,384,384,384,384,384,384,3]. We training for 600,000 epochs and one step per epoch. Loss weights settings are shown in ~\Cref{tab: w settings euler}. For PINN, all settings are the same as RelaxNN without the extra neural network.}
\label{fig: slice euler sod type1}
\end{figure}

\begin{figure}[htbp]
\centering
\subfloat[$t=0.0$]{\includegraphics[width = 0.4\textwidth]{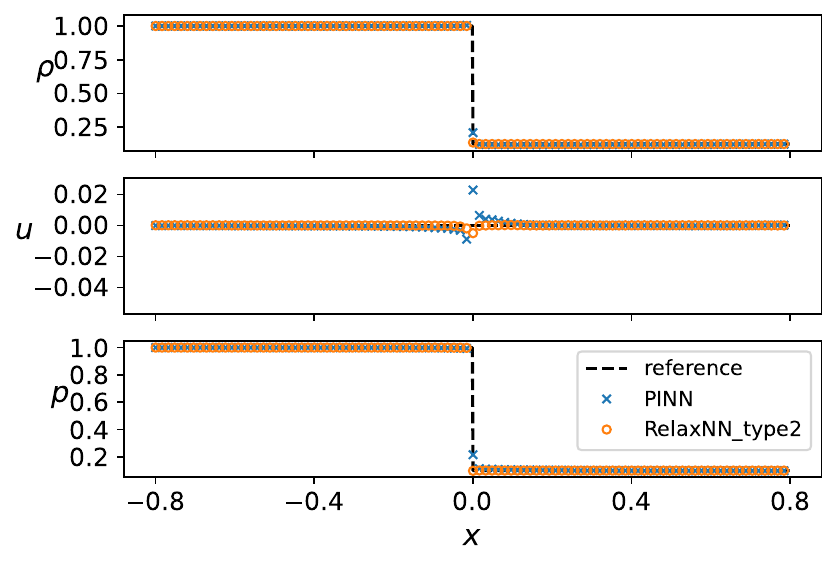}}
\subfloat[$t=0.08$]{\includegraphics[width = 0.4\textwidth]{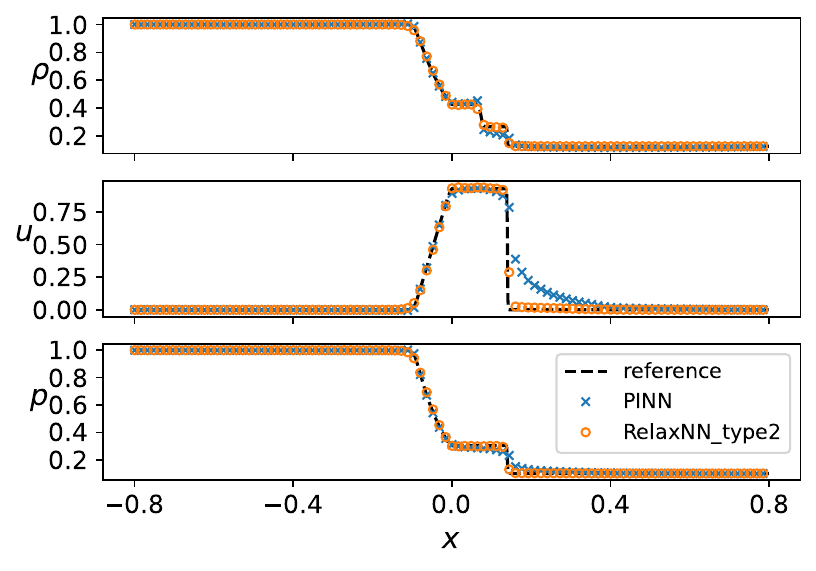}}
\\
\subfloat[$t=0.16$]{\includegraphics[width = 0.4\textwidth]{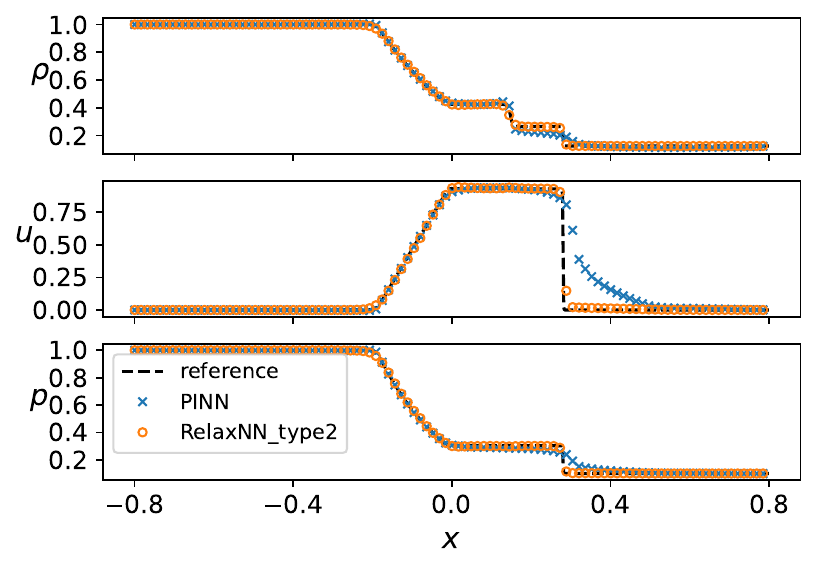}}
\subfloat[$t=0.24$]{\includegraphics[width = 0.4\textwidth]{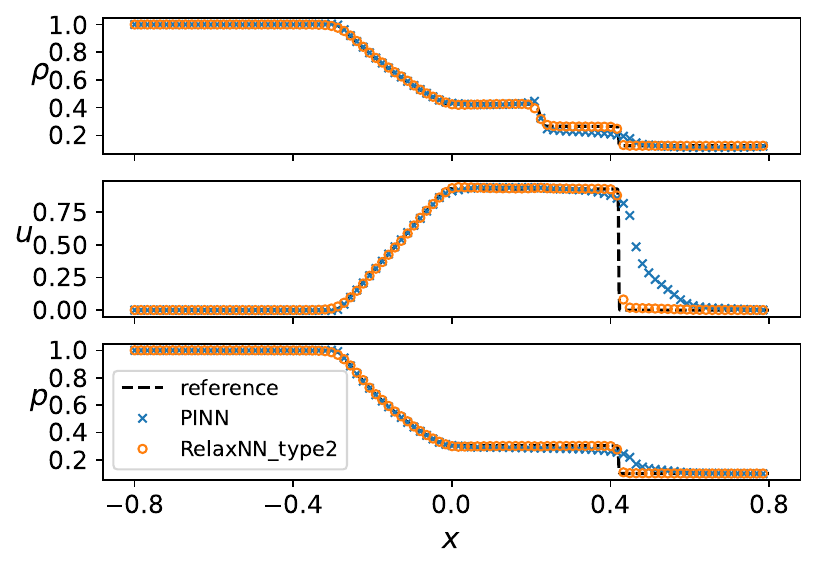}}
\\
\subfloat[$t=0.32$]{\includegraphics[width = 0.4\textwidth]{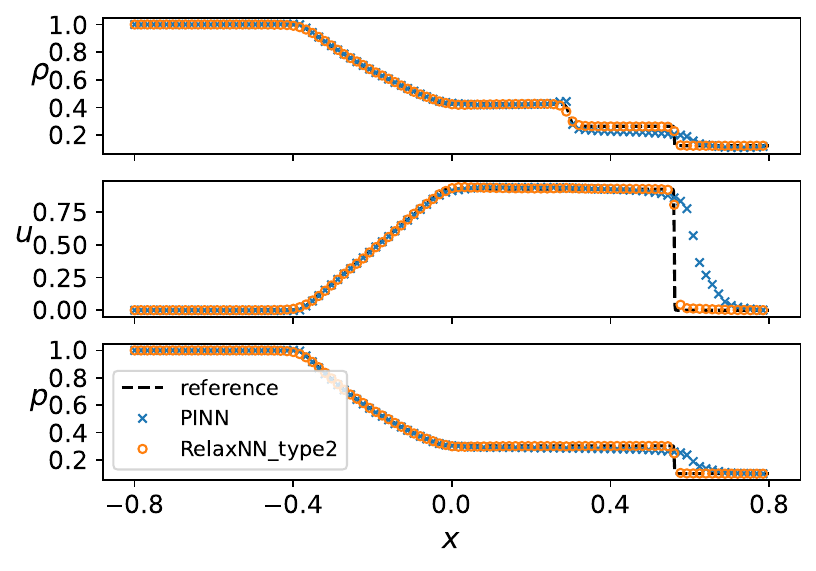}}
\subfloat[$t=0.40$]{\includegraphics[width = 0.4\textwidth]{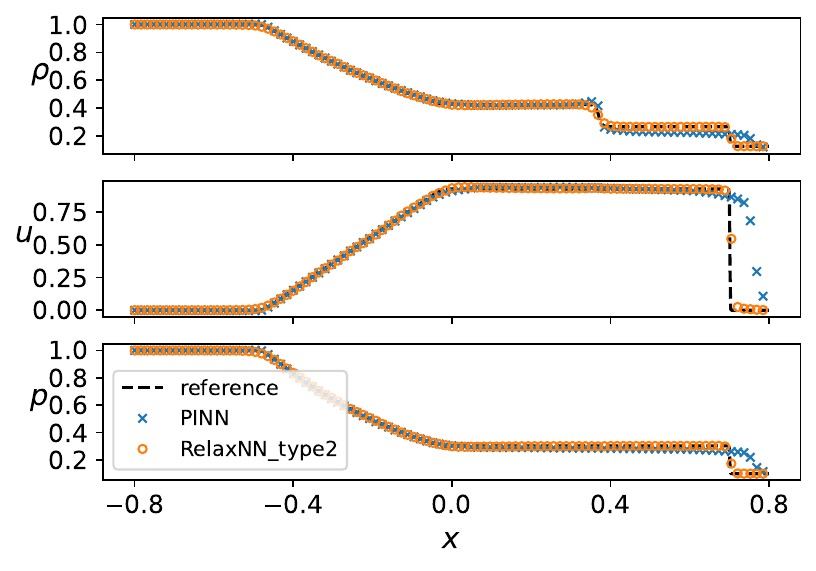}}
\caption{\textit{Euler equations with sod shock tube initial condition}: Comparison among the final epoch prediction of RelaxNN(type2), PINN, and the reference solution spatially at different specific moments. The resulting relative $L^{2}$ error of PINN, RelaxNN are $1.43 \times 10^{-2}$, $9.76 \times 10^{-4}$. For RelaxNN, the configuration of $\bm{u}_{\bm{\theta_{1}}}^{\text{NN}}$, $\bm{v}_{\bm{\theta_{2}}}^{\text{NN}}$ are [2,384,384,384,384,384,384,3] and [2,256,256,256,256,256,256,2]. We training for 600,000 epochs and one step per epoch. Loss weights settings are shown in ~\Cref{tab: w settings euler}. For PINN, all settings are the same as RelaxNN without the extra neural network.}
\label{fig: slice euler sod type2}
\end{figure}

\begin{figure}[htbp]
\centering
\subfloat[$t=0.0$]{\includegraphics[width = 0.4\textwidth]{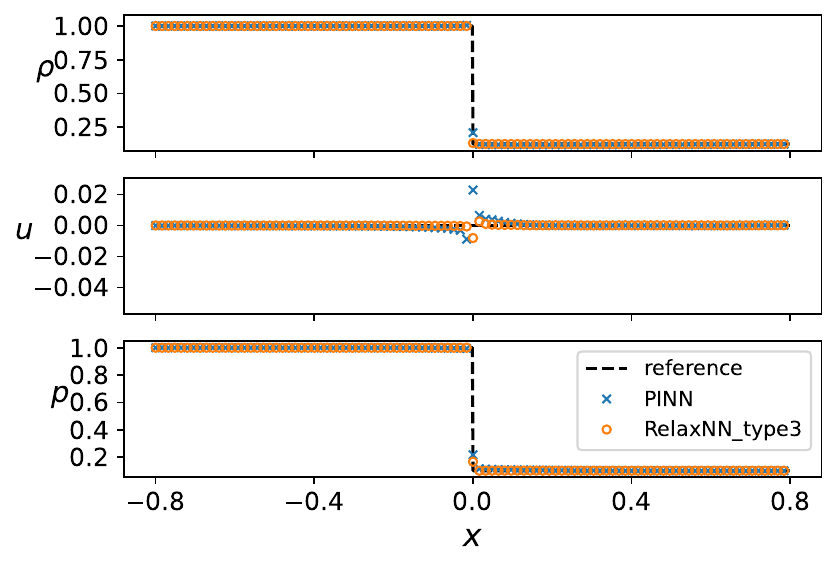}}
\subfloat[$t=0.08$]{\includegraphics[width = 0.4\textwidth]{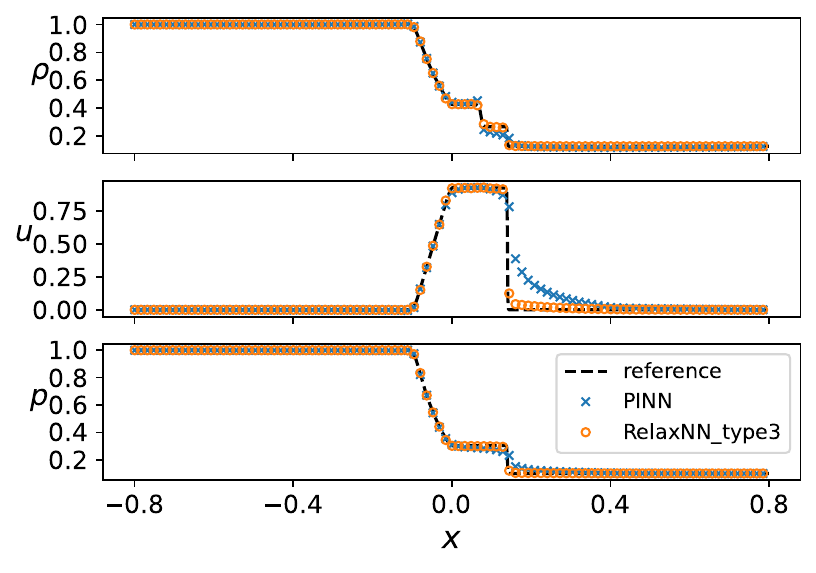}}
\\
\subfloat[$t=0.16$]{\includegraphics[width = 0.4\textwidth]{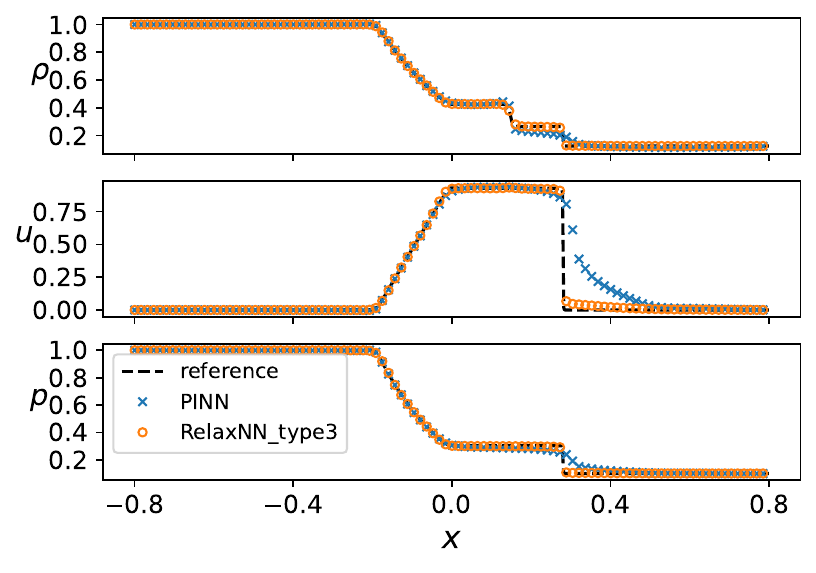}}
\subfloat[$t=0.24$]{\includegraphics[width = 0.4\textwidth]{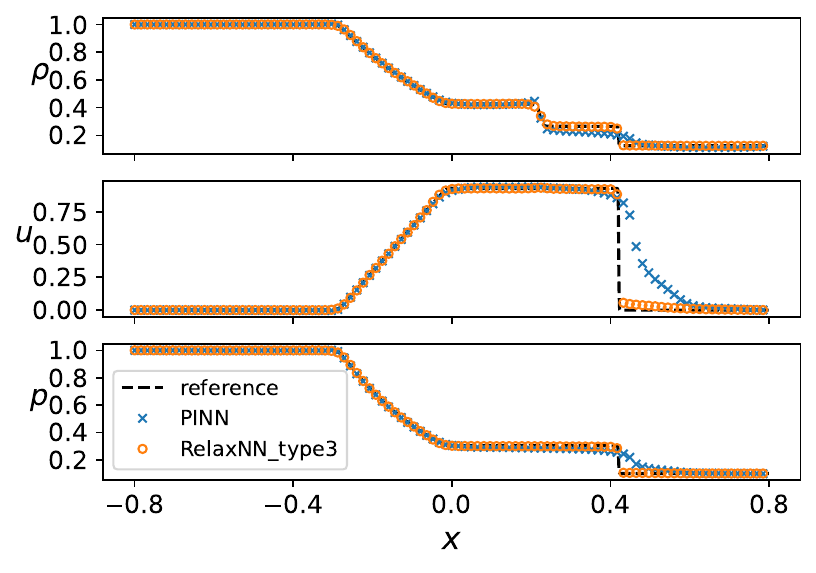}}
\\
\subfloat[$t=0.32$]{\includegraphics[width = 0.4\textwidth]{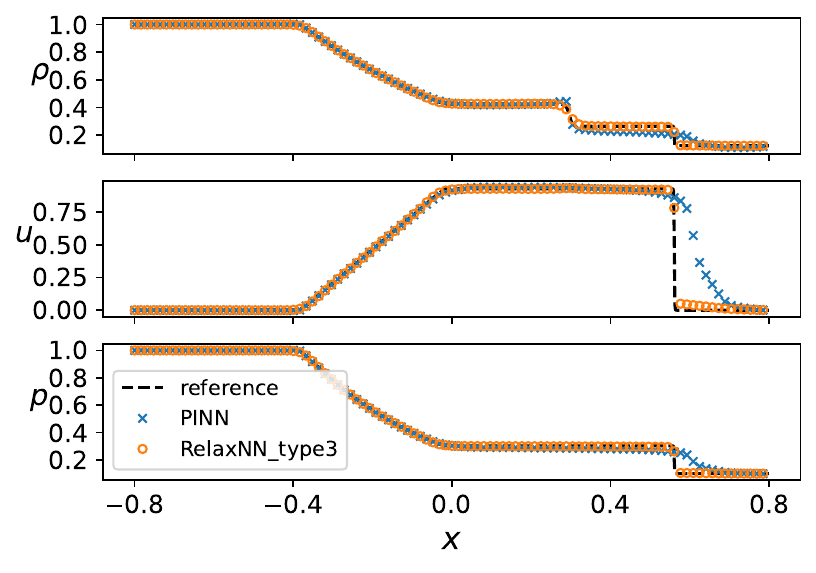}}
\subfloat[$t=0.40$]{\includegraphics[width = 0.4\textwidth]{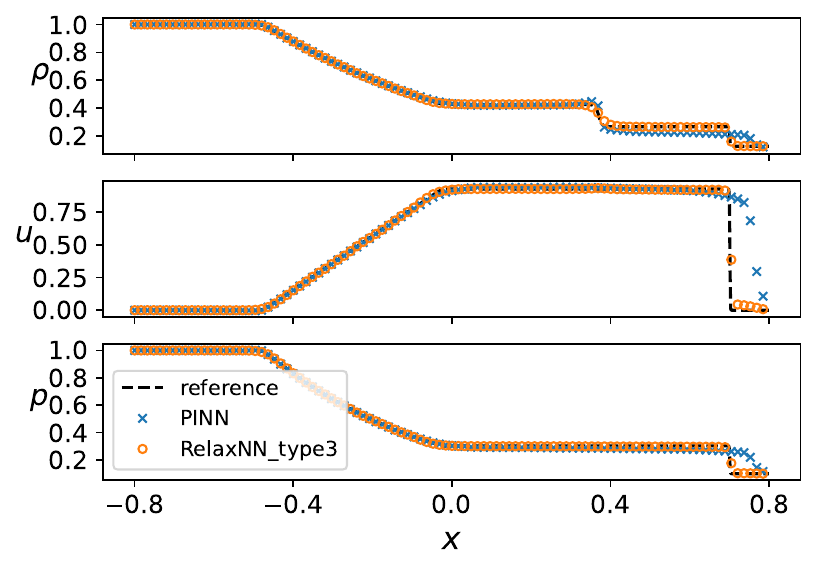}}
\caption{\textit{Euler equations with sod shock tube initial condition}: Comparison among the final epoch prediction of RelaxNN(type3), PINN, and the reference solution spatially at different specific moments. The resulting relative $L^{2}$ error of PINN, RelaxNN are $1.43 \times 10^{-2}$, $5.91 \times 10^{-4}$. For RelaxNN, the configuration of $\bm{u}_{\bm{\theta_{1}}}^{\text{NN}}$, $\bm{v}_{\bm{\theta_{2}}}^{\text{NN}}$ are [2,384,384,384,384,384,384,3] and [2,128,128,128,128,128,128,1]. We training for 600,000 epochs and one step per epoch. Loss weights settings are shown in ~\Cref{tab: w settings euler}. For PINN, all settings are the same as RelaxNN without the extra neural network.}
\label{fig: slice euler sod type3}
\end{figure}

\subsubsection{Lax problem}
\begin{figure}[htbp]
\centering
\subfloat[$t=0.0$]{\includegraphics[width = 0.4\textwidth]{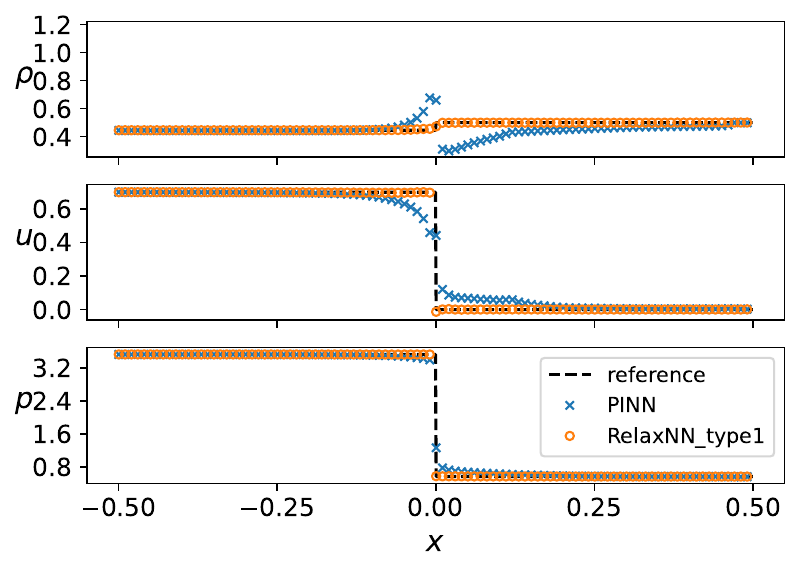}}
\subfloat[$t=0.032$]{\includegraphics[width = 0.4\textwidth]{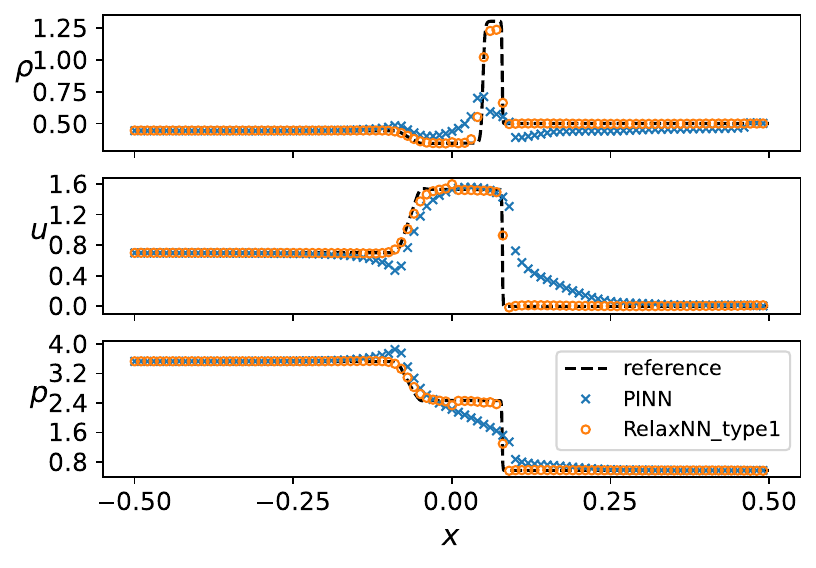}}
\\
\subfloat[$t=0.064$]{\includegraphics[width = 0.4\textwidth]{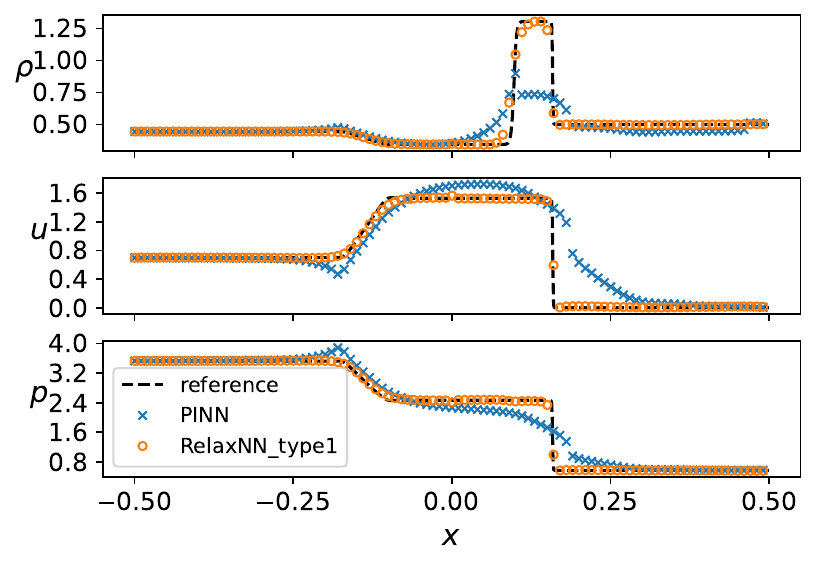}}
\subfloat[$t=0.096$]{\includegraphics[width = 0.4\textwidth]{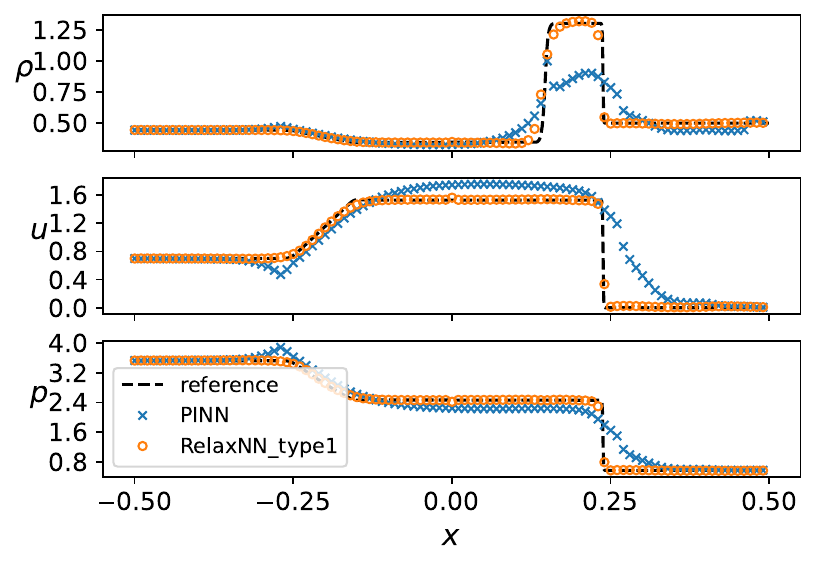}}
\\
\subfloat[$t=0.128$]{\includegraphics[width = 0.4\textwidth]{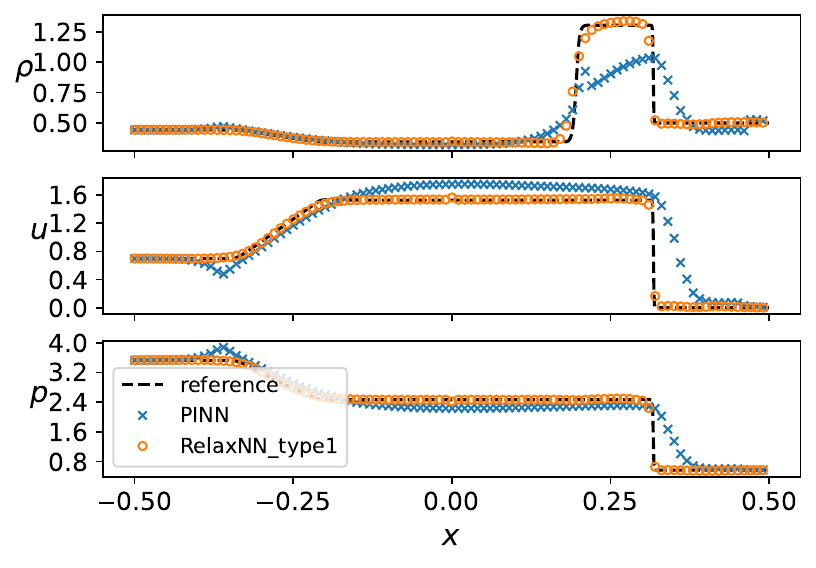}}
\subfloat[$t=0.16$]{\includegraphics[width = 0.4\textwidth]{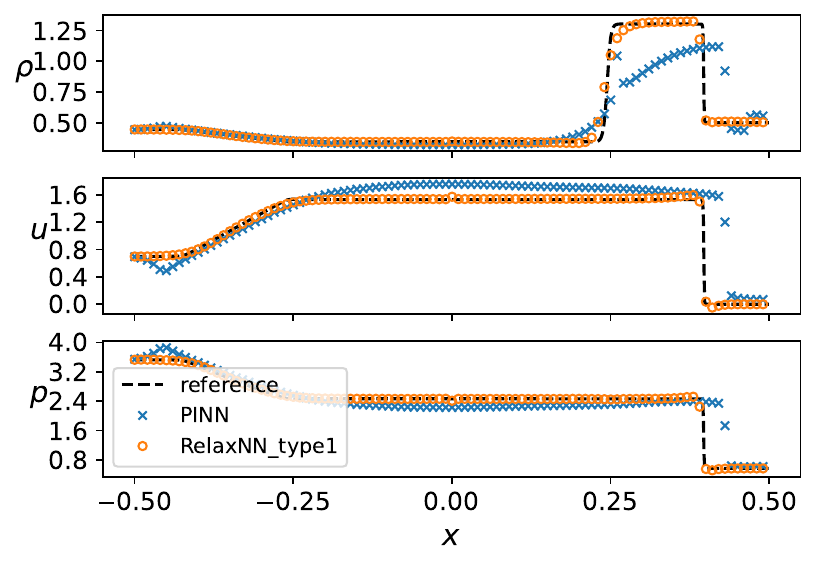}}
\caption{\textit{Euler equations with lax shock tube initial condition}: Comparison among the final epoch prediction of RelaxNN(type1), PINN, and the reference solution spatially at different specific moments. The resulting relative $L^{2}$ error of PINN, RelaxNN are $1.99 \times 10^{-2}$, $6.94 \times 10^{-4}$. For RelaxNN, the configuration of $\bm{u}_{\bm{\theta_{1}}}^{\text{NN}}$, $\bm{v}_{\bm{\theta_{2}}}^{\text{NN}}$ are [2,384,384,384,384,384,384,3] and [2,384,384,384,384,384,384,3]. We training for 600,000 epochs and one step per epoch. Loss weights settings are shown in ~\Cref{tab: w settings euler}. For PINN, all settings are the same as RelaxNN without the extra neural network.}
\label{fig: slice euler lax type1}
\end{figure}

\begin{figure}[htbp]
\centering
\subfloat[$t=0.0$]{\includegraphics[width = 0.4\textwidth]{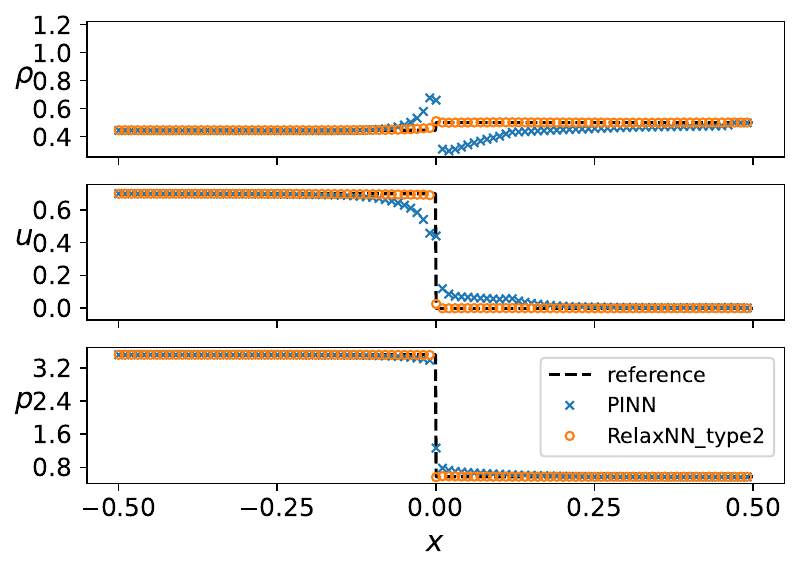}}
\subfloat[$t=0.032$]{\includegraphics[width = 0.4\textwidth]{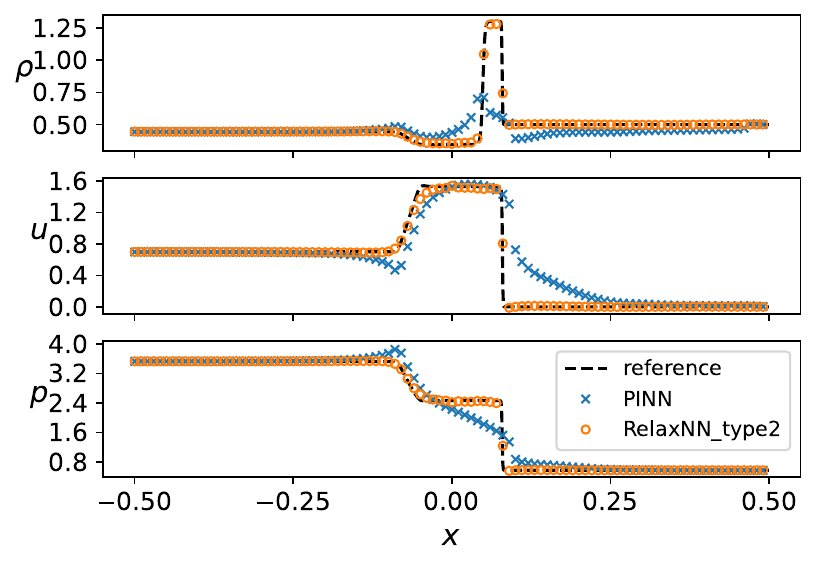}}
\\
\subfloat[$t=0.064$]{\includegraphics[width = 0.4\textwidth]{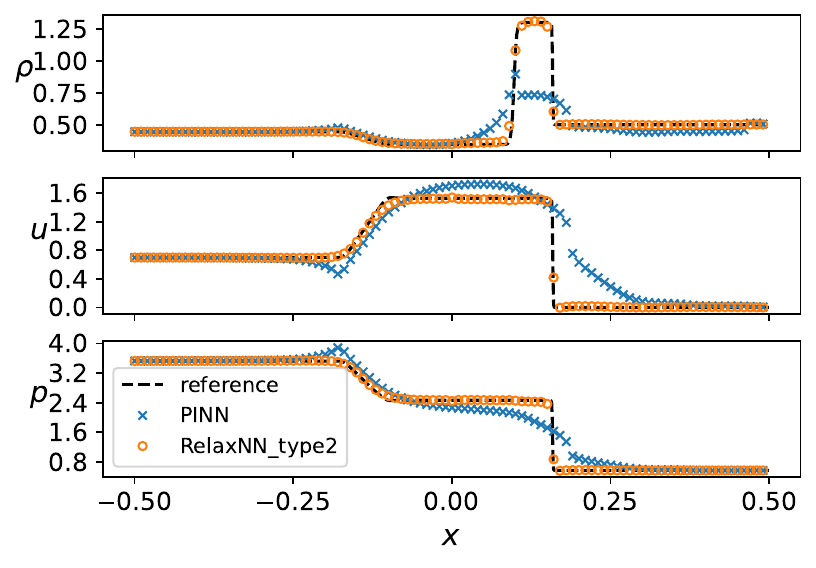}}
\subfloat[$t=0.096$]{\includegraphics[width = 0.4\textwidth]{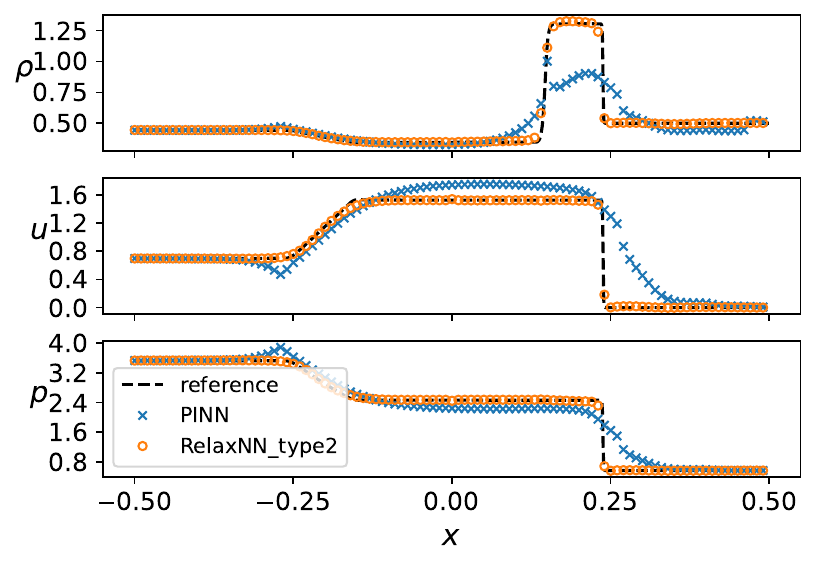}}
\\
\subfloat[$t=0.128$]{\includegraphics[width = 0.4\textwidth]{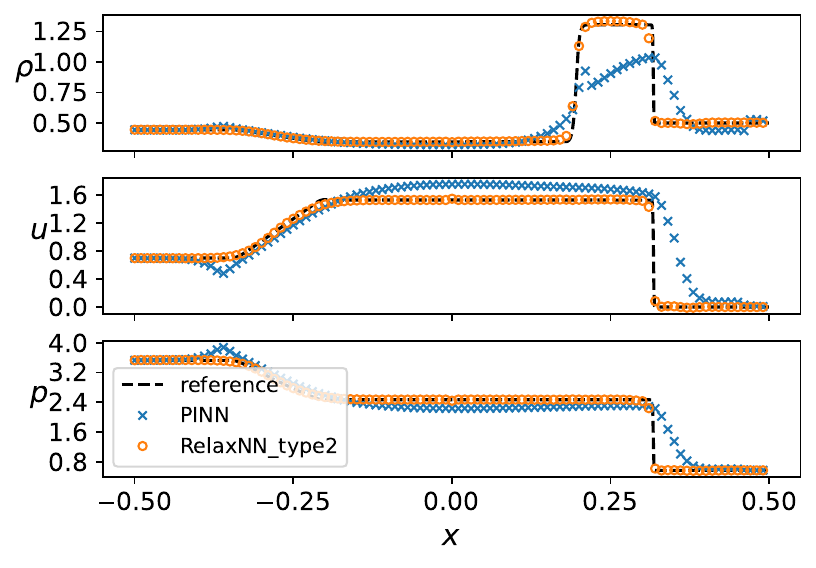}}
\subfloat[$t=0.16$]{\includegraphics[width = 0.4\textwidth]{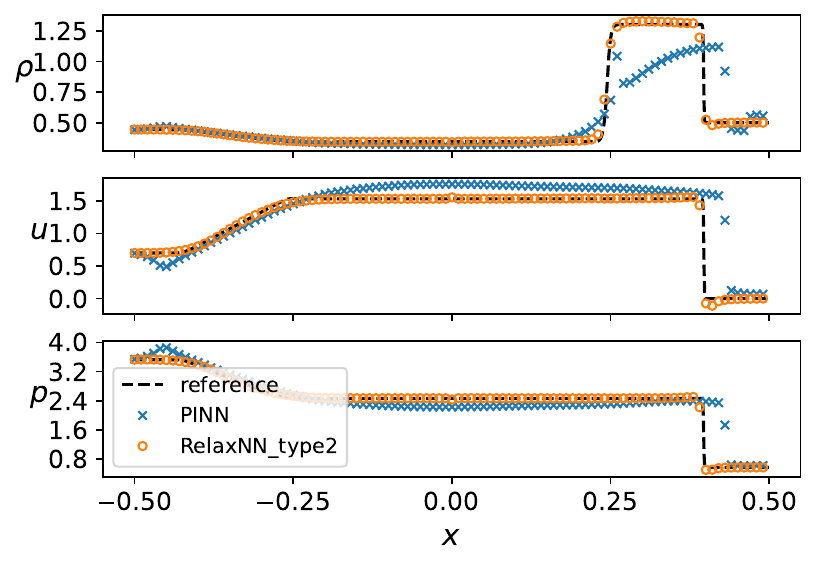}}
\caption{\textit{Euler equations with lax shock tube initial condition}: Comparison among the final epoch prediction of RelaxNN(type1), PINN, and the reference solution spatially at different specific moments. The resulting relative $L^{2}$ error of PINN, RelaxNN are $1.99 \times 10^{-2}$, $7.34 \times 10^{-4}$. For RelaxNN, the configuration of $\bm{u}_{\bm{\theta_{1}}}^{\text{NN}}$, $\bm{v}_{\bm{\theta_{2}}}^{\text{NN}}$ are [2,384,384,384,384,384,384,3] and [2,256,256,256,256,256,256,2]. We training for 600,000 epochs and one step per epoch. Loss weights settings are shown in ~\Cref{tab: w settings euler}. For PINN, all settings are the same as RelaxNN without the extra neural network.}
\label{fig: slice euler lax type2}
\end{figure}

\begin{figure}[htbp]
\centering
\subfloat[$t=0.0$]{\includegraphics[width = 0.4\textwidth]{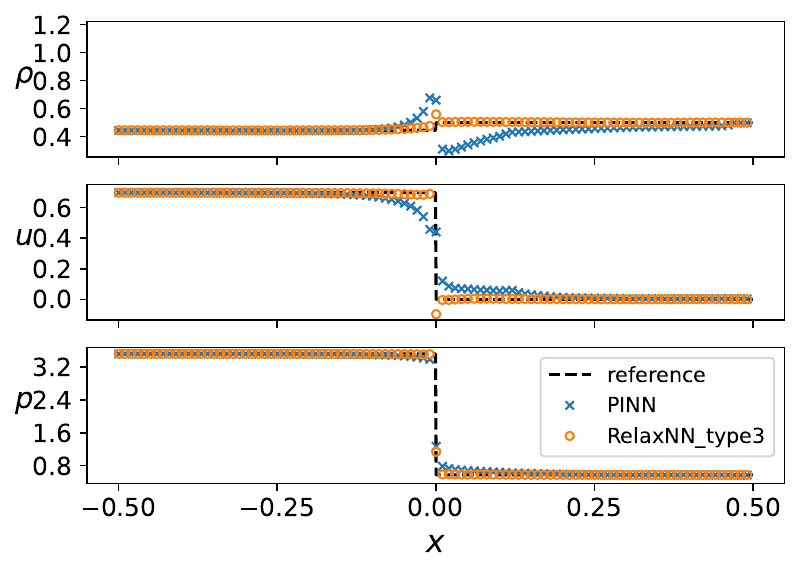}}
\subfloat[$t=0.032$]{\includegraphics[width = 0.4\textwidth]{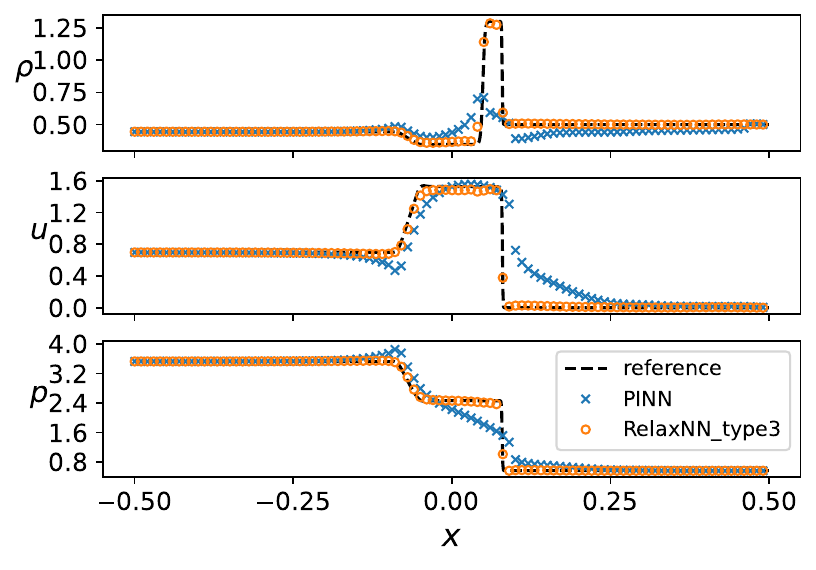}}
\\
\subfloat[$t=0.064$]{\includegraphics[width = 0.4\textwidth]{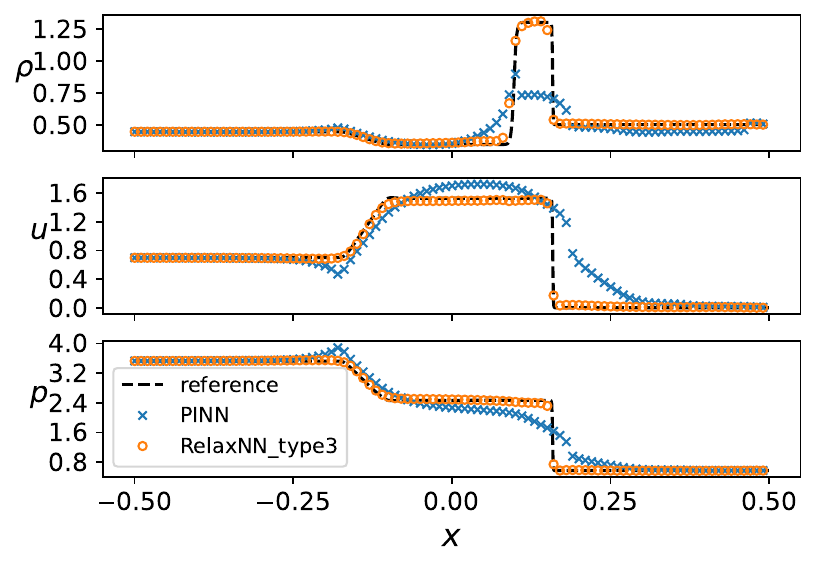}}
\subfloat[$t=0.096$]{\includegraphics[width = 0.4\textwidth]{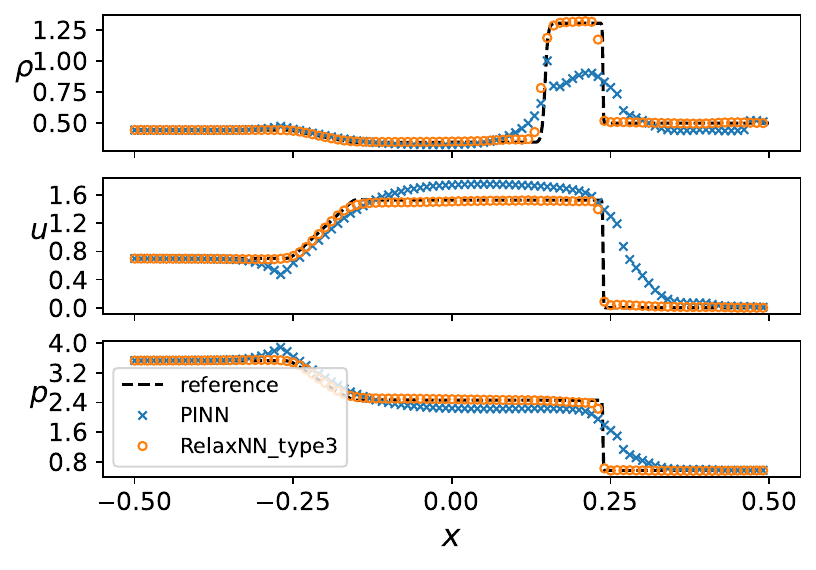}}
\\
\subfloat[$t=0.128$]{\includegraphics[width = 0.4\textwidth]{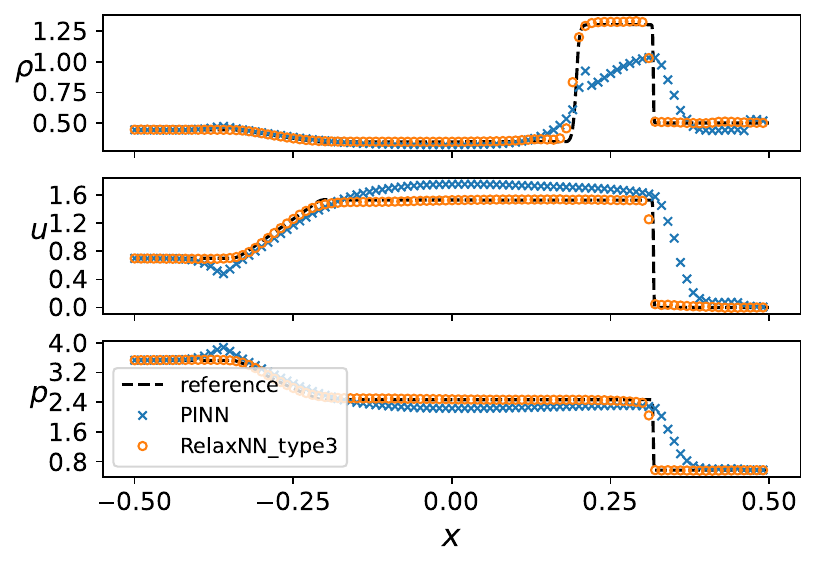}}
\subfloat[$t=0.16$]{\includegraphics[width = 0.4\textwidth]{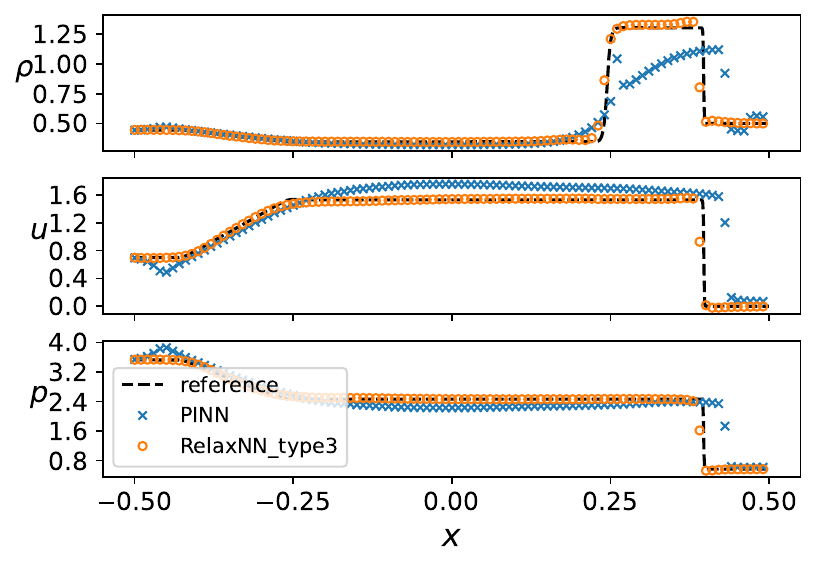}}
\caption{\textit{Euler equations with lax shock tube initial condition}: Comparison among the final epoch prediction of RelaxNN(type1), PINN, and the reference solution spatially at different specific moments. The resulting relative $L^{2}$ error of PINN, RelaxNN are $1.99 \times 10^{-2}$, $1.68 \times 10^{-3}$. For RelaxNN, the configuration of $\bm{u}_{\bm{\theta_{1}}}^{\text{NN}}$, $\bm{v}_{\bm{\theta_{2}}}^{\text{NN}}$ are [2,384,384,384,384,384,384,3] and [2,128,128,128,128,128,128,1]. We training for 600,000 epochs and one step per epoch. Loss weights settings are shown in ~\Cref{tab: w settings euler}. For PINN, all settings are the same as RelaxNN without the extra neural network.}
\label{fig: slice euler lax type3}
\end{figure}

We consider the euler equations in domain $ \Omega = \{ (t,x)\} = [ 0.0,0.16]\times [-0.5,0.5]$ with lax shocktube initial condition as
\begin{equation}
    \bm{u}_{0}(0,x) = \left( 
    \begin{array}{cc}
         0.445 \\ 
         0.698 \\
         3.528
    \end{array}
    \right) , \quad \text{if} \quad -0.5 \leq  x \leq 0.0 ;
    \quad
    \bm{u}_{0}(0,x) = \left( 
    \begin{array}{cc}
         0.5 \\
         0.0 \\
         0.571
    \end{array}
    \right) , \quad \text{if} \quad 0.0 < x \leq 0.5 ;
\end{equation}
 Networks' configurations are presented in \Cref{tab: Networks Config} and weights settings are shown in \Cref{tab: w settings euler}. With these settings, we train RelaxNN for $600,000$ epochs and one step for every epoch.  For the PINN framework, settings are the same without the extra neural network $\bm{v}_{\bm{\theta_{2}}}^{\text{NN}}$. In ~\Cref{fig:lax euler total loss}, We compare the total loss between the PINN and RelaxNN. In ~\Cref{fig: slice euler lax type1}, we compare the prediction of type1 RelaxNN and PINN at the final epoch for some specific times. Also for type2 RelaxNN in ~\Cref{fig: slice euler lax type2} and type3 RelaxNN in ~\Cref{fig: slice euler lax type3}.

\subsection{Uncertainty quantification problems}
In this part, we consider stochastic perturbation on the initial condition for the equations we discussed above, well known as uncertainty quantification problems. For simplicity, we do not write down the loss term or the loss weights again. All we need is to incorporate the random variables into the inputs of our neural networks. The outstanding performances reveal the potential of our method to solve the high-dimensional problems.

\subsubsection{UQ problems in inviscid Burgers' equation with Riemann initial condition}
\begin{figure}[htbp]
\centering
\subfloat[$t=0.1$]{\includegraphics[width = 0.4\textwidth]{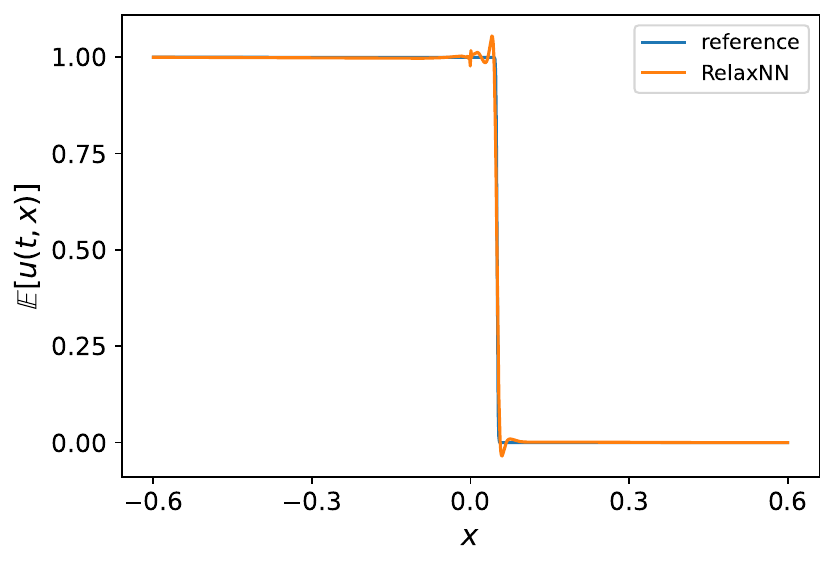}}
\subfloat[$t=0.1$]{\includegraphics[width = 0.4\textwidth]{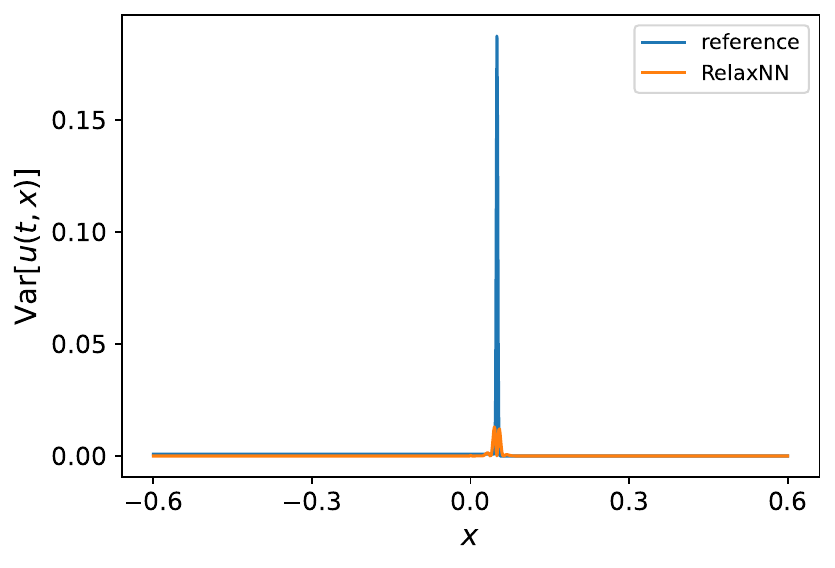}}
\\
\subfloat[$t=0.5$]{\includegraphics[width = 0.4\textwidth]{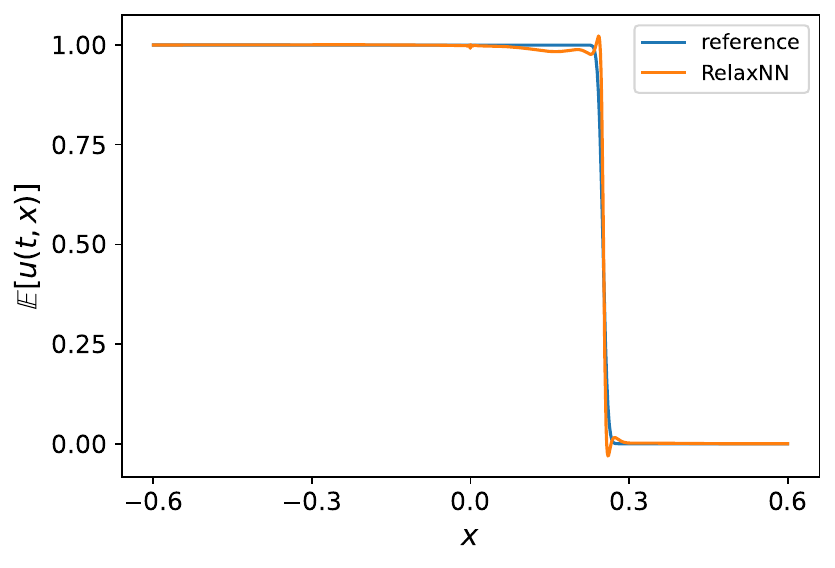}}
\subfloat[$t=0.5$]{\includegraphics[width = 0.4\textwidth]{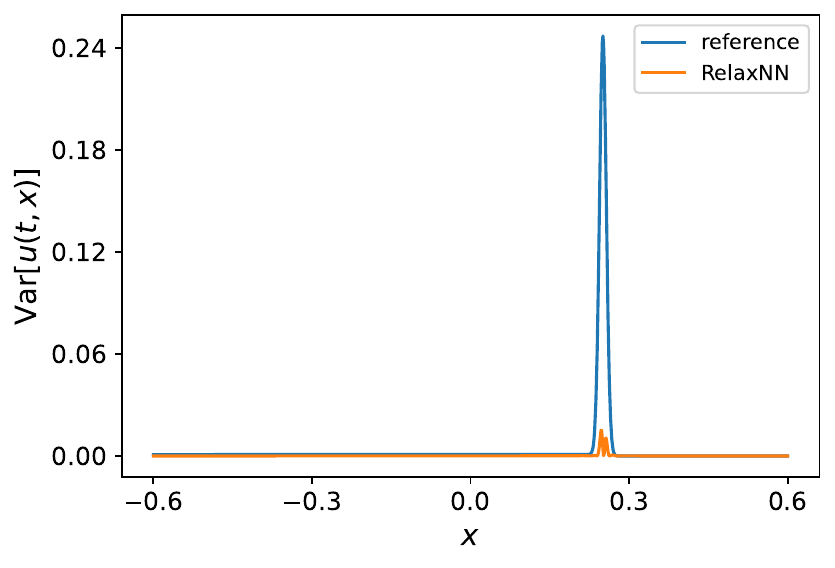}}
\\
\subfloat[$t=1.0$]{\includegraphics[width = 0.4\textwidth]{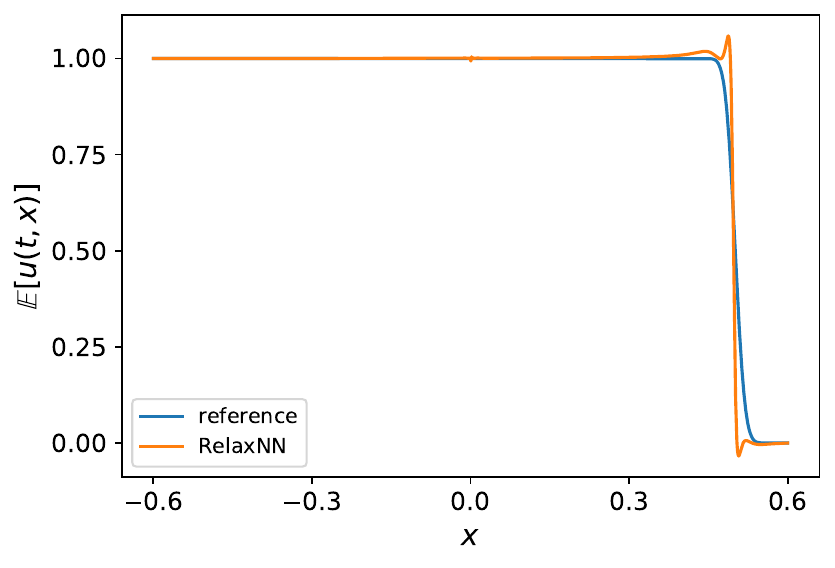}}
\subfloat[$t=1.0$]{\includegraphics[width = 0.4\textwidth]{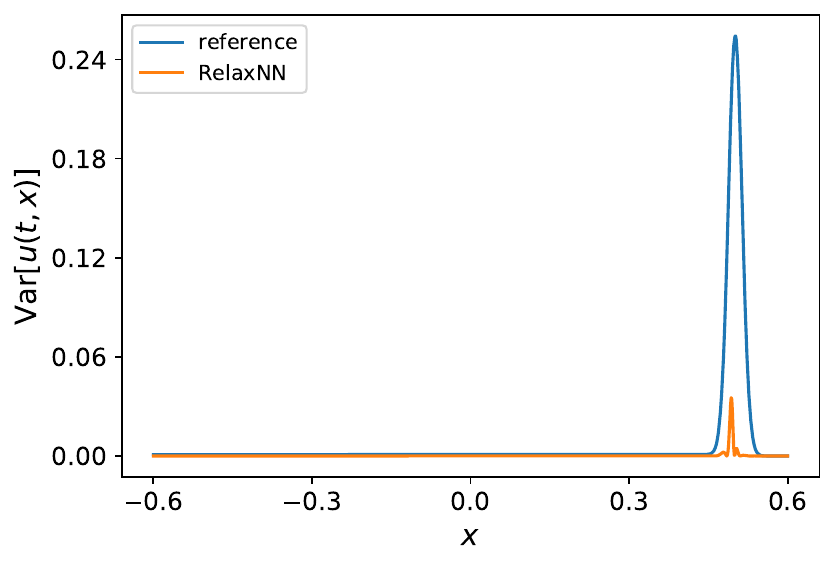}}
\caption{\textit{Burgers' equation with stochastic Riemann initial condition}: Comparison among the final epoch prediction of RelaxNN(type3), and the reference solution spatially at different specific moments. The relative $L^{2}$ error of expectation, variance are $2.40 \times 10^{-2}$, $3.42 \times 10^{-2}$. For RelaxNN, the configuration of $\bm{u}_{\bm{\theta_{1}}}^{\text{NN}}$, $\bm{v}_{\bm{\theta_{2}}}^{\text{NN}}$ are [102,128,128,128,1] and [102,64,64,64,1]. We training for 300,000 epochs and one step per epoch. }
\label{fig: slice stochastic burgers equation}
\end{figure}

We consider the inviscid Burgers' equation \eqref{eq: burgers equation} in domain $ \Omega = \{ (t,x)\} = [ 0.0,1.0]\times [-0.6,0.6]$
with stochastic Riemann initial condition as
\begin{equation}
    u_{0}(0,x,\mathbf{z}) = \left\{
    \begin{aligned}
        & 1.0 + \varepsilon \sum_{i=1}^{s} z_{i}  \, , \quad -0.6 \leq x \leq 0 \\ 
        & 0.0 , \quad 0 < x \leq 0.6
    \end{aligned}
    \right. ,
\end{equation}
where $\mathbf{z} = (z_{1}, z_{2}, \dots, z_{s}) \sim \mathcal{U}([-1,1]^{s}), \; \varepsilon = 0.005, \; s=100$. We then obtain the reference variance of the solution by the Monte Carlo method for $100,000$ random experiments. The mean and variance of the solution of our RelaxNN are also obtained by the Monte Carlo method for $1,000,000$ random experiments. Our results are shown in ~\Cref{fig: slice stochastic burgers equation}.

\subsubsection{UQ problems in shallow water equations with two shocks initial condition}
\begin{figure}[htbp]
\centering
\subfloat[$t=0.1$]{\includegraphics[width = 0.4\textwidth]{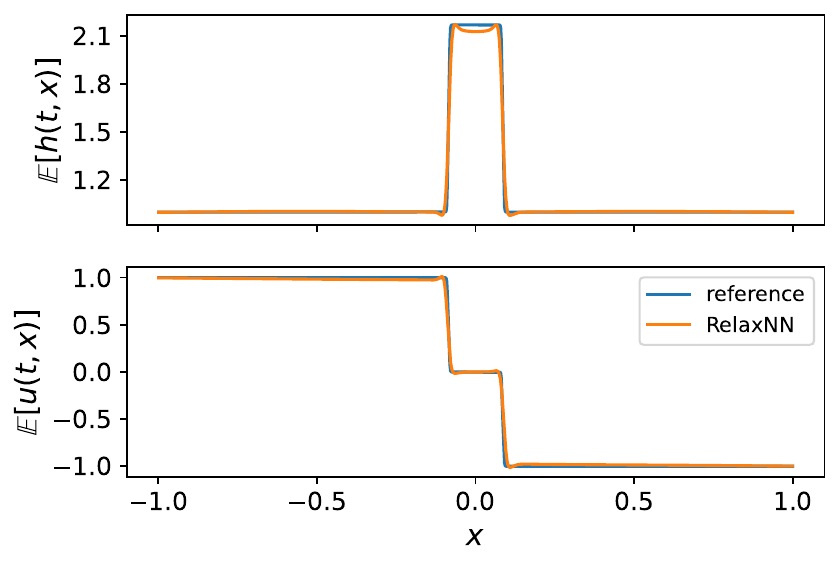}}
\subfloat[$t=0.1$]{\includegraphics[width = 0.4\textwidth]{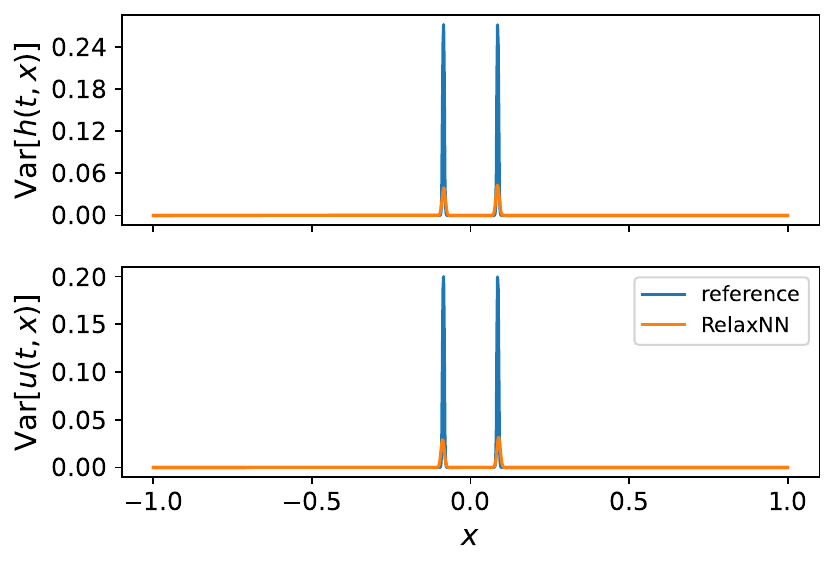}}
\\
\subfloat[$t=0.5$]{\includegraphics[width = 0.4\textwidth]{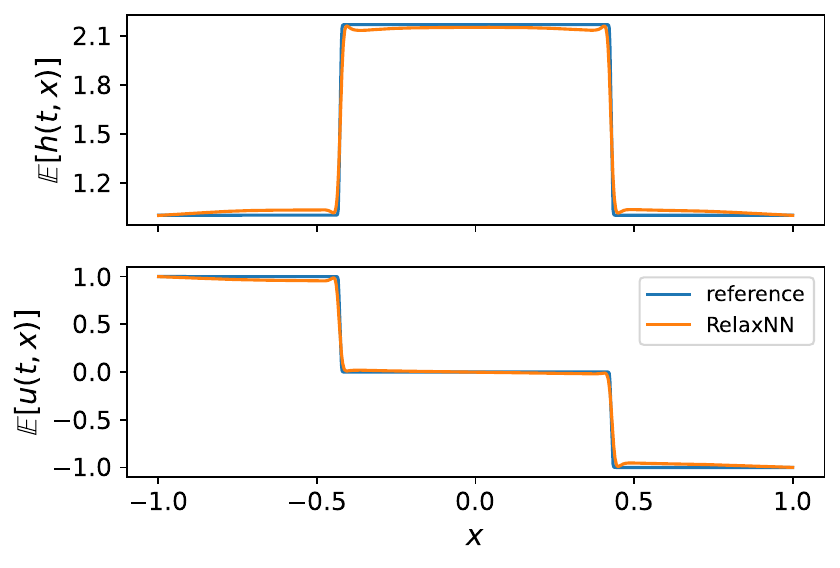}}
\subfloat[$t=0.5$]{\includegraphics[width = 0.4\textwidth]{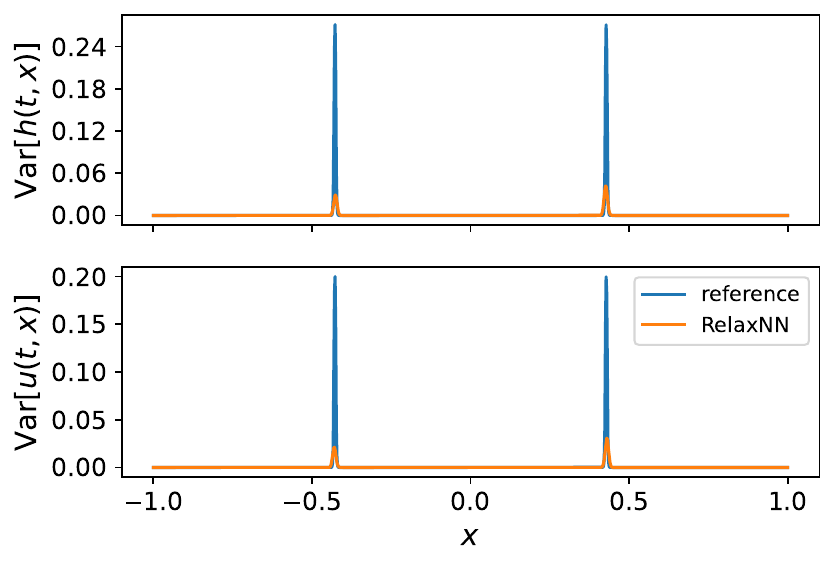}}
\\
\subfloat[$t=1.0$]{\includegraphics[width = 0.4\textwidth]{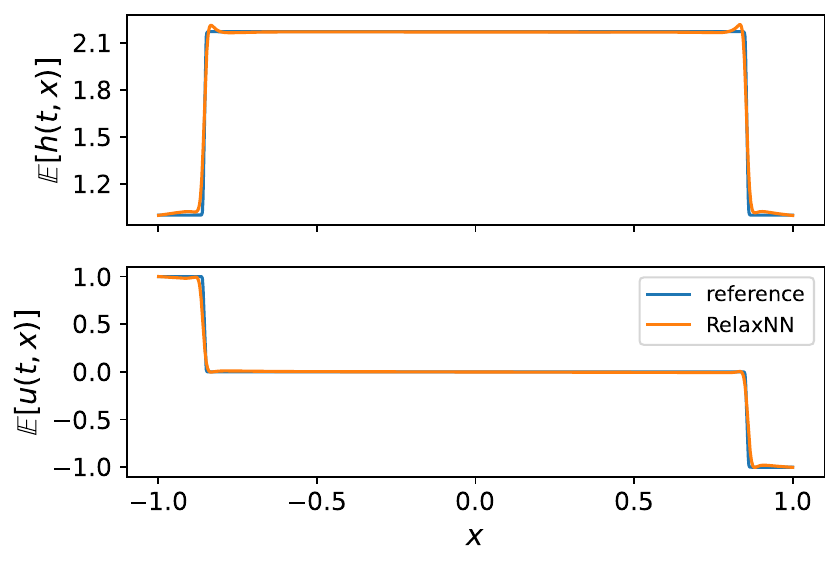}}
\subfloat[$t=1.0$]{\includegraphics[width = 0.4\textwidth]{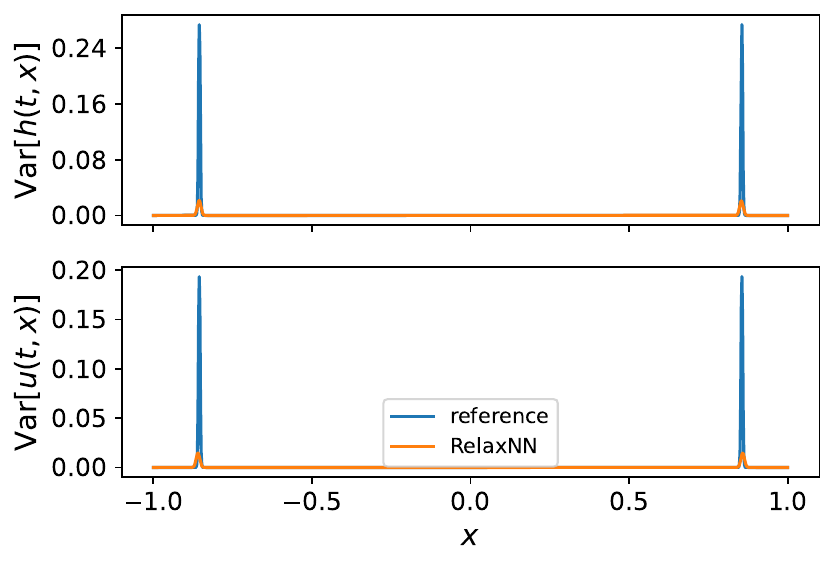}}
\caption{\textit{Shallow water equations with stochastic two shock initial conditions}: Comparison among the final epoch prediction of RelaxNN(type3) and the reference solution spatially at different specific moments. The resulting relative $L^{2}$ error of mean and variance are $2.24 \times 10^{-2}$, $2.60 \times 10^{-2}$. For RelaxNN, the configuration of $\bm{u}_{\bm{\theta_{1}}}^{\text{NN}}$, $\bm{v}_{\bm{\theta_{2}}}^{\text{NN}}$ are [7,256,256,256,256,256,256,2] and [7,128,128,128,128,128,128,1]. We training for 600,000 epochs and one step per epoch. Loss weights settings are shown in ~\Cref{tab: w settings euler}.}
\label{fig: slice stochastic shallow water eq}
\end{figure}

We consider the shallow water equations ~\eqref{eq: shallow water equations} in domain $ \Omega = \{ (t,x)\} = [ 0.0, 1.0]\times [-1.0,1.0]$ with stochastic two shock initial conditions as
\begin{equation}
    \bm{u}_{0}(0,x) = \left( 
    \begin{array}{cc}
         1.0 \\ 
         1.0
    \end{array}
    \right) , \quad \text{if} \quad -1.0 \leq  x \leq \psi(\mathbf{z}) ;
    \quad
    \bm{u}_{0}(0,x) = \left( 
    \begin{array}{cc}
         1.0 \\
         -1.0 
    \end{array}
    \right) ,  \quad \text{if} \quad \psi(\mathbf{z}) < x \leq 1.0 ;
\end{equation}
Where $\psi(\mathbf{z}) = \varepsilon(z_{1}\sigma(z_{2}z_{3} + z_{4}) +z_{5})$, $\sigma(\cdot)$ is the ReLU function and $\varepsilon = 0.005$.We obtain the reference mean and variance of the solution by the Monte Carlo method for $100,000$ random experiments. The mean and variance of the solution of our RelaxNN are obtained by numerical quadrature. Specifically, we used the Gaussian-legendre quadrature with $10$ points in every dimension. In ~\Cref{fig: slice stochastic shallow water eq}, we display our results.

\subsubsection{UQ problems in Euler equations with Sod initial condition}
\begin{figure}[htbp]
\centering
\subfloat[$t=0.04$]{\includegraphics[width = 0.4\textwidth]{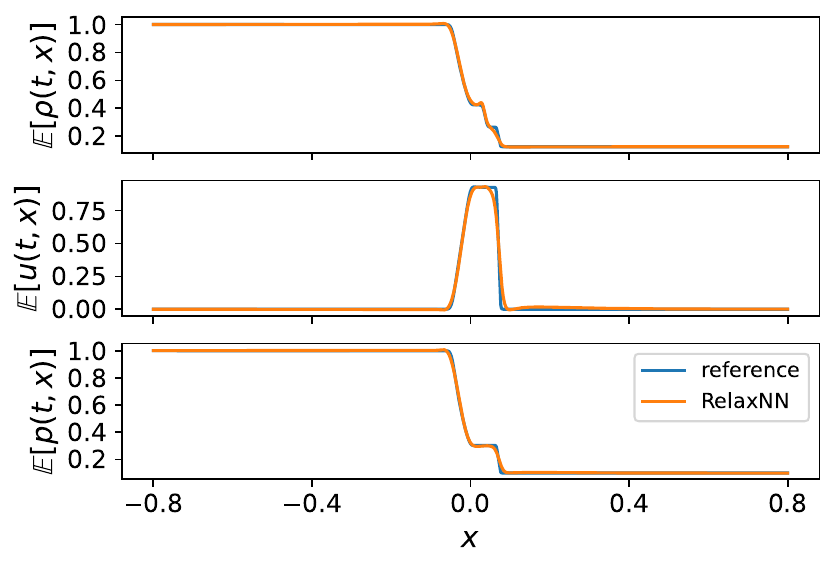}}
\subfloat[$t=0.04$]{\includegraphics[width = 0.4\textwidth]{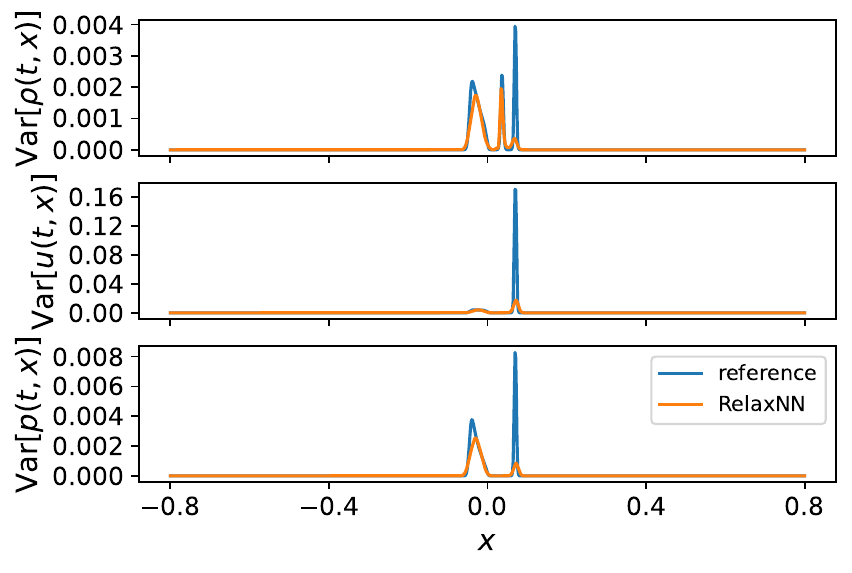}}
\\
\subfloat[$t=0.20$]{\includegraphics[width = 0.4\textwidth]{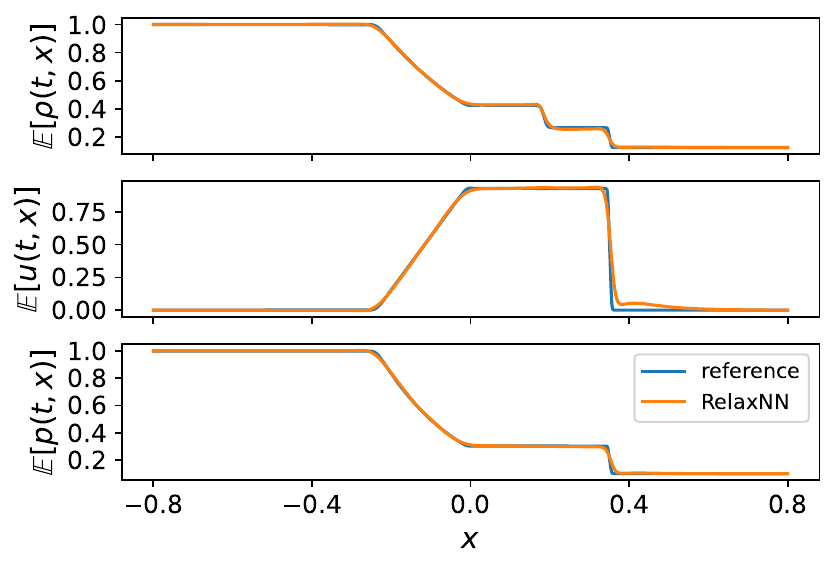}}
\subfloat[$t=0.20$]{\includegraphics[width = 0.4\textwidth]{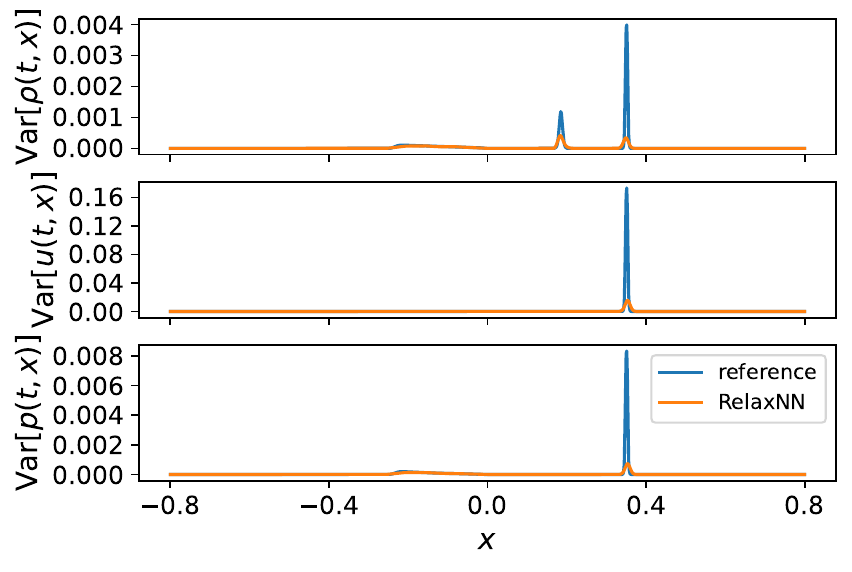}}
\\
\subfloat[$t=0.40$]{\includegraphics[width = 0.4\textwidth]{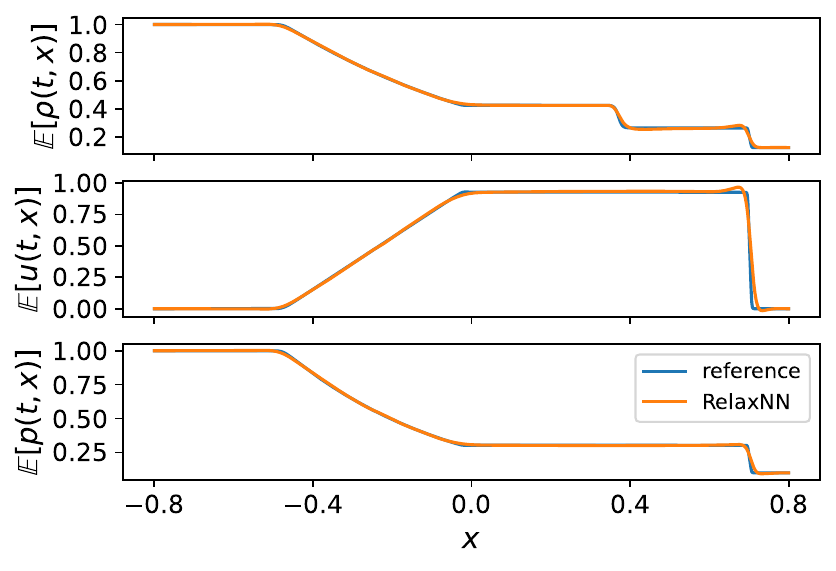}}
\subfloat[$t=0.40$]{\includegraphics[width = 0.4\textwidth]{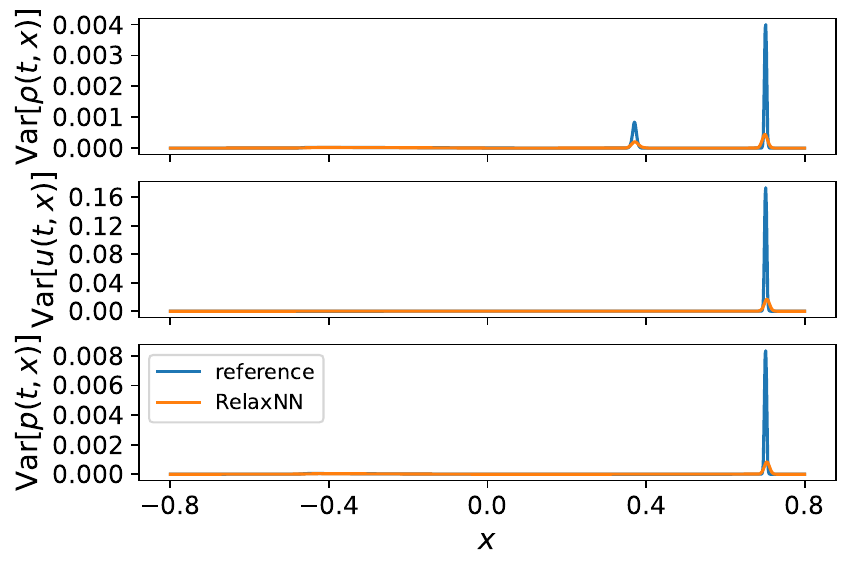}}
\caption{\textit{Euler equations with stochastic Sod shock tube initial condition}: Comparison among the final epoch prediction of RelaxNN(type3) and the reference solution spatially at different specific moments. The resulting relative $L^{2}$ error of mean and variance are $2.26 \times 10^{-2}$, $1.71 \times 10^{-2}$. For RelaxNN, the configuration of $\bm{u}_{\bm{\theta_{1}}}^{\text{NN}}$, $\bm{v}_{\bm{\theta_{2}}}^{\text{NN}}$ are [7,384,384,384,384,384,384,3] and [7,128,128,128,128,128,128,1]. We training for 600,000 epochs and one step per epoch. Loss weights settings are shown in ~\Cref{tab: w settings euler}.}
\label{fig: slice stochastic Euler eq}
\end{figure}

We consider the Euler equations ~\eqref{eq: Euler equations} in domain $ \Omega = \{ (t,x)\} = [ 0.0,0.4]\times [-0.8,0.8]$ with stochastic Sod shock tube initial condition as
\begin{equation}
    \bm{u}_{0}(0,x) = \left( 
    \begin{array}{cc}
         1.0 \\ 
         0.0 \\
         1.0
    \end{array}
    \right) , \quad \text{if} \quad -0.8 \leq x \leq \psi(\mathbf{z}) ;
    \quad
    \bm{u}_{0}(0,x) = \left( 
    \begin{array}{cc}
         0.125 \\
         0.0 \\
         0.1
    \end{array}
    \right) , \quad \text{if} \quad \psi(\mathbf{z}) \leq x \leq 0.8 ;
\end{equation}
Where $\psi(\mathbf{z}) = \varepsilon(z_{1}\sigma(z_{2}z_{3} + z_{4}) +z_{5})$, $\sigma(\cdot)$ is the ReLU function and $\varepsilon = 0.005$.We obtain the reference mean and variance of the solution by the Monte Carlo method for $100,000$ random experiments. The mean and variance of the solution of our RelaxNN are obtained by numerical quadrature. Also, we used the Gaussian-legendre quadrature with $10$ points in every dimension. In ~\Cref{fig: slice stochastic Euler eq}, we display our results.

\bibliographystyle{unsrtnat}
\bibliography{ref.bib}

\end{document}